\documentclass[reqno]{amsart}
\usepackage{graphicx}
\usepackage{pictex}
\usepackage{mathrsfs}
\usepackage[colorlinks]{hyperref}

\begin{document}

\newtheorem{thm}{Theorem}[section]
\newtheorem{lem}[thm]{Lemma}
\newtheorem{cor}[thm]{Corollary}
\newtheorem{add}[thm]{Addendum}
\newtheorem{prop}[thm]{Proposition}

\theoremstyle{definition}
\newtheorem{defn}[thm]{Definition}

\theoremstyle{remark}
\newtheorem{rmk}[thm]{Remark}

\newcommand{\SLtwoC}{\mathrm{SL}(2,{\mathbf C})}
\newcommand{\SLtwoR}{\mathrm{SL}(2,{\mathbf R})}
\newcommand{\PSLtwoC}{\mathrm{PSL}(2,{\mathbf C})}
\newcommand{\PSLtwoR}{\mathrm{PSL}(2,{\mathbf R})}
\newcommand{\SLtwoZ}{\mathrm{SL}(2,{\mathbf Z})}
\newcommand{\PSLtwoZ}{\mathrm{PSL}(2,{\mathbf Z})}
\newcommand{\CmodTwoPiIZ}{{\mathbf C}/2\pi i {\mathbf Z}}
\newcommand{\MCG}{\mathcal{MCG}}
\newcommand{\Cnozero}{{\mathbf C}\backslash \{0\}}
\newcommand{\Cinfty}{{\mathbf C}_{\infty}}
\newcommand{\HH}{{\mathbb H}^2}
\newcommand{\HHH}{{\mathbb H}^3}

\def\square{\hfill${\vcenter{\vbox{\hrule height.4pt \hbox{\vrule width.4pt
height7pt \kern7pt \vrule width.4pt} \hrule height.4pt}}}$}

\newenvironment{pf}{\noindent {\it Proof.}\quad}{\square \vskip 10pt}

\title{Generalized Markoff Maps and McShane's Identity }
\author{Ser Peow Tan, Yan Loi Wong, and Ying Zhang}
\address{Department of Mathematics \\ National University of Singapore \\
2 Science Drive 2 \\ Singapore 117543} \email{mattansp@nus.edu.sg;
matwyl@nus.edu.sg; scip1101@nus.edu.sg}
\address{and the third author}
\address{Department of Mathematics \\Yangzhou University \\Yangzhou 225002 \\P. R. China}
\email{yingzhang@yzu.edu.cn}

\thanks{The authors are partially supported by the National University
of Singapore academic research grant R-146-000-056-112. The third
author is also partially supported by the National Key Basic
Research Fund (China) G1999075104.}

%
%

\begin{abstract}
Following Bowditch \cite{bowditch1998plms}, we study
representations of the free group on two generators into ${\rm
SL}(2,\mathbf C)$, and the connection with generalized Markoff
maps. We show that Bowditch's Q-conditions for generalized Markoff
maps are sufficient for the generalized McShane identity to hold
for the corresponding representations. These conditions are very
close to being necessary as well, and a large class of
representations arising from important and interesting geometric
constructions satisfy these conditions. We also show that the
subset of representations satisfying these conditions is open in
the relative character variety, and it is the largest open subset
on which the mapping class group acts properly discontinuously.
Moreover we generalize Bowditch's results on variations of
McShane's identity for complete, finite volume hyperbolic
3-manifolds which fiber over the circle, with the fiber a
punctured-torus, to identities for incomplete hyperbolic
structures on such manifolds, hence obtaining identities for
closed hyperbolic 3-manifolds which are obtained by doing
hyperbolic Dehn surgery on such manifolds.
\end{abstract}

\maketitle

\vskip 20pt
\section{{\bf Introduction}}\label{s:intro}
\vskip 20pt

In \cite{mcshane1991thesis} McShane proved the following
remarkable identity concerning the lengths of simple closed
geodesics on a once-punctured torus, ${\mathbb T}$, with a
complete, finite area hyperbolic structure:
\begin{eqnarray}\label{eqn:mcshane identity}
\sum_{\gamma}\frac{1}{1+e^{l(\gamma)}}=\frac{1}{2},
\end{eqnarray}
where $\gamma$ ranges over all simple closed geodesics on
${\mathbb T}$, and $l(\gamma)$ is the hyperbolic length of
$\gamma$ under the given hyperbolic structure on ${\mathbb T}$.
This identity is independent of the hyperbolic structure on the
torus, that is, it holds for all points in the Teichm\"uller space
${\mathcal T}$ of the punctured torus. This result was later
generalized to more general hyperbolic surfaces with cusps by
McShane himself \cite{mcshane1998im}, to hyperbolic surfaces with
cusps and/or geodesic boundary components by Mirzakhani
\cite{mirzakhani2004preprint}, to hyperbolic surfaces with cusps,
geodesic boundary and/or conical singularities, as well as to
classical Schottky groups by the authors in
\cite{tan-wong-zhang2004cone-surfaces},
\cite{tan-wong-zhang2004schottky}.

\vskip 10pt

On the other hand, Bowditch in \cite{bowditch1996blms} gave an
alternative proof of (\ref{eqn:mcshane identity}) via Markoff
maps, and extended it in \cite{bowditch1998plms} to
type-preserving representations of the once-punctured torus group
into ${\rm SL}(2, \mathbf C)$ satisfying certain conditions which
we call here the BQ-conditions (Bowditch's Q-conditions). He also
obtained in \cite{bowditch1997t} a variation of (\ref{eqn:mcshane
identity}) which applies to hyperbolic once-punctured torus
bundles. Subsequently, Akiyoshi-Miyachi-Sakuma
\cite{akiyoshi-miyachi-sakuma2004cm355}
\cite{akiyoshi-miyachi-sakuma2004preprint} refined Bowditch's
results in \cite{bowditch1997t} and generalized them to those
which apply to hyperbolic punctured surface bundles.

\vskip 10pt

This paper was first motivated by our attempt to understand the
extent to which the generalized McShane identities for the
one-hole/one-cone torus obtained in
\cite{tan-wong-zhang2004cone-surfaces} hold for general
representations of the torus group into ${\rm SL}(2, \mathbf C)$,
and also the possible generalizations of the McShane-Bowditch
identities in \cite{bowditch1997t} for complete, hyperbolic
punctured torus bundles over the circle to closed hyperbolic
3-manifolds obtained by hyperbolic Dehn surgery on such manifolds.
In the process, we were also led to study some of the fundamental
properties of ``marked'' two generator subgroups of $\SLtwoC$ via
the (generalized) Markoff maps studied by Bowditch in
\cite{bowditch1998plms}, and to determine which of these
properties had generalizations and extensions to the generalized
Markoff maps.

\vskip 10pt

A proper statement of the results requires a fair bit of notation,
and will be deferred to the next section.  We give here first a
brief outline. The first main result (Theorem \ref{thm:TWZ}) is
that the generalized McShane identity obtained in
\cite{tan-wong-zhang2004cone-surfaces} for the one-hole/one-cone
torus holds for representations of $\Gamma$ into ${\rm SL}(2,
\mathbf C)$ which satisfy the BQ-conditions, which we will
describe in detail in the next section. Roughly speaking, these
conditions say that there are no elliptics arising from the simple
closed curves, and that the number of simple closed curves with
length below a certain fixed bound is finite. There are many
interesting examples of representations which satisfy these
conditions arising from some basic and important geometric
constructions. For example, representations arising from the
one-cone/one-hole hyperbolic torus, the hyperbolic three-holed
sphere (often called pair of pants), and more generally, classical
Schottky groups with two generators all give rise to such
representations. However, the class of representations satisfying
the BQ-conditions is in general much larger, and an interesting
problem is to classify this class geometrically. In the case of
type-preserving representations studied in
\cite{bowditch1998plms}, Bowditch has conjectured that the class
coincides with the quasifuchsian representations (Conjecture A in
\cite{bowditch1998plms}). The general case is less clear since the
representations are in general not discrete, or not faithful. The
basic strategy for proving Theorem \ref{thm:TWZ} follows very
closely that used in \cite{bowditch1998plms}. The problem is
reformulated in terms of generalized Markoff maps and it is shown
that if the BQ-conditions are satisfied, then the generalized
Markoff map has Fibonacci growth. This is sufficient to obtain the
absolute convergence of the series on the left hand side of the
identity. To obtain the actual value, in the type-preserving case,
Bowditch made ingenious use of the quantities $\frac{x}{yz}$,
$\frac{y}{zx}$ and $\frac{z}{xy}$ associated to a Markoff triple
$(x,y,z)$ and certain properties of these quantities.
By exploring the geometric interpretation of these quantities, we
were able to find analogous quantities, $\Psi(y,z,x), \Psi(z,x,y)$
and $\Psi(x,y,z)$(for explicit expressions see \S \ref{s:proof of
theorem mu}), associated to a generalized Markoff triple, and
follow through the proof of Bowditch to complete the proof of
Theorem \ref{thm:TWZ}.

\vskip 5pt

The mapping class group of the torus acts naturally on the space
of equivalence classes of representations, preserving the subspace
consisting of representations with fixed trace of the commutator
$[a,b]$ of a generating pair $a,b$ for $\Gamma$, also called the
relative character variety. The case where the trace of the
commutator is $-2$ was the case studied by Bowditch in
\cite{bowditch1998plms}, see also the work of Minsky in
\cite{minsky1999am} where the problem was studied within the
context of Kleinian groups and the ending lamination conjecture.
For real representations, Goldman has studied the dynamics of this
action in \cite{goldman2003gt} and obtained quite complete
results. Our next result is that for each subspace, the mapping
class group acts properly discontinuously on the subset of
representations satisfying the BQ-conditions (Theorem
\ref{thm:proper}). Furthermore, this is  the largest open subset
for which this is true (Proposition
\ref{prop:accummulationpoint}).

\vskip 5pt

 Bowditch also gave variations of McShane's
identity for complete, finite volume hyperbolic 3-manifolds which
fiber over the circle, with fiber the once-punctured torus
\cite{bowditch1997t}. These complete structures can be deformed to
incomplete structures, as shown by Thurston in
\cite{thurston1978notes}, and in certain cases, one can perform
hyperbolic Dehn surgery to obtain closed (complete) hyperbolic
3-manifolds. The next set of results (Theorems \ref{thm:A'},
\ref{thm:B'}) is that a further variation of the McShane-Bowditch
identity holds for these deformations satisfying a variation of
the BQ-conditions, and hence for the closed hyperbolic 3-manifolds
obtained by hyperbolic Dehn surgery on such manifolds.
An interesting aspect of this result
is that it is easier to obtain generalizations of McShane's
identity for this class of closed (complete) hyperbolic
3-manifolds than it is for closed hyperbolic surfaces. The only
example in the surface case which we know of is an identity for
the genus two surface 
(see \cite{mcshane1998preprint} and
\cite{tan-wong-zhang2004cone-surfaces}).

All the results above are proven by transforming to equivalent
statements for generalized Markoff maps and working in that
framework. We remark that besides these connections with
hyperbolic geometry, generalized Markoff maps are also closely
related to dynamical systems (see \cite{goldman2003gt},
\cite{brooks-matelski1981}), algebraic number theory (see
\cite{cusick-flahive1989msm30}, and \cite{silverman1990jnt}) and
mathematical physics (see \cite{roberts1996physa}). There also
seem to be a close connection
 to the pair-of-pants complex introduced by  Hatcher and Thurston in
\cite{hatcher-thurston1980t}, and developed in connection with the
Weil-Petersson metric on  moduli space in the recent work by Brock
in \cite{brock2003jams} (see also \cite{mirzakhani2004preprint}
for another connection of McShane's identity with the
Weil-Petersson volume). A further exploration of these ideas may
lead to the correct framework of Markoff maps for general
surfaces.

\vskip 10pt

The rest of this paper is organized as follows. In \S
\ref{s:notation+results} we give the definitions and state the
results in terms of representations. In \S \ref{s:generalized
Markoff maps} we reformulate the results in terms of generalized
Markoff maps and also state and prove some basic and fundamental
results for generalized Markoff maps, including Theorem
\ref{thm:proper} and Proposition \ref{prop:accummulationpoint}. In
\S \ref{s:proof of theorem mu} we introduce the $\Psi$ function
and give the proof of Theorem \ref{thm:mcshane mu-markoff}, the
reformulation of our first main theorem. In \S \ref{s:torus
bundles} we generalize Bowditch's variations of McShane's identity
to identities for certain incomplete hyperbolic torus bundles over
the circle (Theorems \ref{thm:A'}, \ref{thm:B'}), with
applications to hyperbolic Dehn surgery. Finally, in the
Appendices, we give the geometric meaning of two of the important
functions $h$ and $\Psi$ used in this paper (Appendix A), and
  explain how to draw the gaps in the extended complex plane to
visualize the generalized McShane's identity (\ref{eqn:mcshane
mu-markoff}) (Appendix B).

\vskip 5pt

\noindent {\it Acknowledgements.} We would like to thank Bill
Goldman, Caroline Series and Makoto Sakuma for their
encouragement, helpful conversations/correspondence and comments.

 \vskip 20pt
\section{{\bf Notation and statements of results}}\label{s:notation+results}%
\vskip 10pt

In this section, we set some basic notation and definitions and
give precise statements of our results in terms of representations
of the once-punctured torus group, $\Gamma$, into ${\rm SL}(2,
\mathbf C)$.
As much of the paper is influenced by \cite{bowditch1998plms} and
\cite{bowditch1997t}, we will borrow the notation and definitions
from them as much as possible to avoid confusion.

\vskip 10pt

\noindent {\bf The punctured torus group.}\,\, Let ${\mathbb T}$
be a (topological) once-punctured punctured torus and let $\Gamma
$ be its fundamental group. Note that algebraically, $\Gamma \cong
{\mathbf Z} \ast {\mathbf Z}=\langle a,b \rangle$, the free group
on two generators.

We define an equivalence relation, $\sim$, on $\Gamma$ such that
$g \sim h$ if and only if $g$ is conjugate to $h$ or $h^{-1}$.
Note that $\Gamma/\sim$ can be identified with the set of free
homotopy classes of unoriented closed curves on ${\mathbb T}$.

\vskip 10pt

\noindent {\bf Simple closed curves on the torus.}\,\, Let
${\mathscr C}$ be the set of free homotopy classes of non-trivial,
non-peripheral simple closed curves on ${\mathbb T}$ and let $\hat
\Omega \subset \Gamma/\sim$ be the subset corresponding to
${\mathscr C}$.

\vskip 5pt



Note that $\hat \Omega$ can be identified with ${\mathbf Q}\cup
\{\infty\}$ by considering the ``slope'' of $[g] \in \hat \Omega$
as follows. Fix a pair of generators $a$ and $b$ of $\Gamma$. Then
each class $[g]$ of $\hat \Omega$ has a representative $g=W(a^{\pm
1}, b)$ which is a cyclically reduced word in $\{a^{\pm 1},b\}$,
and the exponents of $a$ in each word is either all positive or
all negative. The word is unique up to cyclic permutations, and
the slope is the quotient of the sum of the exponents of $a$ with
the sum of the exponents of $b$ in the word; see
\cite{series1985mi} for details. Hence, $a$ is identified with
$\infty$, $b$ with $0$, $ab$ with $1$, and so on. Note that $\hat
\Omega$ inherits a cyclic ordering from the cyclic ordering of
${\mathbf Q}\cup \{\infty\}$ induced from the standard embedding
into ${\mathbf R}\cup \{\infty\} \cong S^1$.

\vskip 10pt

%

\noindent {\bf The $\tau$-representations.}\,\, A representation
$\rho: \Gamma \rightarrow {\rm SL}(2, \mathbf C)$ is said to be a
$\tau$-{\it representation}, where $\tau \in \mathbf C$, if for
some (hence every) pair of free generators $a, b \in \Gamma$,
${\rm tr}\rho([a,b]) = \tau$, where $[a,b]=aba^{-1}b^{-1}$. The
space of $\tau$-representations, modulo conjugation, is denoted by
${\mathcal X}_{\tau}$. In particular, ${\mathcal X}_{-2}$ is the
space of type-preserving representations. The mapping class group
of ${\mathbb T}$,  $\MCG$ acts on ${\mathcal X}_{\tau}$.

\vskip 5pt

Note that we can also talk about $\tau$-representations for
representations $\rho$  into ${\rm PSL}(2, \mathbf C)$, since the
trace of the commutator is well-defined. Most geometric
constructions give rise to representations into ${\rm PSL}(2,
\mathbf C)$. It will, however, be more convenient for us to work
with representations into ${\rm SL}(2, \mathbf C)$, as the
statements of the results are neater in this case. This will not
affect the validity of our results since the results stated will
be independent of the lift chosen, and the results could have been
stated in terms of representations into ${\rm PSL}(2, \mathbf C)$
as well.

\vskip 5pt

It is well known that the image of the representation $\rho$ is
non-elementary ($\rho$ is irreducible) if and only if $\tau \neq
2$, we will only be interested in $\tau$-representations with
$\tau \neq 2$ in this paper. The case where $\tau=-2$ is
particularly interesting. In this case, the representation $\rho$
is said to be {\it type-preserving}, and except when ${\rm
tr}\rho(a)={\rm tr}\rho(b)={\rm tr}\rho(ab)=0$, $\rho([a,b])$ is
always a parabolic element for any pair of generating elements
$a,b$ in $\Gamma$. This is the case usually considered in the
study of Kleinian groups, (see \cite{minsky1999am}), and is also
the case extensively studied in \cite{bowditch1998plms}.

\vskip 10pt

\noindent {\bf The BQ-conditions.}\,\, For a fixed $\tau \neq 2$,
a $\tau$-representation $\rho: \Gamma \rightarrow {\rm SL}(2,
\mathbf C)$ (or ${\rm PSL}(2, \mathbf C)$) is said to satisfy the
{\it BQ-conditions} if

\vskip 2pt

(BQ1) ${\rm tr}\rho(g) \not\in [-2,2]$ for all $[g] \in \hat
\Omega$; and

(BQ2) $|{\rm tr}\rho(g)| \le 2$ for only finitely many (possibly
none) $[g]\in \hat \Omega$.

\vskip 2pt

\noindent We also call such a representation $\rho$ a {\it
BQ-representation}, or Bowditch representation, and the space of
such representations the Bowditch representation space, denoted by
$({\mathcal X}_{\tau})_Q$. If we replace (BQ1) by ${\rm tr}\rho(g)
\not\in (-2,2)$ for all $[g] \in \hat \Omega$, we call the
resulting space the extended Bowditch representation space (see
\cite{tan-wong-zhang2004necsuf}).

\vskip 5pt

Note that ${\rm tr}\rho(g_1)={\rm tr}\rho(g_2)$ if $[g_1]=[g_2]$
(since $g_1$ is conjugate to $g_2$ or its inverse by definition);
so the conditions (BQ1) and (BQ2) make sense.

\vskip 10pt

\noindent {\bf Some conventions.} \,\, Throughout this paper we
assume that
\begin{itemize}
\item
the function ${\bf cosh}^{-1}$ has images with nonnegative real
parts and with imaginary parts in $(-\pi, \pi]$;

\item
the function ${\bf log}$ has images with imaginary parts in
$(-\pi, \pi]$; while

\item
the function ${\bf tanh}^{-1}$ has images with imaginary parts in
$(-\pi/2, \pi/2]$.
\end{itemize}

\vskip 10pt

\noindent {\bf The functions $l/2$ and $l$.}\,\, For $x \in
{\mathbf C}$, let $l(x)/2 \in {\mathbf C}/2 \pi i {\mathbf Z}$ be
defined by
\begin{eqnarray}\label{eqn:l(x)/2=}
l(x)/2=\cosh^{-1}(x/2).
\end{eqnarray}
Hence $l(x)=2\cosh^{-1}(x/2)=\cosh^{-1}(x^2/2-1) \in {\mathbf C}/2
\pi i {\mathbf Z}$. In particular, $\Re \big(l(x)/2\big) \ge 0$
and if $\Re \big(l(x)/2\big) = 0$ then $\Im \big(l(x)/2\big) \ge
0$.

\vskip 5pt

Note that $e^{l(x)/2} \in {\mathbf C}$ is well-defined. It can be
shown (see Lemma \ref{lem:e^-l(x)/2}) that if $x \notin [-2,2]$
then $e^{-l(x)/2}=xh(x)$ where $h(x)$ is defined below. Hence
\begin{eqnarray}\label{eqn:e^l(x)=}
e^{l(x)}=x^{-2}h(x)^{-2}=h(x)^{-1}-1.
\end{eqnarray}

\vskip 5pt

For $A \in \SLtwoC$, we define its half translation length $l(A)/2
\in \CmodTwoPiIZ$ by
\begin{eqnarray}\label{eqn:l(A)/2=}
\frac{l(A)}{2}=\frac{l({\rm tr}A)}{2}=\cosh^{-1}\Big(\frac{{\rm
tr}A}{2}\Big).
\end{eqnarray}

\vskip 10pt

\noindent {\bf The function $h$.}\,\, We define an even function
$h: {\mathbf C} \backslash \{0\} \rightarrow {\mathbf C}$ by
\begin{eqnarray}\label{eqn:h(x)=}
h(x)=\frac{1}{2}\bigg(1-\sqrt{1-\frac{4}{x^2}}\,\bigg).
\end{eqnarray}
It is easy to check that
\begin{eqnarray}\label{eqn:h(x)^2-h(x)}
h(x)^2-h(x)+x^{-2}=0.
\end{eqnarray}
In fact, if $x \notin [-2,2]$ then $h(x)$ is the root of the
quadratic equation which has smaller real part. By
(\ref{eqn:e^l(x)=}), we also have
\begin{eqnarray}
h(x)=\frac{1}{1+e^{l(x)}}.
\end{eqnarray}

\vskip 10pt

\noindent {\bf The function ${\mathfrak h}={\mathfrak
h}_\tau$.}\,\, For $\tau \in {\mathbf C}$, set
$\nu=\cosh^{-1}(-\tau/2)$. We define a function

\centerline{${\mathfrak h}={\mathfrak h}_\tau: {\mathbf C}
\backslash \{\pm \sqrt{\tau+2}\} \rightarrow {\mathbf C}$}

\noindent by
\begin{eqnarray}
{\mathfrak h}(x)
&=&2 \tanh^{-1}\bigg(\frac{\sinh\nu}{\cosh\nu+e^{l(x)}}\bigg)
\label{eqn:frak h(x)=2tanh^-1} \\
&=&\log\frac{e^{\nu}+e^{l(x)}}{e^{-\nu}+e^{l(x)}} \label{eqn:frak h(x)=log I} \\
&=&\log\frac{1+(e^{\nu}-1)\,h(x)}{1+(e^{-\nu}-1)\,h(x)},
\label{eqn:frak h(x)=log II}
\end{eqnarray}
where (\ref{eqn:frak h(x)=log II}) follows from (\ref{eqn:frak
h(x)=log I}) by (\ref{eqn:e^l(x)=}). Note that ${\mathfrak h}(0)$
is well-defined if $\tau \neq -2$. In fact,
\begin{eqnarray*}
{\mathfrak h}(0)=\nu + \pi i.
\end{eqnarray*}
Note also that $e^{\pm \nu} + e^{l(x)} = 0$ if and only if
\begin{eqnarray*}
\pm \nu + \pi i = l(x),
\end{eqnarray*}
which implies that $x^2=\tau +2$.

\vskip 10pt

\noindent {\bf Geometric interpretation of $h$ and ${\mathfrak h}$
in the real case.}\,\, When the representation arises from a real
one-cone/cusp/holed torus, then $h$ and ${\mathfrak h}$ correspond
to the gaps on the boundary of the torus arising from intersecting
geodesics, see \cite{mcshane1998im} or
\cite{tan-wong-zhang2004cone-surfaces}, we give a more detailed
description in Appendix A.

\vskip 15pt

\noindent {\bf Bowditch's extension of McShane's identity.}\,\,
With the above notation and definitions, Bowditch's extension and
reformulation of McShane's identity (\ref{eqn:mcshane identity})
can be stated as follows.

\vskip 5pt

\begin{thm}\label{BM} {\rm (Theorem 3 in \cite{bowditch1998plms})}
Let $\rho: \Gamma \rightarrow {\rm SL}(2, \mathbf C)$ be a
type-preserving representation which satisfies the BQ-conditions.
Then
\begin{eqnarray}\label{eqn:mcshane}
\sum_{[g] \in \hat \Omega}h\big({\rm tr}\rho(g)\big)=\frac{1}{2},
\end{eqnarray}
where the sum converges absolutely.
\end{thm}

\vskip 10pt

\noindent{\it Remarks.}\,\,
\begin{itemize}

\item[(i)] The identity (\ref{eqn:mcshane}) is also true for
BQ-representations into ${\rm PSL}(2, \mathbf C)$ since $h(x)$ is
an even function and the trace of an element of ${\rm PSL}(2,
\mathbf C)$ is well-defined up to sign.

\item[(ii)] For an element $A \in {\rm PSL}(2, \mathbf C)$, its {\it complex
translation length} $$l(A) =l({\rm tr}A)=\cosh^{-1}
(\textstyle\frac12\,{\rm tr}^2A-1) \in {\mathbf C} / 2\pi i
{\mathbf Z}$$ is well-defined with $\Re\, l(A) \ge 0$ and $h({\rm
tr}A)={1}/(1+e^{l(A)})$.
\end{itemize}

\vskip 10pt

\noindent {\bf Our first main theorem}, which is a generalization
of  Theorem \ref{BM} of Bowditch above, and Theorem 1.4 of
\cite{tan-wong-zhang2004cone-surfaces}, is as follows:

\vskip 10pt

\begin{thm}\label{thm:TWZ}
Let $\rho: \Gamma \rightarrow {\rm SL}(2, \mathbf C)$ be a
$\tau$-representation {\rm(}where $\tau \neq 2${\rm)} satisfying
the BQ-conditions. Set $\nu=\cosh^{-1}(-\tau/2)$. Then
\begin{eqnarray}\label{eqn:TWZ}
\sum_{[g] \in \hat \Omega}{\mathfrak h} \big( {\rm
tr}\rho(g)\big)=\nu  \mod 2 \pi i,
\end{eqnarray}
where the sum converges absolutely.
\end{thm}

\vskip 8pt

\noindent {\it Remarks.} \begin{itemize}

\item[(i)] Here $\nu$ is a specific
choice of {\it half} the complex translation length of the
commutator $[\rho(a),\rho(b)]$.

\item[(ii)] Since ${\mathfrak h}$ is an even function, Theorem \ref{thm:TWZ} also holds for
representations into ${\rm PSL}(2, \mathbf C)$ satisfying the
BQ-conditions.

\item[(iii)] Note that for $x \neq 0$, when $\tau
\rightarrow -2$, or equivalently, $\nu \rightarrow 0$, we have
\begin{eqnarray*}
{\mathfrak h}(x)&=&\log
\frac{1+(e^{\nu}-1)h(x)}{1+(e^{-\nu}-1)h(x)} \\
&=&\log
\bigg(1+\frac{2\sinh\nu \,h(x)}{1+(e^{-\nu}-1)h(x)}\bigg) \\
&\sim& 2 \nu h(x).
\end{eqnarray*}
Hence the identity (\ref{eqn:mcshane}) can be `obtained' by
considering the first order infinitesimal terms of
(\ref{eqn:TWZ}).

\item[(iv)] Theorem \ref{thm:TWZ} has an equivalent formulation as Theorem
\ref{thm:mcshane mu-markoff} in terms of $\mu$-Markoff maps.
\end{itemize}

%

\vskip 10pt

Next, consider the action of the mapping class group of $\mathbb
T$, $\mathcal{MCG} \cong {\rm SL}(2,{\mathbf Z})$, and its induced
action  on $\Gamma$ and on ${\mathcal X}_{\tau}$, the space of
$\tau$-representations. (Recall that $\mathcal{MCG}$ is the group
of isotopy classes of diffeomorphisms of ${\mathbb T}$ which fixes
the puncture, that is, there exists a neighborhood $N$ of the
puncture which is pointwise fixed by the elements of
$\mathcal{MCG}$.) Any $H \in \mathcal{MCG}$ induces an
automorphism $H_{\ast}$ of $\Gamma$, and $H$ acts on ${\mathcal
X}_{\tau}$ by
$$H(\rho)(g)=\rho(H_{\ast}(g)),$$
for any $\rho \in {\mathcal X}_{\tau}$ and $g \in \Gamma$. We
shall see in the next section that the set $({\mathcal
X}_{\tau})_Q$ of Bowditch representations is open in ${\mathcal
X}_{\tau}$ (Theorem \ref{thm:B3.16}). It is also clear that
$\mathcal{MCG}$ acts on $({\mathcal X}_{\tau})_Q$. We have the
following result:

\vskip 8pt

\begin{thm}\label{thm:proper}
$\mathcal{MCG}$ acts properly discontinuously on $({\mathcal X}_{\tau})_Q$. %
\end{thm}

We also have a partial converse as follows:

\vskip 5pt

\begin{prop}\label{prop:accummulationpoint}
Suppose that $\rho \in {\mathcal X}_{\tau}$ does not satisfy the
extended BQ-conditions, that is, either ${\rm tr}\rho(g) \in
(-2,2)$ for some $[g ]\in \hat\Omega$ or $|{\rm tr}\rho(g)| \le 2$
for infinitely many $[g] \in \hat\Omega$. Then there exists a
sequence of distinct elements $H_i \in \mathcal{MCG}$ and $\rho_0
\in {\mathcal X}_{\tau}$ such that $H_i(\rho)$ converges to
$\rho_0 \in {\mathcal X}_{\tau}$.
\end{prop}

From Proposition \ref{prop:accummulationpoint} it is easy to
deduce that $({\mathcal X}_{\tau})_Q$ is the largest open subset
of ${\mathcal X}_{\tau}$ on which $\MCG$ acts properly
discontinuously. The proofs of Theorem \ref{thm:proper} and
Proposition \ref{prop:accummulationpoint} will be deferred to the
end of \S \ref{s:generalized Markoff maps}.

\vskip 5pt

Bowditch \cite{bowditch1997t} studied representations $\rho \in
{\mathcal X}_{tp}={\mathcal X}_{-2}$ stabilized by a cyclic
subgroup $\langle H \rangle < \mathcal{MCG}$ generated by a
hyperbolic element and proved a variation of the McShane's
identity. This result can be generalized for
$\tau$-representations as follows. The case $\tau=-2$ with
$\mathfrak h$ replaced by $h$ is Bowditch's variation (Theorem A
in \cite{bowditch1997t}).

\begin{thm}\label{thm:BowditchT}
Suppose that a $\tau$-representation $\rho: \Gamma \rightarrow
\SLtwoC$, where $\tau \neq 2$, is stabilized by a hyperbolic
element $H \in {\rm SL}(2,{\mathbf Z}) \cong \mathcal{MCG}$ and
$\rho$ satisfies the BQ-conditions on $\hat \Omega/\langle
H_{\ast}\rangle$, that is,
\begin{itemize}

\item[(i)] ${\rm tr}\rho([g]) \notin [-2,2]$ for all classes $[g] \in \hat
\Omega/\langle H_{\ast}\rangle$, and

\item[(ii)] $|{\rm tr}\rho([g])| \le 2$ for only finitely many {\rm(}possibly
no{\rm)} classes $[g] \in \hat \Omega/\langle H_{\ast}\rangle$.
\end{itemize} Then
\begin{eqnarray}\label{eqn:hyperbolic stab}
\sum_{[g] \in \hat \Omega/\langle H_{\ast}\rangle} {\mathfrak
h}\big( {\rm tr}\rho([g])\big)=0 \mod 2\pi i,
\end{eqnarray}
where the sum converges absolutely.
\end{thm}

\vskip 5pt

Theorem \ref{thm:BowditchT} will be reformulated and proved as
Theorem \ref{thm:stabilized} in \S \ref{s:torus bundles}, together
with a more interesting version where the sum is taken only over
the left partition (Theorem \ref{thm:longitude}), with
interpretations in terms of (in)complete hyperbolic structures on
torus bundles (Theorems \ref{thm:A'} and \ref{thm:B'}).


\vskip 20pt
\section{{\bf Generalized Markoff maps}}\label{s:generalized Markoff maps} %
\vskip 10pt

In this section we study the basic geometry of generalized Markoff
maps, the connection with the space of equivalence classes of
representations ${\mathcal X}$, and show that many of the
fundamental results obtained by Bowditch for Markoff maps extend
to  this case as well. In particular,  we obtain the Fibonacci
growth for generalized Markoff maps which satisfy Bowditch's
Q-conditions. We shall follow the notation and proofs of Bowditch
\cite{bowditch1998plms} whenever possible. However, instead of
working in the more general combinatorial set-up there, we will
fix a concrete realization of the combinatorial structure involved
to help with the visualization of the maps.


\begin{figure}
\setlength{\unitlength}{1mm} 
\begin{picture}(60,30)
\thicklines \put(20,10){\vector(1,0){20}}
\put(20,10){\line(-1,1){10}} \put(20,10){\line(-1,-1){10}}
\put(40,10){\line(1,1){10}} \put(40,10){\line(1,-1){10}}
\put(30,20){$X$} \put(30,0){$Y$} \put(50,10){$W$} \put(6,10){$Z$}
\put(30,11){$\vec e$}
\end{picture}
\caption{The directed edge $\vec e = (X,Y;Z \rightarrow W)$}
\label{fig:edge oriented}
\end{figure}
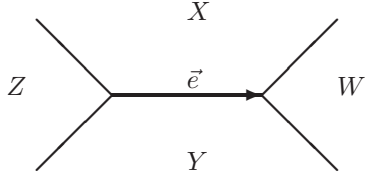

\begin{figure}
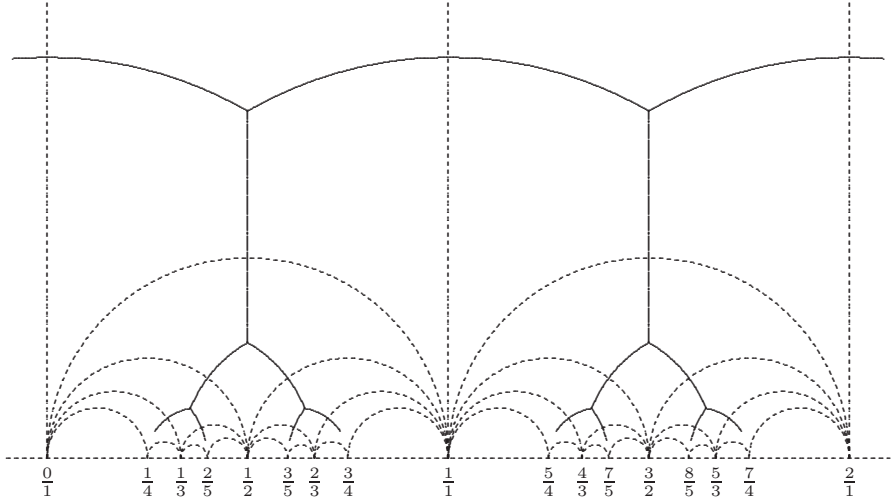


\begin{center}
\mbox{
\beginpicture
\setcoordinatesystem units <2.1in,2.1in>

\setplotarea x from -0.1 to 2.1, y from 0 to 1.2


\plot 0.5 0.86603 0.5 0.28868 /

\circulararc 38.21321 degrees from 0.5 0.28868  center at 0.66667
0

\circulararc -38.21321 degrees from 0.5 0.28868 center at 0.33333
0

\circulararc 49.58256 degrees from 0.35714 0.12372  center at
0.375 0

\circulararc -25.03966 degrees from 0.35714 0.12372 center at 0.2
0

\circulararc 25.03966 degrees from 0.64286 0.12372 center at 0.8 0

\circulararc -49.58256 degrees from 0.64286 0.12372  center at
0.625 0


\plot 1.5 0.86603 1.5 0.28868 /

\circulararc 38.21321 degrees from 1.5 0.28868  center at 1.66667
0

\circulararc -38.21321 degrees from 1.5 0.28868 center at 1.33333
0

\circulararc 49.58256 degrees from 1.35714 0.12372  center at
1.375 0

\circulararc -25.03966 degrees from 1.35714 0.12372 center at 1.2
0

\circulararc 25.03966 degrees from 1.64286 0.12372 center at 1.8 0

\circulararc -49.58256 degrees from 1.64286 0.12372  center at
1.625 0


\circulararc 60 degrees from 1.5 0.86603  center at 1 0

\circulararc 35 degrees from 0.5 0.86603  center at 0 0

\circulararc -35 degrees from 1.5 0.86603  center at 2 0


\put {\mbox{\small $\frac{0}{1}$}} [cb] <0mm,-5mm> at 0 0

\put {\mbox{\small $\frac{1}{1}$}} [cb] <0mm,-5mm> at 1 0

\put {\mbox{\small $\frac{1}{2}$}} [cb] <0mm,-5mm> at 0.5 0

\put {\mbox{\small $\frac{1}{3}$}} [cb] <0mm,-5mm> at 0.33333 0

\put {\mbox{\small $\frac{2}{3}$}} [cb] <0mm,-5mm> at 0.66667 0

\put {\mbox{\small $\frac{1}{4}$}} [cb] <0mm,-5mm> at 0.25 0

\put {\mbox{\small $\frac{2}{5}$}} [cb] <0mm,-5mm> at 0.4 0

\put {\mbox{\small $\frac{3}{5}$}} [cb] <0mm,-5mm> at 0.6 0

\put {\mbox{\small $\frac{3}{4}$}} [cb] <0mm,-5mm> at 0.75 0


\put {\mbox{\small $\frac{2}{1}$}} [cb] <0mm,-5mm> at 2 0

\put {\mbox{\small $\frac{3}{2}$}} [cb] <0mm,-5mm> at 1.5 0

\put {\mbox{\small $\frac{4}{3}$}} [cb] <0mm,-5mm> at 1.33333 0

\put {\mbox{\small $\frac{5}{3}$}} [cb] <0mm,-5mm> at 1.66667 0

\put {\mbox{\small $\frac{5}{4}$}} [cb] <0mm,-5mm> at 1.25 0

\put {\mbox{\small $\frac{7}{5}$}} [cb] <0mm,-5mm> at 1.4 0

\put {\mbox{\small $\frac{8}{5}$}} [cb] <0mm,-5mm> at 1.6 0

\put {\mbox{\small $\frac{7}{4}$}} [cb] <0mm,-5mm> at 1.75 0


\setdashes<1.5pt>

\plot -0.1 0 2.1 0 /

\plot 0 0 0 1.15 / \plot 1 0 1 1.15 /

\circulararc 180 degrees from 1 0  center at 0.5 0

\circulararc 180 degrees from 0.5 0  center at 0.25 0

\circulararc 180 degrees from 1 0  center at 0.75 0

\circulararc 180 degrees from 0.33333 0  center at 0.16667 0

\circulararc 180 degrees from 0.5 0  center at 0.41667 0

\circulararc 180 degrees from 0.66667 0  center at 0.58333 0

\circulararc 180 degrees from 1 0  center at 0.83333 0

\circulararc 180 degrees from 0.25 0  center at 0.125 0

\circulararc 180 degrees from 0.33333 0  center at 0.29167 0

\circulararc 180 degrees from 0.4 0  center at 0.36667 0

\circulararc 180 degrees from 0.5 0  center at 0.45 0

\circulararc 180 degrees from 0.6 0  center at 0.55 0

\circulararc 180 degrees from 0.66667 0  center at 0.63333 0

\circulararc 180 degrees from 0.75 0  center at 0.70833 0

\circulararc 180 degrees from 1 0  center at 0.875 0

\plot 2 0 2 1.15 /

\circulararc 180 degrees from 2 0  center at 1.5 0

\circulararc 180 degrees from 1.5 0  center at 1.25 0

\circulararc 180 degrees from 2 0  center at 1.75 0

\circulararc 180 degrees from 1.33333 0  center at 1.16667 0

\circulararc 180 degrees from 1.5 0  center at 1.41667 0

\circulararc 180 degrees from 1.66667 0  center at 1.58333 0

\circulararc 180 degrees from 2 0  center at 1.83333 0

\circulararc 180 degrees from 1.25 0  center at 1.125 0

\circulararc 180 degrees from 1.33333 0  center at 1.29167 0

\circulararc 180 degrees from 1.4 0  center at 1.36667 0

\circulararc 180 degrees from 1.5 0  center at 1.45 0

\circulararc 180 degrees from 1.6 0  center at 1.55 0

\circulararc 180 degrees from 1.66667 0  center at 1.63333 0

\circulararc 180 degrees from 1.75 0  center at 1.70833 0

\circulararc 180 degrees from 2 0  center at 1.875 0


\endpicture
}\end{center}

\caption{Farey tessellation and the binary tree
$\Sigma$}\label{fig:Farey}
\end{figure}

\vskip 5pt \noindent {\bf The Farey triangulation ${\mathcal F}$
and the dual binary tree $\Sigma$.}\,\, Let ${\mathcal F}$ be the
Farey triangulation of the hyperbolic plane $\HH$, and $\Sigma$ be
the dual binary tree (also called an infinite trivalent graph).
Recall that the Farey triangulation consists of edges which are
complete hyperbolic geodesics joining all pairs $\{\xi,\eta\}
\subset {\mathbf Q} \cup \{\infty\}$ which are Farey neighbors,
and that $\xi=p/q$ and $\eta=r/s$ (where we always assume that
$p,q,r,s \in {\mathbf Z}$ and $(p,q)=(r,s)=1$) are Farey neighbors
if $ps-rq=\pm 1$. See Figure \ref{fig:Farey}.


\vskip 5pt \noindent {\bf Complementary regions.}\,\,A {\it
complementary region} of $\Sigma$ is the closure of a connected
component of the complement.

We denote by $\Omega=\Omega(\Sigma)$ the set of complementary
regions of $\Sigma$. Similarly, we use $V(\Sigma)$, $E(\Sigma)$
for the set of vertices and edges, respectively.


\vskip 5pt

\noindent {\bf Action of $\PSLtwoZ$.}\,\,It is clear from the
construction that there is a natural action of ${\rm PSL}(2,
\mathbf Z)$ (in fact, ${\rm PGL}(2, \mathbf Z)$) on $\Sigma$, and
there is a natural correspondence of $\Omega(\Sigma)$ with
${\mathbf Q}\cup\{\infty\}$ (the vertices of ${\mathcal F}$),
together with the induced cyclic ordering.

\vskip 5pt

We use the letters $X,Y,Z,W, \ldots$ to denote the elements of
$\Omega$, and also introduce the notation $X(p/q)$ to indicate
that $X \in \Omega$ corresponds to $p/q \in {\mathbf Q} \cup
\{\infty\}$. We also use the notation $e \leftrightarrow
(X,Y;Z,W)$ to indicate that $e=X \cap Y$ and $e\cap Z$ and $e \cap
W$ are the endpoints of $e$; see Figure \ref{fig:edge oriented},
ignoring the direction of $e$ there. Note also that there is a
correspondence of the ends of $\Sigma$ with ${\mathbf R} \cup
\{\infty\} \cong S^1$ which is one-to-one from ends to the
irrationals and two-to-one from ends to the rational numbers, much
like the decimal expansions for the real numbers although we shall
not use it in this paper.

\vskip 5pt

\noindent {\bf Directed edges, the sets $\Omega^0(\vec e)$,
$\Omega^\pm(\vec e)$ and $\Omega^{0\pm}(\vec e).$}\,\, Denote by
$\vec E(\Sigma)$ (or just $\vec E$) the set of directed edges of
$\Sigma$ where the direction is always taken to be from from the
tail to the head (as in the direction of the arrow), and for a
directed edge $\vec e \in \vec E(\Sigma)$, we use $\vec e = (X,Y;Z
\rightarrow W)$ (or just $\vec e = (X,Y; \rightarrow W)$), to
indicate that $e \cap W$ is the head of $\vec e$, that is, $\vec
e$ is the directed edge from $Z$ to $W$, see Figure \ref{fig:edge
oriented}. We also use $e$ to indicate the underlying undirected
edge corresponding to the directed edge $\vec e$, and $-\vec e$ to
indicate the directed edge in the opposite direction of $\vec e$.
Associated to each $\vec e \in \vec E(\Sigma)$ is a partition of
$\Omega$ which we will use repeatedly later. For $e
\leftrightarrow (X,Y;Z,W)\in E$, define $\Omega^0(e)=\{X,Y\}$. For
$\vec e = (X,Y;Z \rightarrow W)\in \vec E$, $\Sigma \setminus
\mbox{int}(e)$ consists of two components, denoted by
$\Sigma^+(\vec e)$ and $\Sigma^-(\vec e)$, where $\Sigma^+(\vec
e)$ is the component containing the head of $\vec e$ and
$\Sigma^-(\vec e)$ is the component containing the tail of $\vec
e$ (so that $Z$ meets $\Sigma^-(\vec e)$ and $W$ meets
$\Sigma^+(\vec e)$). Define $\Omega^{+}(\vec e)$ (resp.
$\Omega^{-}(\vec e)$) to be the set of regions in $\Omega$ whose
boundaries lie in $\Sigma^{+}(\vec e)$ (resp. $\Sigma^{-}(\vec
e)$). Hence $\Omega = \Omega^{+}(\vec e) \sqcup \Omega^0(e) \sqcup
\Omega^{-}(\vec e)$. Define $\Omega^{0+}(\vec e)=\Omega^0(e) \cup
\Omega^{+}(\vec e)$ (resp. $\Omega^{0-}(\vec e)=\Omega^0(e) \cup
\Omega^{-}(\vec e)$).

\vskip 5pt

\noindent {\bf $\mu$-Markoff triples.}\,\, For a complex number
$\mu$, a $\mu$-{\it Markoff triple} is an ordered triple $(x,y,z)$
of complex numbers satisfying the $\mu$-Markoff equation:
\begin{eqnarray}\label{eqn:x^2+y^2+z^2=xyz}
x^2+y^2+z^2-xyz=\mu.
\end{eqnarray}
Thus Markoff triples are just $0$-Markoff triples. The set of
$\mu$-Markoff triples is also called the {\it relative character
variety}, see \cite{goldman2003gt}.

\vskip 10pt

It is easily verified that if $(x,y,z)$ is a $\mu$-Markoff triple,
so are $(x,y,xy-z)$, $(x,xz-y,z)$, $(yz-x,y,z)$ and the
permutation triples of each of them.

\vskip 10pt

\noindent {\bf $\mu$-Markoff maps.}\,\, A $\mu$-{\it Markoff map}
is a function $\phi : \Omega \rightarrow \mathbf C$ such that
\begin{itemize}

\item[(i)] for every vertex $v \in V(\Sigma)$, the triple
$(\phi(X), \phi(Y), \phi(Z))$ is a $\mu$-Markoff triple, where
$X,Y,Z \in \Omega$ are the three regions meeting $v$; and

\item[(ii)] for every edge $e \in E(\Sigma)$ such that
$e \leftrightarrow (X,Y;Z,W)$, we have
\begin{eqnarray}\label{eqn:xy=z+w}
xy=z+w,
\end{eqnarray}
where $x=\phi(X), y=\phi(Y)$, $z=\phi(Z)$ and
$w=\phi(W)$.\end{itemize}

We shall use ${\bf \Phi}_{\mu}$ to denote the set of all
$\mu$-Markoff maps and lower case letters to denote the $\phi$
values of the regions, for example, $\phi(X)=x, ~\phi(Y)=y$.
Recall that in \cite{bowditch1998plms} the set of all Markoff maps
is denoted by ${\bf \Phi}$; while here it is denoted by ${\bf
\Phi}_{0}$ in our notation. We use ${\bf \Phi}$ to denote the set
of all generalized Markoff maps.

\vskip 10pt

As in the case of Markoff maps, if the edge relation
(\ref{eqn:xy=z+w}) is satisfied along all edges, then it suffices
that the vertex relation (\ref{eqn:x^2+y^2+z^2=xyz}) be satisfied
at a single vertex. In fact one may establish a bijective
correspondence between $\mu$-Markoff maps and $\mu$-Markoff
triples, by fixing three regions $X, Y, Z$  which meet at some
vertex $v_0$, say $X(\infty), Y(0)$ and $Z(1)$. This process may
be inverted by constructing a tree of $\mu$-Markoff triples as
Bowditch did in \cite{bowditch1998plms} for Markoff triples: given
a triple $(x,y,z)$, set $\phi(X)=x, \phi(Y)=y, \phi(Z)=z$, and
extend over $\Omega$ as dictated by the edge relations. In this
way one obtains an identification of ${\bf \Phi}_{\mu}$ with the
relative character variety in $\mathbf C^3$ given by the
$\mu$-Markoff equation. In particular, ${\bf \Phi}_{\mu}$ gets a
nice topology as a subset of $\mathbf C^3$.

\vskip 5pt

The natural action of ${\rm PSL}(2, \mathbf Z)$ on $\Sigma$
induces an action on ${\bf \Phi}_{\mu}$, given by
$$H(\phi)(X)=\phi(H_*(X)),$$
where $H \in {\rm PSL}(2, \mathbf Z)$, $\phi \in {\bf
\Phi}_{\mu}$, $X$ is any element of $\Omega$ and $H_*$ is the
induced action of $H$ on $\Omega$. There is also an action of the
Klein-four group, $\mathbf Z_{2}^{2}$, on ${\bf \Phi}_{\mu}$
obtained by changing two of the signs in a $\mu$-Markoff triple,
for example, $(x,y,z) \mapsto (-x,-y,z)$. (We get the same action,
up to automorphisms of $\mathbf Z_{2}^{2}$, no matter at which
vertex we choose to perform this operation.) This action is free
and properly discontinuous on ${\bf \Phi}_{\mu} \backslash \{(\pm
\sqrt{\mu},0,0) , (0,\pm \sqrt{\mu}, 0), (0,0,\pm \sqrt{\mu})\}$.

\vskip 10pt

\noindent {\bf Natural correspondence: ${\mathcal X}_{\mu-2}
\equiv {\bf \Phi}_{\mu}$.}\,\, There is also a natural
correspondence between conjugacy classes of
$(\mu-2)$-representations $\rho$ and $\mu$-Markoff maps, obtained
by fixing a generating pair $a,b$ (and $ab$) for $\Gamma$. First
note that by fixing $a$ and $b$, we get a correspondence between
$\Omega$ and $\hat \Omega$ via the correspondence with ${\mathbf
Q}\cup \{\infty\}$ (see \S \ref{s:notation+results} and ealier
discussion in this section). Note that $X(p/q), Y(r/s) \in \Omega$
share an edge if and only if $p/q$ and $r/s$ are Farey neighbors,
if and only if the corresponding elements $[g_1], [g_2] \in\hat
\Omega$ correspond to simple closed curves on ${\mathbb T}$ which
have geometric intersection number one, or equivalently, have
representatives $g_1$ and $g_2$ which are a pair of free
generators for $\Gamma$.

\vskip 5pt

A $(\mu -2)$-representation $\rho: \Gamma \rightarrow {\rm SL}(2,
\mathbf C)$ naturally gives rise to a $\mu$-Markoff map $\phi \in
{\bf \Phi}_{\mu}$ via the correspondence between $\Omega$ and
$\hat \Omega$ by taking the traces of $[g] \in \hat \Omega$, that
is, by $\phi(X)={\rm tr}\rho(g)$ where $[g] \in \hat \Omega$
represents the simple closed curve corresponding to $X \in
\Omega$. The edge and vertex relations follow from the trace
identities in ${\rm SL}(2, \mathbf C)$:
\begin{eqnarray}\label{eqn:trAtrB=}
{\rm tr}A \,{\rm tr}B={\rm tr}AB+{\rm tr}AB^{-1},
\end{eqnarray}
\begin{eqnarray}\label{eqn:2+tr[A,B]=}
2+{\rm tr}[A,B]=({\rm tr}A)^2+({\rm tr}B)^2+({\rm tr}AB)^2-{\rm
tr}A \,{\rm tr}B \,{\rm tr}AB.
\end{eqnarray}
Representations conjugate in ${\rm SL}(2, \mathbf C)$ give rise to
the same $\mu$-Markoff map.

\vskip 3pt

Conversely, given any $\mu$-Markoff map $\phi$, we can recover the
$(\mu -2)$-representation $\rho: \Gamma \rightarrow {\rm SL}(2,
\mathbf C)$ up to conjugacy as follows. Consider the three regions
$X(\infty)$, $Y(0)$ and $Z(1)$ which meet at the vertex $v_0$ and
consider the $\mu$-Markoff triple $(x,y,z)$, where $x=\phi(X)$,
$y=\phi(Y)$ and $z=\phi(Z)$. Then one can find $A,B \in {\rm
SL}(2, \mathbf C)$, unique up to simultaneous conjugacy, such that
${\rm tr}A=x$, ${\rm tr}B=y$ and ${\rm tr}AB=z$; a specific choice
is given by Bowditch in \S 4 of \cite{bowditch1998plms}, see also
\cite{goldman2003gt} for a more natural choice. This gives a
$(\mu-2)$-representation $\rho$ with $\rho(a)=A$, $\rho(b)=B$ and
hence ${\rm tr}\rho(a)=x,{\rm tr}\rho(b)=y,{\rm tr}\rho(ab)=z$.

\vskip 5pt

In this way we can identify ${\bf \Phi}_{\mu}$ naturally with the
set ${\mathcal X}_{\mu-2}$ of conjugacy classes of $(\mu
-2)$-representations of $\Gamma$ into ${\rm SL}(2, \mathbf C)$,
and study the latter via the former. Also, there is a natural
identification of the mapping class group ${\mathcal MCG}$ with
$\PSLtwoZ$, and the action of ${\mathcal MCG}$ on ${\mathcal
X}_{\mu-2}$ corresponds exactly to the action of $\PSLtwoZ$ on
${\bf \Phi}_{\mu}$.

A key observation is that a $\mu$-Markoff map gives a tremendous
amount of information about the corresponding representation
$\rho$, since it incorporates in its definition the action of the
mapping class group on $\rho$, as well as the structure of the
mapping class group itself (which is encoded in the tree
$\Sigma$).

\vskip 10pt

\vskip 5pt

\noindent {\bf The subsets $\Omega_{\phi}(k) \subset \Omega$.}\,\,
Given $\phi \in {\bf \Phi}_{\mu}$ and $k \ge 0$, the set
${\Omega}_{\phi} (k) \subseteq \Omega$ is defined by %
\begin{eqnarray*}
{\Omega}_{\phi}(k)=\{ X \in \Omega \mid |\phi(X)| \le k \}.
\end{eqnarray*}

\vskip 5pt We now state Bowditch's Q-conditions for $\mu$-Markoff
maps.

\vskip 5pt \noindent {\bf BQ-conditions for $\mu$-Markoff
maps.}\,\, A generalized Markoff map $\phi \in {\bf \Phi}_{\mu}$
is said to satisfy the {\it BQ-conditions} if

\vskip 3pt

(BQ1)\,\, $\phi^{-1}([-2,2])=\emptyset$; and

(BQ2)\,\, $\Omega_\phi (2)$ is finite (possibly empty).

\vskip 5pt

We denote by $({\bf \Phi}_{\mu})_{Q}$ the set of all generalized
$\mu$-Markoff maps which satisfy the BQ-conditions, and call the
set of all such maps the Bowditch $\mu$-Markoff maps. It is clear
that $({\bf \Phi}_{\mu})_{Q}$ corresponds to $({\mathcal
X}_{\mu-2})_Q$.

\vskip 10pt

In the rest of this section we  examine the results and arguments
given in \S 1--\S 4 of \cite{bowditch1998plms} about Markoff maps
and point out those which  hold for generalized Markoff maps. A
key observation here is that most of the basic results derived
there only depend on the edge relation (\ref{eqn:xy=z+w}) and not
on the vertex relation (\ref{eqn:x^2+y^2+z^2=xyz}) and hence holds
for generalized Markoff maps as well.

\vskip 5pt

\noindent {\bf Results for generalized Markoff maps.}\,\, We first
state the results for generalized Markoff maps corresponding to
those in \S\S 1--3 of \cite{bowditch1998plms} for Markoff maps.

\vskip 5pt

\begin{thm}\label{thm:B1}If $\phi \in {\bf \Phi}_\mu$ then

{\rm (1)} there exists a constant $m=m(\mu)>0$ {\rm(}depending on
$\mu$ but not on $\phi${\rm)} such that ${\Omega}_{\phi}(m)$ is
non-empty; and

{\rm (2)} for any $k \ge 2$, the union $\bigcup {\Omega}_{\phi}(k)
:= {\bigcup}_{X \in {\Omega}_{\phi}(k)}X$ is connected {\rm (}as a
subset of the hyperbolic plane{\rm )}. In particular, $\bigcup
{\Omega}_{\phi}(2)$ is connected.
\end{thm}

\vskip 5pt Part (2) of Theorem \ref{thm:B1}  is a fundamental
result and is crucial to the study of the action of $\PSLtwoZ$ on
${\bf \Phi}_{\mu}$.


\vskip 5pt We have seen that $\Phi_{\mu}$ admits a nice topology
from its identification with the relative character variety
$\{(x,y,z)\in {\mathbf C}^3 \mid x^2+y^2+z^2-xyz=\mu\}$. We have
the following:

\vskip 5pt
\begin{thm}\label{thm:B3.16}
The set $({\bf \Phi}_{\mu})_{Q}$ is an open subset of ${\bf
\Phi}_{\mu}$.
\end{thm}


%

\vskip 5pt An important innovation in \cite{bowditch1998plms} was
to study the growth rate of a $\mu$-Markoff map $\phi$ by
comparing it to a simple function $F_e:\Omega \mapsto {\mathbf N}$
called the Fibonacci function, defined for a  fixed edge $e \in
E(\Sigma)$. Recall that Bowditch defined the term
`\,${\log}^{+}|\phi|$ has Fibonacci growth' in \S 2.1
\cite{bowditch1998plms}, where ${\log}^{+}(x)=\max\{0, \log(x)\}$
for $x>0$. The definition of Fibonacci growth, lower and upper
Fibonacci bound, and the Fibonacci function $F_e$, will be given
in the later part of this section. The power of this notion lies
in that if ${\log}^{+}|\phi|$ has lower Fibonacci bound then
$\sum_{X \in \Omega}|\phi(X)|^{-t}$ converges for any $t > 0$; in
particular, $\Omega_\phi (k)$ is finite for all $k$ and so
$\phi(\Omega) \subseteq {\mathbf C}$ is discrete.

\vskip 10pt

\begin{thm}\label{thm:B2}
Suppose $\phi \in ({\bf \Phi}_{\mu})_{Q}$ and $\mu \neq 4$. Then
${\log}^{+}|\phi|$ has Fibonacci growth.
\end{thm}

\vskip 5pt


\vskip 5pt

Exactly as in \cite{bowditch1998plms}, the argument for proving
the lower Fibonacci bound in Theorem \ref{thm:B2} can be applied
to a single branch of the tree $\Sigma$.

\vskip 5pt

\begin{prop}\label{prop:B3.9} {\rm (Proposition 3.9 in
\cite{bowditch1998plms})} Suppose $\vec e$ is a directed edge of
$\Sigma$ such that $\Omega^{0-}(\vec e) \,\cap \,\Omega_\phi(2)$
is finite and $\Omega^{0-}(\vec e) \,\cap\,
\phi^{-1}[-2,2]=\emptyset$. Then $\log^{+}|\phi|$ has  lower
Fibonacci bound on $\Omega^{0-}(\vec e)$.
\end{prop}

\vskip 5pt

For $\phi \in ({\bf \Phi}_{\mu})_{Q}$, we have the following
version of the generalized McShane's identity, which is a
reformulation of our first main result, Theorem \ref{thm:TWZ} in
terms of generalized Markoff maps.

\vskip 10pt

\begin{thm}\label{thm:mcshane mu-markoff}
If $\phi \in ({\bf \Phi}_{\mu})_{Q}$ {\rm(}$\mu \neq 0,4${\rm)}
then
\begin{eqnarray}\label{eqn:mcshane mu-markoff}
\sum_{X \in \Omega} {\mathfrak h}(\phi(X))=\nu \mod 2 \pi i,
\end{eqnarray}
where $\nu=\cosh^{-1}(1-\mu/2)$ and ${\mathfrak h}={\mathfrak
h}_\tau:{\mathbf C} \backslash \{\pm \sqrt{\mu}\} \rightarrow
{\mathbf C}$ is defined as in \S {\rm\ref{s:notation+results}}
with $\tau = \mu - 2$. Moreover, the sum converges absolutely.
\end{thm}

\vskip 5pt

%

There is a version, Proposition \ref{prop:B3.13 half}, of Theorem
\ref{thm:mcshane mu-markoff} applicable to a single branch of
$\Sigma$, by defining the edge weight $\psi(\vec e)=\psi_\phi(\vec
e)$ and a specific half $\hat{\mathfrak h}(x)$ of $\mathfrak h(x)$
as follows.

%
%
%

\vskip 5pt

\noindent {\bf The function $\hat{\mathfrak h}$.}\, Given $\mu \in
\mathbf C$ with $\mu \neq 0,4$, define $\hat{\mathfrak h}:
{\mathbf C}\backslash\{0,\pm\sqrt{\mu}\} \rightarrow {\mathbf C}$
by
\begin{eqnarray}
\hat{\mathfrak h}(x)=\log
\frac{1+(e^{\nu}-1)h(x)}{\sqrt{1-\mu/x^2}} \in {\mathbf C},
\end{eqnarray}
where $\nu=\cosh^{-1}(1-\mu/2)$. It is easy to check that for all
$ x \in {\mathbf C} \backslash \{0, \pm \sqrt{\tau+2}\}$,
\begin{eqnarray}\label{eqn:2 hat frak h}
2\, \hat{\mathfrak h}(x) = {\mathfrak h}(x) \mod 2 \pi i.
\end{eqnarray}

\vskip 5pt

\noindent {\bf The edge weight $\psi(\vec e)$.}\,\, Recall that we
denote a directed edge ${\vec e} \in \vec E(\Sigma)$ by ${\vec e}
= (X,Y; \rightarrow Z)$ if $X, Y$ and $Z$ are regions such that $e
= X \cap Y$ and ${\vec e}$ points towards $Z$. For a fixed $\phi
\in {\bf \Phi}_{\mu}$ and a directed edge ${\vec e} = (X,Y;
\rightarrow Z)$ with $x,y \neq 0,\pm\sqrt{\mu}$, we define the
$\phi$-{\it weight} $\psi(\vec e)$ by
\begin{eqnarray}\label{eqn:psi(vec e)}
\psi(\vec e)={\psi}_\phi(\vec e)= %
\log \frac{1+(e^{\nu}-1)(z/xy)}{\sqrt{1-\mu/x^2}\sqrt{1-\mu/y^2}}. %
\end{eqnarray}
Note that $\psi(\vec e)=\Psi(x,y,z)$ where the function $\Psi$ is
defined later in \S \ref{s:proof of theorem mu}.

\vskip 5pt

\begin{prop}\label{prop:B3.13 half} Suppose ${\vec e} = (X,Y; \rightarrow Z)$
is a directed edge of $\Sigma$ such that $\Omega^{0-}(\vec e) \cap
\Omega_\phi(2)$ is finite and $\Omega^{0-}(\vec e) \cap
\phi^{-1}[-2,2] = \emptyset$. Then
\begin{eqnarray}\label{mu-mcshane branch half}
\psi(\vec e)=\sum_{X \in \Omega^{0}(e)} \hat{\mathfrak h}(\phi(X))
+ \sum_{X \in \Omega^{-}(\vec e)} {\mathfrak h}(\phi(X)) \mod 2 \pi i, %
\end{eqnarray}
where the infinite sum converges absolutely. \square
\end{prop}

\vskip 5pt

%
Some other dynamical properties of ${\bf \Phi}_{\mu}$
corresponding to those in \S\S 4--5 of \cite{bowditch1998plms} for
Markoff maps will be discussed in a future paper
\cite{tan-wong-zhang2004endinvariants}.

\vskip 5pt

 \noindent {\bf Examination and extension of Bowditch's
arguments.}\,\, For the rest of this section, we fix on one
$\mu$-Markoff map, $\phi$, where $\mu \neq 0,4$.

\vskip 10pt

We make the {\bf assumption} that $\phi^{-1}(0) = \emptyset$,
which is not essential, to simplify the exposition. Note that the
assumption is true for $\phi \in ({\bf \Phi}_{\mu})_{Q}$.

\vskip 5pt

We shall also adopt the following {\bf convention} of Bowditch
\cite{bowditch1998plms}: We use upper case latin letters for
elements of $\Omega$, and the corresponding lower case letters for
the values assigned to them by $\phi$; that is, $x=\phi(X), \,
y=\phi(Y)$ etc.

\vskip 10pt

The above convention and assumption will allow us to write the
edge relation, (\ref{eqn:xy=z+w}), in the following convenient
form:
\begin{eqnarray}\label{eqn:z/xy+w/xy=1}
\frac{z}{xy}+\frac{w}{xy}=1.
\end{eqnarray}

%

\vskip 5pt

\noindent {\bf Arrows assigned by a $\mu$-Markoff map.}\,\, As
Bowditch did in \cite{bowditch1998plms}, we may use $\phi \in {\bf
\Phi}_{\mu}$ to assign to each undirected edge, $e$, a particular
directed edge, $\alpha_\phi (e)$, with underlying edge $e$.
Suppose $e \leftrightarrow (X,Y;Z,W)$. If $|z|>|w|$ then the arrow
on $e$ points towards $W$; in other words, $\alpha_\phi (e) =
(X,Y;Z \rightarrow W)$. Note that the statement is equivalent to
$\Re \big(\frac{z}{xy}\big) > \frac{1}{2} $. In particular, it
implies that $2|z|>|xy|$. If $|z|<|w|$, we put an arrow on $e$
pointing towards $Z$, that is, $\alpha_\phi (e) = (X,Y;W
\rightarrow Z)$. If it happens that $|z|=|w|$ then we choose
$\alpha_\phi(e)$ arbitrarily.

\vskip 5pt

The following elementary lemma, which is Lemma 3.2(3) of
\cite{bowditch1998plms}, is fundamental. It generalizes to
$\mu$-Markoff maps as the proof only uses the edge relation
(\ref{eqn:xy=z+w}).


\begin{lem}\label{lem:B3.2(3)} {\rm (Lemma 3.2(3),
\cite{bowditch1998plms})} Suppose $X,Y,Z \in \Omega$ meet at a
vertex $v \in V(\Sigma)$, and that the arrows on the edges $X \cap
Y$ and $X \cap Z$ both point away from $v$. Then $|x| \le 2$.

\end{lem}

\begin{pf}
Let $y'$ be the value of $\phi$ on the region opposite $Y$. Then
from the direction of the arrow, $2|y|\ge |y|+|y'|\ge|y+y'|=|xz|$.
Similarly, $2|z| \ge |xy|$, from which the inequality follows.
\end{pf}

\noindent {\bf Proof of Theorem \ref{thm:B1}(2).} \, This is now a
simple exercise. Assume that the result is false, choose a minimal
path on $\Sigma$ joining two connected components. If the path
consists of only one edge, we get a contradiction using
(\ref{eqn:xy=z+w}). If there is more than one edge on this
connecting path, we observe that the two ends point outwards,  and
we derive a contradiction using Lemma \ref{lem:B3.2(3)}, see
\cite{bowditch1998plms} for details. \square

\vskip 5pt

Note that Lemma 3.2(1), \cite{bowditch1998plms} (which states that
there are no sources for $\phi \in {\bf \Phi}_0 \backslash \{{\bf
0}\}$) will {\it no longer} hold for general generalized Markoff
maps. For example, if $x \in [-2,2]$, then the generalized Markoff
map which corresponds to the triple $(x,x,x)$ at vertex $v$ has a
source at $v$. On the other hand, Lemma 3.2(2) of
\cite{bowditch1998plms} can be modified as follows.

\vskip 5pt

\begin{lem}\label{lem:B3.2(2)} {\rm (Lemma 3.2(2) in
\cite{bowditch1998plms})} There is a constant $m(\mu)>0$ such that
if three regions $X,Y,Z$ meet at a sink, then
\begin{eqnarray}
\min\{|x|,|y|,|z|\} \le m(\mu).
\end{eqnarray}
\end{lem}

\vskip 2pt

\begin{pf} We show that if $|x|,|y|,|z|$ are all sufficiently large
then the vertex, $v$, that the regions $X,Y,Z$ meet cannot be a
sink. We may assume $x,y,z \neq 0$ and $|x| \le |y| \le |z|$. We
can rewrite (\ref{eqn:x^2+y^2+z^2=xyz}) as:
\begin{eqnarray}
\frac{z}{xy}+\frac{x}{yz}+\frac{y}{zx}=1+\frac{\mu}{xyz},
\end{eqnarray}
hence
\begin{eqnarray}
\Re\Big(\frac{z}{xy}\Big)+\Re\Big(\frac{x}{yz}\Big)+\Re\Big(\frac{y}{zx}\Big)=1+\Re\Big(\frac{\mu}{xyz}\Big),
\end{eqnarray}
Since $|x|,|y|,|z|$ are all sufficiently large and $\mu$ is fixed,
$\big| \frac{x}{yz} \big|, \big| \frac{y}{zx} \big|$ and $\big|
\frac{\mu}{xyz} \big|$ are all sufficiently small and so are
$\big|\Re \big(\frac{x}{yz}\big) \big|, \big|\Re
\big(\frac{y}{zx}\big) \big|$ and $\big|\Re
\big(\frac{\mu}{xyz}\big) \big|$. It follows that
$\Re\big(\frac{z}{xy}\big) > \frac12$. Hence the arrow on the edge
$X \cap Y$ is directed away from $v$ and therefore $v$ is not a
sink.
\end{pf}

It seems difficult to determine the optimal $m(\mu)$ for general
$\mu$, although it was proven in Lemma \ref{lem:B3.2(2)},
\cite{bowditch1998plms} that the optimal $m(0)=3$.

%

\vskip 10pt


\noindent {\bf Neighbors around a region.}\,\, For each $X \in
\Omega$, its boundary $\partial X$ is a bi-infinite path
consisting of a sequence of edges of the form $X \cap Y_n$, where
$(Y_n)_{n \in \mathbf Z}$ is a bi-infinite sequence of
complementary regions. Let $x=\lambda + \lambda^{-1}$ where
$|\lambda| \ge 1$. If $x=2$, then the vertex relation tells us
that $y_{n+1}=y_n \pm \sqrt{\mu-4}$, and the edge relation tells
us that $y_{n+1}-y_n=y_n-y_{n-1}$, hence the $\pm$ sign is
constant in $n$. Similarly, if $x=-2$, then $y_{n+1}=-y_n \pm
\sqrt{\mu-4}$, but this time, $y_{n+1}-y_n=-(y_n-y_{n-1})$, hence
the $\pm$ sign alternates in $n$. If $x=\pm \sqrt{\mu}$, then
$y_{n+1}= (1/2)(x \pm \sqrt{\mu-4})y_n$ where the $\pm$ sign is
constant in $n$. If $x \notin \{\pm 2, \pm \sqrt{\mu}\}$ then
there are constants $A,B \in \mathbf C \backslash \{0\}$ with
$AB=(x^2-\mu)/(x^2-4)$ such that $y_n=A {\lambda}^n + B
{\lambda}^{-n}$. Note that $|\lambda|=1$ if and only if $x \in
[-2,2] \subseteq \mathbf R$. Hence we deduce that the following
holds. (This is Corollary 3.3 in \cite{bowditch1998plms} in the
case $\mu=0$.)

\vskip 10pt

\begin{lem}\label{lem:B3.3} Suppose that $X\in \Omega$ has
neighboring regions $Y_n$, $n \in {\mathbf Z}$. Then
\begin{enumerate}
\item If $x \notin [-2,2] \cup \{\pm \sqrt{\mu}\}$, then $|y_n|$
grows exponentially as $n \rightarrow \infty$ and as $n
\rightarrow -\infty$.

\item If $x \in (-2,2)$, then $|y_n|$ remains bounded.

\item If $x=2$, then either $y_n=y_0+n\sqrt{\mu-4}$ for all $n$,
or $y_n=y_0-n\sqrt{\mu-4}$ for all $n$.

\item If $x=-2$, then either $y_n=(-1)^{n}(y_0+n\sqrt{\mu-4})$ for
all $n$, or $y_n=(-1)^{n}(y_0-n\sqrt{\mu-4})$ for all $n$.

\item If $x=\pm \sqrt{\mu}$, then either
$y_n=y_0[(x+\sqrt{\mu-4})/2]^n$ for all $n$, or
$y_n=y_0[(x-\sqrt{\mu-4})/2]^n$ for all $n$.

\end{enumerate}
\square
\end{lem}

\vskip 5pt

Using part (5) of the Lemma above, we can show:\vskip 8pt

\begin{lem}\label{lem:x^2=mu}
If $\phi \in {\bf \Phi}_{\mu}$ and $\phi(X)=\pm \sqrt{\mu}$ for
some $X \in \Omega$, then $\phi \notin ({\bf \Phi}_{\mu})_{Q}$.
\end{lem}

\vskip 2pt

\begin{pf}
Without loss of generality, we may assume that
$y_n=y_0[(x+\sqrt{\mu-4})/2]^n$ for all $n$. If $y_0=0$ then $\phi
\notin ({\bf \Phi}_{\mu})_{Q}$. Suppose $y_0 \neq 0$ and let
$\lambda =(x+\sqrt{\mu-4})/2$. Then it is easy to see that
$\lambda \neq 0$. If $|\lambda|=1$, then ${\lambda}^{-1}
=(x-\sqrt{\mu-4})/2$. Hence $x=\lambda + {\lambda}^{-1} \in
[-2,2]$, again $\phi \notin ({\bf \Phi}_{\mu})_{Q}$. Now suppose
$|\lambda| \neq 1$. If $|\lambda| < 1$ then $|y_n| \rightarrow 0$
as $n \rightarrow \infty$, whereas if $|\lambda| > 1$ then $|y_n|
\rightarrow 0$ as $n \rightarrow -\infty$. Hence
${\Omega}_{\phi}(2)$ is infinite and $\phi \notin ({\bf
\Phi}_{\mu})_{Q}$.
\end{pf}

The following lemma, together with Lemma \ref{lem:B3.2(3)}, will
imply that if $\Omega_\phi(2)=\emptyset$ then there will be a
sink.

\vskip 8pt

\begin{lem}\label{lem:B3.4} {\rm (Lemma 3.4 in
\cite{bowditch1998plms})} Suppose $\beta$ is an infinite ray in
$\Sigma$ consisting of a sequence, $(e_n)_{n \in \mathbb N}$, of
edges of $\Sigma$ such that the arrow on each $e_n$ assigned by
$\phi$ is directed towards $e_{n+1}$. Then $\beta$ meets at least
one region $X$ with $|\phi(X)|<2$.
\end{lem}

 \vskip 2pt

\begin{pf}
The same  proof used for Lemma 3.4 of \cite{bowditch1998plms}
works here  with a slight refinement for its last part as follows.
Suppose that all regions $X$ which are incident to $\beta$ have
$|\phi(X)| \ge 2$. Then by the argument there, for $n$
sufficiently large, there are directed edges ${\alpha_\phi}(e_n) =
(X,Y;Z \rightarrow W)$ in $\beta$ with $X \cap Z \subseteq \beta$
and $Y \cap W \subseteq \beta$ so that (where $\simeq$ means
`arbitrarily close to') $z/xy \simeq 1/2, x^2 \simeq 4$ and $y^2
\simeq 4$. Hence $z^2 \simeq 4$ and $xyz \simeq 2z^2 \simeq 8$.
Thus $\mu = x^2 + y^2 + z^2 -xyz \simeq 4$ which is impossible
since we have assumed that $\mu \neq 4$.
\end{pf}

\vskip 5pt

Note that if the path $\beta$ in Lemma \ref{lem:B3.4} above does
not eventually lie in the boundary of a fixed region, then $\beta$
meets infinitely many regions $X$ with $|\phi(X)|<2$.

\vskip 5pt

 \noindent {\bf Proof of Theorem \ref{thm:B1}(1).} \, If
$\Omega_\phi(2)=\emptyset$, then Lemmas \ref{lem:B3.2(3)} and
\ref{lem:B3.4} tell us that there must be a sink. We now apply
Lemma \ref{lem:B3.2(2)}. This proves Theorem \ref{thm:B1}(1).
\square

\vskip 15pt

\noindent {\bf Fibonacci function.}\,\, We recall now the
definition given in \cite{bowditch1998plms} for the Fibonacci
function $F_e$ associated to an edge $e \in E(\Sigma)$, and the
definition for a function $f:\Omega \rightarrow [0, \infty)$ to
have an upper/lower Fibonacci bound.

\vskip 10pt

%
%
%

Let $\vec e \in \vec E(\Sigma)$ with underlying edge $e$. For $X
\in \Omega^{0-}(\vec e)$ we define $d(X)=d_{\vec e}(X)$ to be the
number of edges in the shortest path joining the head of $\vec e$
to $X$. Given any $Z \in \Omega^{-}(\vec e)$, there are precisely
two regions $X,Y \in \Omega^{0-}(\vec e)$ meeting $Z$ and
satisfying $d(X)<d(Z)$ and $d(Y)<d(Z)$. Note that $X,Y,Z$ all meet
in a vertex.

\vskip 5pt

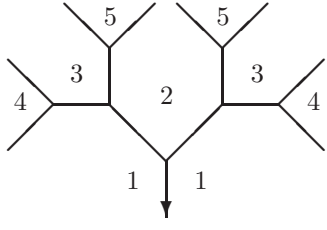
\begin{figure}
\setlength{\unitlength}{0.75mm}
\begin{picture}(80,40)
\thicklines \put(40,10){\vector(0,-1){10}}
\put(40,10){\line(1,1){10}} \put(40,10){\line(-1,1){10}}
\put(50,20){\line(0,1){10}} \put(50,20){\line(1,0){10}}
\put(30,20){\line(0,1){10}} \put(30,20){\line(-1,0){10}}
\put(50,30){\line(1,1){8}} \put(50,30){\line(-1,1){8}}
\put(30,30){\line(1,1){8}} \put(30,30){\line(-1,1){8}}
\put(60,20){\line(1,1){8}} \put(60,20){\line(1,-1){8}}
\put(20,20){\line(-1,-1){8}} \put(20,20){\line(-1,1){8}}
\put(33,5){$1$} \put(45,5){$1$} \put(39,20){$2$} \put(23,24){$3$}
\put(55,24){$3$} \put(29,34){$5$} \put(49,34){$5$}
\put(13,19){$4$} \put(65,19){$4$}
\end{picture}
\caption{Some values of the Fibonacci function $F_e$}
\label{fig:Fe}
\end{figure}

Now we can define the {\it Fibonacci function} $F_e: \Omega
\rightarrow {\mathbf N}$ with respect to an edge $e$ as follows.
We orient $e$ arbitrarily as $\vec e$ and define $F_{\vec e}:
\Omega^{0-}(\vec e) \rightarrow {\mathbf N}$ by $F_{\vec e}(W)=1$
for $W \in \Omega^{0}(e)$ and $F_{\vec e}(Z)= F_{\vec
e}(X)+F_{\vec e}(Y)$ for $Z \in \Omega^{-}(\vec e)$ where $X,Y \in
\Omega^{0-}(\vec e)$ are the two regions meeting $Z$ and
satisfying $d(X)<d(Z)$ and $d(Y)<d(Z)$. Now we define
$F_e(X)=F_{\vec e}(X)$ for $X \in \Omega^{0-}(\vec e)$ and
$F_e(X)=F_{-\vec e}(X)$ for $X \in \Omega^{+}(\vec e)$. See Figure
\ref{fig:Fe} where the values $F_e(X)$ are given for some $X \in
\Omega$.


\vskip 10pt

The functions $F_e$ provide a means for measuring the growth rates
of functions defined on subsets of $\Omega$. The following lemma
can be easily proved by induction. Its corollary shows that the
concept of upper and lower Fibonacci bound is independent of the
edge $e$ used.

\begin{lem}\label{lem:B2.1.1}{\rm (Lemma 2.1.1 in \cite{bowditch1998plms})}
Suppose $f: \Omega^{0-}(\vec e) \rightarrow [0, \infty)$ where
$\Omega^{0}(e)=\{X_1,X_2\}$.

{\rm (1)} If $f$ satisfies $f(Z) \le f(X)+f(Y)+c$ for some fixed
constant $c$ and arbitrary $X, Y, Z \in \Omega^{0-}(\vec e)$
meeting at a vertex and satisfying $d(X)<d(Z)$ and $d(Y)<d(Z)$,
then $f(X) \le (M+c)F_e(X)-c$ for all $X \in \Omega^{0-}(\vec e)$,
where $M=\max\{f(X_1), f(X_2)\}$.

{\rm (2)} If $f$ satisfies $f(Z) \ge f(X)+f(Y)-c$ for some fixed
constant $c$ where $0<c<m=\min\{f(X_1), f(X_2)\}$ and arbitrary
$X,Y,Z$ as in part {\rm (1)}, then $f(X) \ge (m-c)F_e(X)+c$ for
all $X \in \Omega^{0-}(\vec e)$. \square
\end{lem}

\vskip 8pt

\begin{cor}\label{cor:B2.1.2}{\rm (Corollary 2.1.2 in \cite{bowditch1998plms})}
Suppose $f: \Omega \rightarrow [0, \infty)$ satisfies an
inequality of the form $f(Z) \le f(X)+f(Y)+c$ for some fixed
constant $c$, whenever $X,Y,Z \in \Omega$ meet at a vertex. Then
for any given edge $e \in E(\Sigma)$, there is a constant $K>0$,
such that $f(X) \le KF_e(X)$ for all $X \in \Omega$. \square
\end{cor}

\vskip 8pt

Note that for any edge $e' \in E(\Sigma)$, $f=F_{e'}$ satisfies
the hypotheses of Corollary \ref{cor:B2.1.2}, with $c=0$. Thus for
any two edges $e,e' \in E(\Sigma)$, there is some constant
$K=K(e,e')>0$ such that
\begin{eqnarray*}
K^{-1}F_e(X) \le F_{e'}(X) \le KF_e(X)
\end{eqnarray*}
for all $X \in \Omega$. Hence the properties in the following
definitions are independent of the choices the edge $e \in
E(\Sigma)$.

\vskip 10pt

\noindent {\bf Fibonacci bounds.}\,\, Suppose $f: \Omega
\rightarrow [0, \infty)$, and $\Omega^{\prime} \subseteq \Omega$.
We say that \vskip 3pt

$f$ has an {\it upper Fibonacci bound on} $\Omega^{\prime}$ if
there is some constant $\kappa > 0$ such that $f(X) \le \kappa\,
F_e(X)$ for all $X \in \Omega^{\prime}$;

$f$ has an {\it lower Fibonacci bound on} $\Omega^{\prime}$ if
there is some constant $\kappa > 0$ such that  $f(X) \ge \kappa\,
F_e(X)$ for all but finitely many $X \in \Omega^{\prime}$;

$f$ has {\it Fibonacci growth} on $\Omega^{\prime}$ if it has both
upper and lower Fibonacci bounds on $\Omega^{\prime}$; and

$f$ has {\it Fibonacci growth} if $f$ has Fibonacci growth on {\it
all} of $\Omega$.

\vskip 10pt

Note that if $\Omega^{\prime}$ is the union of a finite set of
subsets $\Omega_1, \cdots, \Omega_m \subseteq \Omega$, then $f$
has an upper (lower) Fibonacci bound if and only if it has an
upper (lower) Fibonacci bound on each of $\Omega_i$.

\vskip 10pt

\noindent {\bf Upper Fibonacci bounds.}\,\, The following lemma
tells us that for an arbitrary $\mu$-Markoff map $\phi$, the
function $\log^{+}|\phi|$ always has an upper Fibonacci bound on
$\Omega$. Hence we only need to consider criteria for it to have a
lower Fibonacci bound on certain branches of the binary tree
$\Sigma$.

\vskip 10pt

\begin{lem}\label{lem:upper F bound}
If $\phi \in {\Phi}_{\mu}$, then $\log^{+}|\phi|$ has an upper
Fibonacci bound on $\Omega$.
\end{lem}

\vskip 5pt

\begin{pf} By Corollary \ref{cor:B2.1.2} we only need to show that for
an arbitrary $\mu$-Markoff map $(x,y,z)$
\begin{eqnarray}\label{eqn:upper F bound}
\log^{+}|z| \le \log 4 + \log^{+}|\mu| + \log^{+}|x| +
\log^{+}|y|.
\end{eqnarray}

If $|z| \le 2|x|$ or $|z| \le 2|y|$ then (\ref{eqn:upper F bound})
holds already. So we suppose $|z| \ge 2|x|$ and $|z| \ge 2|y|$.
Then since $\mu + xyz = x^2 + y^2 + z^2$ we have
\begin{eqnarray*}
|\mu|+|xyz| &\ge& |z|^2-|x|^2-|y|^2 \\
&=& |z|^2/2 + (|z|^2/4-|x|^2) + (|z|^2/4-|y|^2) \\
&\ge& |z|^2/2.
\end{eqnarray*}
Hence $|z|^2 \le 4|xyz|$ or $|z|^2 \le 4|\mu|$ according to
whether $|\mu| \le |xyz|$ or $|xyz| \le |\mu|$. Thus
\begin{eqnarray*}
|z| \le 4|xy| \,\,\, {\rm or} \,\,\, |z| \le |z|^2 \le 4|\mu|
\end{eqnarray*}
(since we may assume $|z| \ge 1$) from which (\ref{eqn:upper F
bound}) follows easily.
\end{pf}

\vskip 8pt

\noindent {\bf Lower Fibonacci bounds.}\,\, The lower Fibonacci
bounds are more interesting since, as the following proposition
shows, they give the convergence of certain series, in particular,
the absolute convergence of series on the left hand side of
(\ref{eqn:mcshane mu-markoff}) since we have $|{\mathfrak
h}(x)|=O(|x|^{-2})$ as $|x| \rightarrow \infty$.

\vskip 5pt

\begin{prop}\label{prop:B2.1.4}{\rm (Proposition 2.1.4 in \cite{bowditch1998plms})}
If $f: \Omega \rightarrow [0, \infty)$ has a lower Fibonacci
bound, then ${\sum}_{X\in\Omega}f(X)^{-s}$ converges for all $s>2$
{\rm (}after excluding a finite subset of $\Omega$ on which $f$
takes the value $0${\rm )}. \square
\end{prop}

\vskip 5pt

\begin{cor}\label{cor:sum phi^-2 convergence}
If $\phi \in ({\Phi}_{\mu})_{Q}$, then for any $t>0$, the series
$\sum_{X \in \Omega}|{\phi}(X)|^{t}$ converges absolutely. \square
\end{cor}

\vskip 5pt

The following lemma and corollary of \cite{bowditch1998plms} hold
with the same proofs there, giving a criterion for
$\log^{+}|\phi|$ to have a lower Fibonacci bound on certain
branches of $\Sigma$.

\vskip 5pt

\begin{lem}\label{lem:B3.5} {\rm (Lemma 3.5 and Corollary 3.6 in
\cite{bowditch1998plms})} Suppose $\vec e \in \vec E(\Sigma)$ is
such that $\alpha_\phi(e)=\vec e$ and $\Omega^0(e) \cap
\Omega_\phi(2)=\emptyset$. Then $\Omega^{0-}(\vec e) \cap
\Omega_\phi(2) = \emptyset$ and the arrow on each edge of
$\Sigma^{-}$ is directed towards $e$.

Furthermore, $\log |\phi(X)| \ge (m - \log 2)F_e(X)$ for all $X
\in \Omega^{0-}(\vec e)$, where $m=\min \{\log |\phi(X)| \mid X
\in \Omega^0(e)\} > \log 2$. \square
\end{lem}

\vskip 3pt

\begin{cor}\label{cor:B3.7} {\rm (Corollary 3.7 in
\cite{bowditch1998plms})} If $\Omega_\phi(2)=\emptyset$, then
there is a unique sink, and $\log^{+}|\phi|$ has a lower Fibonacci
growth. \square
\end{cor}

\vskip 5pt

To prove Theorem \ref{thm:B2} we expand somewhat on Lemma
\ref{lem:B3.5} and consider the case where $\Omega^{-}(\vec e)
\cap \Omega_\phi(2)=\emptyset$ and exactly one of the two regions
in $\Omega^0(e)$ has norm no greater than $2$.

\vskip 5pt

\begin{lem}\label{lem:B3.8}{\rm (Lemma 3.8 in \cite{bowditch1998plms})}
Suppose $\vec e \in \vec E(\Sigma)$ is such that
$\alpha_\phi(e)=\vec e$ and $\Omega^{0-}(\vec e) \cap
\Omega_\phi(2)=\{X_0\}$ where $X_0 \in \Omega^0(\vec e)$ with $x_0
\notin [-2,2]$. Then $\log^+|\phi|$ has a lower Fibonacci bound on
$\Omega^{0-}(\vec e)$.
\end{lem}

\vskip 3pt

\begin{pf} The proof of Lemma 3.8 in \cite{bowditch1998plms}
applies. We repeat it here. Let $(\vec e_n )_{n=0}^{\infty}$ be
the sequence of directed edges lying in the boundary of $X_0$ and
in $\Omega^{0-}(\vec e)$ so that $\vec e_0=\vec e$ and $\vec e_n$
is directed away from $\vec e_{n+1}$. For $n \ge 1$, let $v_n$ be
the vertex incident on both $e_{n-1}$ and $e_{n}$, and let $\vec
\varepsilon_n$ be the third edge (distinct from $e_{n-1}$ and
$e_{n}$) incident on $v_n$ and directed towards $v_n$. For $n \ge
0$, let $Y_n$ be the region such that $Y_n \cap X_0 = e_n$. See
Figure \ref{fig:3.8} for an illustration. Thus $\Omega^{0-}(\vec
e)=\{X_0\} \cup \bigcup_{n=1}^{\infty}\Omega^{0-}(\vec
\varepsilon_n)$.

By Lemma \ref{lem:B3.3}, $|y_n|$ grows exponentially as $n
\rightarrow \infty$, and so $\log|y_n| \ge cn$ for some $c>0$.
Hence we have $\log^+|\phi(X)| \ge cnF_{\vec \varepsilon_n}(X)$
for all $n \ge 1$ and for all $X \in \Omega^{0-}(\vec
\varepsilon_n)$. Thus it follows easily (using Lemma
\ref{lem:B2.1.1}(2)) that $\log^+|\phi|$ has a lower Fibonacci
bound on $\Omega^{0-}(\vec e)$.
\end{pf}

\begin{figure}
\begin{center}\mbox{
\hskip -0.05in
\beginpicture

\setcoordinatesystem units <0.55in,0.55in> \setplotarea x from
-0.5 to 8.5, y from -0.3 to 2.8

\plot 0 1 8 1 /

\arrow <6pt> [.16,.6] from 0.5 1  to 0.4 1

\arrow <6pt> [.16,.6] from 1.5 1  to 1.4 1

\arrow <6pt> [.16,.6] from 2.5 1  to 2.4 1

\arrow <6pt> [.16,.6] from 3.5 1  to 3.4 1

\arrow <6pt> [.16,.6] from 4.5 1  to 4.4 1

\arrow <6pt> [.16,.6] from 5.5 1  to 5.4 1

\arrow <6pt> [.16,.6] from 6.5 1  to 6.4 1

\arrow <6pt> [.16,.6] from 7.5 1  to 7.4 1

\plot 1 2 1 1 /  \plot 2 2 2 1 /

\plot 3 2 3 1 /  \plot 4 2 4 1 /

\plot 5 2 5 1 /  \plot 6 2 6 1 /

\plot 7 2 7 1 /

\arrow <6pt> [.16,.6] from 1 1.5 to 1 1.4

\arrow <6pt> [.16,.6] from 2 1.5 to 2 1.4

\arrow <6pt> [.16,.6] from 3 1.5 to 3 1.4

\arrow <6pt> [.16,.6] from 4 1.5 to 4 1.4

\arrow <6pt> [.16,.6] from 5 1.5 to 5 1.4

\arrow <6pt> [.16,.6] from 7 1.5 to 7 1.4

\arrow <6pt> [.16,.6] from 6 1.5 to 6 1.4

\plot 0.7 2.3 1 2 1.3 2.3 /

\plot 1.7 2.3 2 2 2.3 2.3 /

\plot 2.7 2.3 3 2 3.3 2.3 /

\plot 3.7 2.3 4 2 4.3 2.3 /

\plot 4.7 2.3 5 2 5.3 2.3 /

\plot 6.7 2.3 7 2 7.3 2.3 /

\plot 5.7 2.3 6 2 6.3 2.3 /

\put {\mbox{\Huge $\cdot$}} [cc] <0mm,-0.2mm> at 0 1

\put {\mbox{\Huge $\cdot$}} [cc] <0mm,-0.2mm> at 1 1

\put {\mbox{\Huge $\cdot$}} [cc] <0mm,-0.2mm> at 2 1

\put {\mbox{\Huge $\cdot$}} [cc] <0mm,-0.2mm> at 3 1

\put {\mbox{\Huge $\cdot$}} [cc] <0mm,-0.2mm> at 4 1

\put {\mbox{\Huge $\cdot$}} [cc] <0mm,-0.2mm> at 5 1

\put {\mbox{\Huge $\cdot$}} [cc] <0mm,-0.2mm> at 6 1

\put {\mbox{\Huge $\cdot$}} [cc] <0mm,-0.2mm> at 7 1

\put {\mbox{\small $\vec e$}} [ct] <0mm,-2mm> at 0.5 1

\put {\mbox{\small ${\vec e}_1$}} [ct] <0mm,-2mm> at 1.5 1

\put {\mbox{\small ${\vec e}_2$}} [ct] <0mm,-2mm> at 2.5 1


\put {\mbox{\small ${\vec e}_{n}$}} [ct] <0mm,-2mm> at 5.5 1

\put {\mbox{\small ${\vec e}_{n+1}$}} [ct] <0mm,-2mm> at 6.6 1

\put {\mbox{\small $X_0$}} [cc] <0mm,1mm> at 4 0.30

\put {\mbox{$\cdots$}} [cc] <-1mm,1mm> at 4.6 1.25

\put {\mbox{\small $Y_0$}} [cb] <0mm,1mm> at 0.5 1.15

\put {\mbox{\small $Y_1$}} [cb] <0mm,1mm> at 1.5 1.15

\put {\mbox{\small $Y_2$}} [cb] <0mm,1mm> at 2.5 1.15

\put {\mbox{$\cdots$}} [cc] <-1mm,1mm> at 3.6 1.25

\put {\mbox{\small $Y_{n}$}} [cb] <0mm,1mm> at 5.5 1.15

\put {\mbox{\small $Y_{n+1}$}} [cb] <0mm,1mm> at 6.6 1.15

\put {\mbox{\small ${\vec \varepsilon}_1$}} [lc] <1mm,0mm> at 1
1.85

\put {\mbox{\small ${\vec \varepsilon}_2$}} [lc] <1mm,0mm> at 2
1.85


\put {\mbox{\small ${\vec \varepsilon}_{n}$}} [lc] <1mm,0mm> at 5
1.85

\put {\mbox{\small ${\vec \varepsilon}_{n+1}$}} [lc] <1mm,0mm> at
6 1.85

\endpicture}\end{center}
\vskip -0.25in \caption{}\label{fig:3.8}
\end{figure}

\vskip 5pt

\noindent {\bf Proof of Theorem \ref{thm:B2}.} \, The proof of
Theorem \ref{thm:B2} is then the same as that of Theorem 2 in
\cite{bowditch1998plms}. We sketch it as follows. By Lemma
\ref{lem:upper F bound}, we only need to show that $\log^+|\phi|$
has a lower Fibonacci bound on $\Omega$. If $\Omega_\phi(2)$ has
at most one element, the conclusion follows easily by Corollary
\ref{cor:B3.7} and Lemma \ref{lem:B3.8}. Hence we suppose
$\Omega_\phi(2)$ has at least two elements.

Recall that $\Omega_\phi(2) \subseteq \Omega$ is finite and
$\bigcup\Omega_\phi(2)$ is connected. Let $T$ be the (finite)
subtree of $\Sigma$ {\it spanned} by the set of edges $e$ such
that $\Omega^0(e) \subseteq \Omega_\phi(2)$. Let $C=C(T)$ be the
circular set of directed edges given by $T$. Note that
$\Omega=\bigcup_{\vec e \in C}\Omega^{0-}(\vec e)$. Hence it
suffices to show that $\log^+|\phi|$ has a lower Fibonacci bound
on $\Omega^{0-}(\vec e)$ for every $\vec e \in C$. Then the
conclusion of Theorem \ref{thm:B2} follows by the following claim,
Lemma \ref{lem:B3.5} and Lemma \ref{lem:B3.8}.

\vskip 5pt

{\bf Claim.} \, For each $\vec e \in C$, we have $\vec e =
\alpha_\phi(e)$, $\Omega^-(\vec e) \cap \Omega_\phi(2) =
\emptyset$ and $\Omega^0(e) \cap$ $\Omega_\phi(2)$ has at most one
element.

\vskip 5pt

To prove the claim, let $\vec e = (X,Y;Z \rightarrow W) \in C(T)$.
If one of $X$ and $Y$, say $X$, is in $\Omega_\phi(2)$ then  $Y, Z
\notin \Omega_\phi(2)$ and $W \in \Omega_\phi(2)$ by the
definition of $T$. Hence in this case $\vec e = \alpha_\phi(e)$,
$\Omega^-(\vec e) \cap \Omega_\phi(2) = \emptyset$ and
$\Omega^0(e) \cap \Omega_\phi(2)$ has one element, $X$. Now
suppose neither $X$ nor $Y$ is in $\Omega_\phi(2)$ then $W \in
\Omega_\phi(2)$ and $Z \notin \Omega_\phi(2)$ since
$\bigcup\Omega_\phi(2)$ is connected. Thus in this case $\vec e =
\alpha_\phi(e)$, $\Omega^{0-}(\vec e) \cap \Omega_\phi(2) =
\emptyset$. This proves the claim, completing the proof of the
theorem. \square

\vskip 10pt

\noindent {\bf The definition of $H(x)$}.\, To prove Theorem
\ref{thm:B3.16}, that is, $({\bf \Phi}_{\mu})_{Q}$ is an open
subset of ${\bf \Phi}_{\mu}$, we introduce a continuous function
$H: \mathbf C \backslash \big([-2,2] \cup \{ \pm\sqrt{\mu}\}\big)
\rightarrow \mathbf R_{>0}$ so that for any $\phi \in {\bf
\Phi}_{\mu}$ and $X \in \Omega$, if $(Y_n)_{n \in \mathbf Z}$ is
the bi-infinite sequence of regions meeting $X$, then there are
integers $n_1 \le n_2$ such that $|y_n| \le H(x)$ if and only if
$n_1 \le n \le n_2$, and  $|y_n|$ is monotonically decreasing for
$n \in (-\infty, n _1]$ and  monotonically increasing for $n \in
[n_2,\infty)$.

\begin{lem}\label{lem:H(x)} For $x \neq \pm \sqrt{\mu}$ and $x
\notin [-2,2]$, $H(x)$ can be chosen as
\begin{eqnarray}\label{eqn:H(x)}
H(x)=\sqrt{\Big|\frac{x^2-\mu}{x^2-4}\Big|}\,\frac{2|\lambda(x)|^2}{|\lambda(x)|-1},
\end{eqnarray}
where $$\lambda(x)=\frac{x}{2} \bigg( 1+\sqrt{1-\frac{4}{x^2}}\,
\bigg)=e^{l(x)/2}.$$
\end{lem}

\begin{pf} Let us write $\lambda=\lambda(x)$. Then $|\lambda|>1$
since $x \notin [-2,2]$. As explained before Lemma \ref{lem:B3.3},
there exist $A,B$ with $AB=(x^2-\mu)/(x^2-4)$ such that
$y_n=A\lambda^n+B\lambda^{-n}$. Note that since we can replace
$A,B$ respectively by $A\lambda^m, B\lambda^{-m}$ for $m \in
\mathbf Z$, we may assume that $\sqrt{|AB|}|\lambda|^{-1} \le
|A|,|B| \le \sqrt{|AB|}|\lambda|$, that is,
\begin{eqnarray}\label{eqn:|A|,|B|}
\sqrt{\Big|\frac{x^2-\mu}{x^2-4}\Big|}\,|\lambda|^{-1} \le |A|,|B|
\le \sqrt{\Big|\frac{x^2-\mu}{x^2-4}\Big|}\,|\lambda|.
\end{eqnarray}
(Note that this argument fails for $x=\pm\sqrt{\mu}$ since in that
case we have $AB=0$.)

\vskip 5pt

In the rest of the proof we assume that $n>0$.

\vskip 5pt

{\bf Claim.} If
\begin{eqnarray}\label{eqn:lambda^2n-1}
|\lambda|^{2n-1} \ge \frac{|\lambda|+1}{|\lambda|-1}, %
\end{eqnarray}
then\, $|y_{n+1}| \ge |y_n|$ \,and\, $|y_{-n-1}|\ge|y_{-n}|$.

\vskip 5pt

To prove the claim, note that \vskip 3pt

\centerline{$|y_{n+1}| \ge |A||\lambda|^{n+1} -
|B||\lambda|^{-n-1}$ and\, $|y_{n}| \le |A||\lambda|^{n} +
|B||\lambda|^{-n}$.}

\vskip 3pt \noindent Thus \,$|y_{n+1}|\ge|y_n|$\, if
$$|A||\lambda|^{n+1} - |B||\lambda|^{-n-1} \ge |A||\lambda|^{n} + |B||\lambda|^{-n}$$
or equivalently
$$|\lambda|^{2n+1} \ge \frac{|B|}{|A|}\frac{|\lambda|+1}{|\lambda|-1}.$$
Since $|B|/|A| \le |\lambda|^2$, we have $|y_{n+1}|\ge|y_n|$ if
$|\lambda|^{2n-1} \ge (|\lambda|+1)/(|\lambda|-1)$. The other
inequality $|y_{-n-1}|\ge|y_{-n}|$ can be similarly proved.

\vskip 5pt

We continue the proof of Lemma \ref{lem:H(x)}. Note that
\begin{eqnarray*}
|y_n| &\le& |A||\lambda|^{n} + |B||\lambda|^{-n} \\
&\le& \sqrt{\Big|\frac{x^2-\mu}{x^2-4}\Big|}\,|\lambda|\,(|\lambda|^n+|\lambda|^{-n}) \\
&\le& \sqrt{\Big|\frac{x^2-\mu}{x^2-4}\Big|}\,|\lambda|\,(|\lambda|^n+1). %
\end{eqnarray*}
Hence if
$$|y_n| \ge H(x)=\sqrt{\Big|\frac{x^2-\mu}{x^2-4}\Big|}\,|\lambda|\bigg(\frac{|\lambda|+1}{|\lambda|-1}+1\bigg)$$ %
then $|\lambda|^n \ge \frac{|\lambda|+1}{|\lambda|-1}$. Hence
(\ref{eqn:lambda^2n-1}) holds. It follows from the claim that
$|y_{n+1}|\ge|y_n|$. Similarly, if $|y_{-n}| \ge H(x)$ then
$|y_{-n-1}|\ge|y_{-n}|$.
\end{pf}

\noindent {\it Remark.}\,\, We observe from the proof that we
always have $|y_0| < H(x)$.

\vskip 5pt

On the other hand, for $x \in [-2,2]$ or  $x=\pm\sqrt{\mu}$, we
set $H(x)=\infty$. By replacing $H(x)$ by $\max\{H(x),2\}$, we may
assume that $H(x) \ge 2$. Hence $H$ is continuous on $\mathbf C
\backslash \{\pm\sqrt{\mu}\}$ (but is not continuous at $\pm
\sqrt{\mu}$ if $\mu \notin [0,4]$).

\vskip 10pt

\noindent {\bf The definition of $T(t)$}.\, Now for each $\phi \in
{\bf \Phi}_{\mu}$ and $t \ge 0$, we define a subtree $T(t)$ of
$\Sigma$ as follows. First, for each $X \in \Omega$ with
$x=\phi(X) \notin [-2,2]$ and $x \neq \pm\sqrt{\mu}$, and for $r
\ge H(x) \ge 2$, we set $$J_r(X)=\bigcup \{ X \cap Y_n \mid |y_n|
\le r \}.$$  It is an subarc of $\partial X$ determined by $\phi$
with the property that if $e$ is any edge in $\partial X$ not
lying in $J_r(X)$ then the arrow on $e$ assigned by $\phi$ points
towards $J_r(X)$. Note that $J_r(X)$ contains at least one edge by
the observation we have just made after the proof of Lemma
\ref{lem:H(x)}.

\vskip 5pt

Now, for given $t \ge 0$, let $T(t)=T_\phi(t)$ be the union of all
the arcs $J_{H(x)+t}(X)$ as $X$ varies in
$\Omega(2+t)=\Omega_\phi(2+t)$.

\vskip 5pt

It can be shown (by the argument of Lemma 3.11,
\cite{bowditch1998plms}) that $T(t)$ is connected, hence a subtree
of $\Sigma$. It follows that $T(t)$ has the following nice
property.

\vskip 5pt

\begin{lem}\label{lem:T(t)property}
If $T(t)\neq\emptyset$, then the arrow on any edge not in the
subtree $T(t)$ points towards $T(t)$. \square
\end{lem}

Actually, by Theorem \ref{thm:B1}(1), we have

\vskip 5pt

\begin{lem}\label{lem:T(t)} $T(t) \neq \emptyset$ for any $\phi \in {\bf
\Phi}_\mu$ and for any $t \ge m(\mu)-2$, where $m(\mu) > 2$ is a
constant such that $\Omega_\phi(m(\mu)) \neq \emptyset$ for all
$\phi \in {\bf \Phi}_\mu$. \square
\end{lem}

\vskip 5pt

Alternatively, we can describe $T(t)$ directly in terms of its
edges.

\vskip 5pt

\begin{lem}\label{lem:T(t) by edges}
An edge $e=X \cap Y$ is an edge of $T(t)$ if and only if either
$|x|\le 2+t$ and $|y|\le H(x)+t$ or $|y|\le 2+t$ and $|x|\le
H(y)+t$. \square
\end{lem}

\vskip 5pt

Thus we have the following lemma which gives a finite criterion
for recognizing that a given $\mu$-Markoff map $\phi$ lies in
$({\bf \Phi}_{\mu})_Q$. This generalizes Lemma 3.15 in
\cite{bowditch1998plms}, and the proof there also works here.

\vskip 5pt

\begin{lem}\label{lem:B3.15} For any fixed $t\ge 0$, $\phi \in ({\bf
\Phi}_{\mu})_Q$ if and only if $T_{\phi}(t)$ is finite. \square
\end{lem}

\vskip 5pt

Now we can give the proof of Theorem \ref{thm:B3.16} by showing
that this criterion is an open property.

\vskip 10pt

\noindent {\bf Proof of Theorem \ref{thm:B3.16}.} \, The proof is
essentially the same as that of Theorem 3.16,
\cite{bowditch1998plms}, which states that $({\bf \Phi}_0)_Q$ is
an open subset of ${\bf \Phi}_{0}$. Here we supply a bit more
details.

Fix any $t_1 > m(\mu)-2$. Suppose $\phi \in ({\bf \Phi}_\mu)_Q$
and write $T(t)$ for $T_\phi(t)$. By Lemma \ref{lem:B3.15},
$T(t_1)$ is a finite subtree of $\Sigma$. We may choose $t_2 >
t_1$ large enough so that $T(t_2)$ contains $T(t_1)$ in its
interior, that is, it contains $T(t_1)$, together with all the
edges of the circular set $C(T(t_1))$. Note that $T(t_2)$ is also
a finite subtree of $\Sigma$.

For any given $\phi^{\prime} \in {\bf \Phi}_\mu$, we write $T'(t)$
for $T_{\phi^{\prime}}(t)$.

\vskip 5pt

{\bf Claim.}\, If $\phi^{\prime}$ is sufficiently close to $\phi$,
then $T'(t_1) \cap T(t_2) \subseteq T(t_1)$.

\vskip 5pt

To prove the claim, choose $e =X \cap Y \in T(t_2) \backslash
T(t_1)$. We may assume $|\phi(X)| \le 2+t_2$ and $|\phi(Y)| \le
H(x)+t_2$ since $e \in T(t_2)$. Then either $|\phi(X)| > 2+t_1$ or
$|\phi(Y)| > H(x)+t_1 \ge 2+t_1$ since $e \notin T(t_1)$. Thus if
$\phi^{\prime}$ is sufficiently close to $\phi$, we have either
$|\phi(X)| > 2+t_1$ or $|\phi(Y)| > 2+t_1$, hence $e \notin
T'(t_1)$. This proves the claim since there are only finitely many
edges $e$ in $T(t_2) \backslash T(t_1)$.

\vskip 5pt

Since $T(m(\mu)-2)$ is a non-empty subtree of $T(t_2)$ and $t_1 >
m(\mu)-2$, it follows that $T'(t_1) \cap T(t_2) \supseteq
T(m(\mu)-2) \neq \emptyset$, provided that $\phi^{\prime}$ is
sufficiently close to $\phi$. Since $T'(t_1) \cap T(t_2) \subseteq
T(t_1)$ and $T(t_1)$ is contained in the interior of $T(t_2)$, we
know that $T(t_2)$ contains a connected component of $T'(t_1)$.
Since $T'(t_1)$ is connected, we must have $T'(t_1) \subseteq
T(t_2)$. Therefore $T'(t_1)$ is finite, and so $\phi^{\prime} \in
({\bf \Phi}_\mu)_Q$. This proves Theorem \ref{thm:B3.16}. \square

\vskip 5pt

\noindent {\bf Proof of Theorem \ref{thm:proper}}.\,  This is
equivalent to showing that $\PSLtwoZ$ acts properly
discontinuously on $({\bf \Phi}_{\mu})_Q$, that is, for any
compact subset $K$ of $({\bf \Phi}_{\mu})_Q$, the set $\{H \in
\PSLtwoZ \mid HK\cap K \neq \emptyset \}$ is finite. Suppose not.
Then there exists a sequence of distinct $H_i \in \PSLtwoZ$ and
$\phi_i \in K$ such that $H_i(\phi_i) \in K$. Passing to a
subsequence, by the compactness of $K$, we may assume that $\phi_i
\rightarrow \phi \in K$, $H_i \rightarrow \infty$, and
$H_i(\phi_i) \rightarrow \phi' \in K$. Now, as in the proof of
Theorem \ref{thm:B3.16}, we have the tree $T_{\phi}(t_1)$ of
$\phi$ is finite for some $t_1>0$, and that $T_{\phi_i}(t_1)$ is
contained in the finite tree $T_{\phi}(t_2)$, for some $t_2>t_1$
and for all $i$ sufficiently large. This implies that the same
constant $\kappa$ can be used in the lower Fibonacci bound for all
$\phi_i$ for $i$ sufficiently large, and hence $H_i(\phi_i)
\rightarrow \infty$ as $i \rightarrow \infty$ (where we use the
identification of ${\bf \Phi}_{\mu}$ with the character variety
$\{(x,y,z) \in {\mathbf C}^3 \mid x^2+y^2+z^2-xyz=\mu\}$ to make
sense of $\phi \rightarrow \infty$). This contradicts $H_i(\phi_i)
\rightarrow \phi' \in K$. \square

\vskip 5pt

 \noindent {\bf Proof of Proposition \ref{prop:accummulationpoint}}.\,
Again we prove the equivalent statement for the action of
$\PSLtwoZ$ on ${\bf \Phi}_{\mu}$. Suppose that $\phi \in {\bf
\Phi}_{\mu}$ does not satisfy the extended BQ-conditions, then
either $\phi(X) \in (-2,2)$ for some $X \in \Omega$ or
$\Omega_{\phi}(2)$ is infinite. In the first case, since the
values of $\phi$ on the neighbors of $X$ are bounded by Lemma
\ref{lem:B3.3}(2), we can find an infinite sequence $H_i$ of
parabolic elements in $\PSLtwoZ$ all fixing $X$ such that
$H_i(\phi)$ converges to $\phi_0 \in {\bf \Phi}_{\mu}$ (here it is
possible that $H_i(\phi)=\phi_0$ for all $i$). In the latter case,
there are infinitely many edges $e \in E(\Sigma)$ such that
$\Omega^0(e) \subset \Omega_{\phi}(2)$. Passing to a subsequence
if necessary, this implies that there is an infinite sequence of
vertices $v_i \in V(\Sigma)$ such that if $X_i,Y_i$ and $Z_i$ are
the regions in $\Omega$ meeting at $v_i$, then $(x_i,y_i,z_i)$
converges to $(x_0,y_0,z_0)$ satisfying
(\ref{eqn:x^2+y^2+z^2=xyz}), hence there exists an infinite
sequence of distinct $H_i$ such that again $H_i(\phi)$ converges.
\square

\vskip 10pt

Note that for $\phi$ satisfying the extended BQ-conditions, for
any sequence of distinct elements $H_i$ in $\PSLtwoZ$, $H_i(\phi)
\rightarrow \infty$.

\vskip 20pt
\section{{\bf Proof of Theorem
\ref{thm:mcshane mu-markoff}}}\label{s:proof of theorem mu}
\vskip 10pt

In this section we prove Theorem \ref{thm:mcshane mu-markoff} by
defining the $\phi$-weight $\psi(\vec e)$ for each directed edge
$\vec e \in \vec E(\Sigma)$ and then following through the proof
of Theorem 2 in \cite{bowditch1998plms}. The geometric meaning of
$\psi(\vec e)$ is given in Appendix A, the reader interested in
the geometric motivation for the definition of $\psi$ is advised
to look there first.

\vskip 10pt

\noindent {\bf The functions $l$ and $h$ again.}\,\, For $x \in
{\mathbf C}$, recall that $l(x)/2 \in {\mathbf C}/2 \pi i {\mathbf
Z}$ is defined by $l(x)/2=\cosh^{-1}(x/2)$. In particular, $\Re
\big(l(x)/2\big) \ge 0$.

\vskip 10pt

\begin{lem}\label{lem:e^-l(x)/2} If $x \notin [-2,2]$ then
\begin{eqnarray}\label{eqn:e^-(x)l/2}
e^{-l(x)/2}=xh(x).
\end{eqnarray}
\end{lem}

\vskip 5pt

\begin{pf} Note that $\Re \big(l(x)/2\big) > 0$ since $x \notin [-2,2]$.
Hence $|e^{l(x)/2}| > |e^{-l(x)/2}|$ and
\begin{eqnarray*}
|x^{-1}e^{l(x)/2}| > |x^{-1}e^{-l(x)/2}|.
\end{eqnarray*}
On the other hand, $\cosh\big(l(x)/2\big)=x/2$ implies that
\begin{eqnarray*}
e^{l(x)/2}+e^{-l(x)/2}=x.
\end{eqnarray*}
Hence
\begin{eqnarray*}
\bigg(\frac{e^{l(x)/2}}{x}\bigg)^2-\bigg(\frac{e^{l(x)/2}}{x}\bigg)+\frac{1}{x^2}=0
\end{eqnarray*}
and
\begin{eqnarray*}
\bigg(\frac{e^{-l(x)/2}}{x}\bigg)^2-\bigg(\frac{e^{-l(x)/2}}{x}\bigg)+\frac{1}{x^2}=0.
\end{eqnarray*}
Then $|x^{-1}e^{-l(x)/2}| < |x^{-1}e^{l(x)/2}|$ implies that
\begin{eqnarray*}
x^{-1}e^{-l(x)/2}=\frac{1}{2}\bigg( 1-
\sqrt{1-\frac{4}{x^2}}\bigg)=h(x).
\end{eqnarray*}
This proves Lemma \ref{lem:e^-l(x)/2}.
\end{pf}

\vskip 5pt

\noindent {\bf The function $\Psi$.}\,\,  We define a function
$$\Psi: {\mathbf C}^3 \rightarrow {\mathbf C}$$ as
follows. Given $x,y,z \in {\mathbf C}$, set
$$\mu=x^2+y^2+z^2-xyz \,\,\,\,\,\,\,\,\, {\rm and} \,\,\,\,\,\,\,\,\,
\nu=\cosh^{-1}(1-\mu/2).$$ Then \,$\Psi(x,y,z) \in {\mathbf C}$\,
is defined by
\begin{eqnarray}\label{Eqn:Psi(x,y,z)=log}
\Psi(x,y,z)= \log
\frac{xy+(e^{\nu}-1)z}{(x^2-\mu)^{1/2}(y^2-\mu)^{1/2}},
\end{eqnarray}
or equivalently, by the following two equations:
\begin{eqnarray}
\cosh \Psi(x,y,z)=
\frac{xy-(\mu/2)z}{(x^2-\mu)^{1/2}(y^2-\mu)^{1/2}},
\end{eqnarray}
\begin{eqnarray}
\sinh \Psi(x,y,z)=
\frac{(\sinh\nu)\,z}{(x^2-\mu)^{1/2}(y^2-\mu)^{1/2}}.
\end{eqnarray}

\vskip 10pt

\noindent {\it Remarks.}\,\,
\begin{itemize}

\item[(i)] {\bf Notational convention.}\,\, In this paper, for a
complex number $u \in \mathbf C$, we use the notation $u^{1/2}$
(as opposed to $\sqrt{u}$\,) to mean either (once and for all)
choice of one of the two square roots of $u$.

\item[(ii)] Note that if $x^2, y^2 \neq \mu$ then $xy+(e^{\nu}-1)z=[xy-(\mu/2)z] + [(\sinh\nu)z] \neq
0$ since
\begin{eqnarray}
[xy-(\mu/2)z]^2=(x^2-\mu)(y^2-\mu)+[(\sinh\nu)z]^2.
\end{eqnarray}
Hence $\Psi(x,y,z)$ is defined if $x^2, y^2 \neq \mu$, or
equivalently, $y^2+z^2 \neq xyz$ and $x^2+z^2 \neq xyz$.

\item[(iii)] Note that the value $\Psi(x,y,z)$ depends on the
choices of the square roots in its expression, thus it is only
well-defined modulo $\pi i$ without specifying the choices of the
square roots. However, in Proposition \ref{prop:properties of Psi}
below, the appropriate sums there do not depend on the choices of
square roots and hence are well-defined modulo $2 \pi i$ with
arbitrary (once and for all) choices of the square roots.

\item[(iv)] It can be checked that
\begin{eqnarray}\label{eqn:Psi'=z/xy}
\frac{\partial}{\partial
\nu}\bigg|_{\nu=0}2\,\Psi(x,y,z)=\frac{z}{xy}.
\end{eqnarray}

\end{itemize}

\vskip 10pt

\noindent {\bf Some properties of the function $\Psi$.}\,\, The
function $\Psi$ defined above has the following nice properties.

\vskip 5pt

\begin{prop}\label{prop:properties of Psi} Let
$\mu=x^2+y^2+z^2-xyz$, $\nu=\cosh^{-1}(1-\mu/2)$ and suppose $\mu
\neq 0,4$ in each of the following cases.
\begin{itemize}

\item[{\rm (i)}] For $x,y,z \in {\mathbf C}$ with $x^2,y^2, z^2 \neq
\mu$,
\begin{eqnarray}\label{eqn:Psi+Psi+Psi=nu}
\Psi(y,z,x) + \Psi(z,x,y) + \Psi(x,y,z) = \nu \mod 2 \pi i.
\end{eqnarray}

\item[{\rm (ii)}] For $x,y,z,w \in {\mathbf C}$ with $z+w=xy$ and
$x^2,y^2 \neq \mu$,
\begin{eqnarray}\label{eqn:Psi+Psi=nu}
\Psi(x,y,z) + \Psi(x,y,w) = \nu \mod 2 \pi i.
\end{eqnarray}

\item[{\rm (iii)}] For $x,y \in {\mathbf C} \backslash \{0\}$, and
\begin{eqnarray}\label{eqn:z=}
z=\frac{xy}{2}\bigg( 1-\sqrt{1-4\left(\frac{1}{x^2}+
\frac{1}{y^2}-\frac{\mu}{x^2y^2} \right)}\,\bigg),
\end{eqnarray}
we have
\begin{eqnarray}
\lim_{y \rightarrow \infty} 2\Psi(x, y, z) = {\mathfrak h}(x) \mod
2 \pi i.
\end{eqnarray}

\item[{\rm (iv)}] If $\phi \in ({\bf \Phi_\mu})_{Q}$ then there is a constant
$C=C(\phi) > 0$ such that
\begin{eqnarray}\label{eqn:estimate}
|2\Psi(x,y,z)-{\mathfrak h}(x)| \le C |y|^{-2},
\end{eqnarray}
for all $x=\phi(X), X \in \Omega$; $|y|$ sufficiently large; and
$z$  given by {\rm(\ref{eqn:z=})}.
\end{itemize}
\end{prop}

\begin{pf} One can prove (i)--(iii) by direct calculations. For
details, see \S 2.6 of \cite{zhang2004thesis}.

\vskip 5pt

(iv) The estimate (\ref{eqn:estimate}) follows from the fact that
\vskip 5pt \centerline{$|\log(1+u)\,|\le 2\,|u|$ \hskip 10pt for
$u \in {\mathbf C}$ \, with $|u| \le 1/2$} \vskip 5pt \noindent
and the following calculations:
\begin{eqnarray*}
& &2\Psi(x,y,z)-{\mathfrak h}(x)\\
&=&\log\frac{x^2y^2[1+(e^\nu-1)(z/xy)]^2}{(x^2-\mu)(y^2-\mu)}-\log\frac{x^2[1+(e^\nu-1)h(x)]^2}{x^2-\mu}\\
&=&\log \bigg(\frac{y^2}{y^2-\mu}\bigg[ \frac{1+(e^\nu-1)(z/xy)}{1+(e^\nu-1)h(x)} \bigg]^2 \bigg)\\
&=&\log\bigg( 1+\frac{\mu}{y^2-\mu} \bigg)+2\log\bigg( 1+
\frac{(e^\nu-1)[z/xy - h(x)]}{1+(e^\nu-1)h(x)} \bigg),
\end{eqnarray*}
(note that the above are true equalities in $\mathbf C$ without
modulo $2\pi i$: the first one is because of (\ref{eqn:z=}) and
the other two because each expression after the log symbol is
sufficiently close to $1$ provided that $|y|$ is sufficiently
large)
\begin{eqnarray*}
\big|z/xy - h(x)\big|
&=& \big|\sqrt{1-4/x^2}-\sqrt{1-4/x^2-4/y^2+4\mu/x^2y^2}\,\big|\, \big/2\\
&=& \frac{\big|(4/y^2)(1-\mu/x^2)\big|}{2\big|\sqrt{1-4/x^2}+\sqrt{1-4/x^2-4/y^2+4\mu/x^2y^2}\,\big|}\\
&\le& 2|y|^{-2}\big|1-\mu/x^2\big|\,\big/\sqrt{|1-4/x^2|},
\end{eqnarray*}
since by our convention the square roots here all have nonnegative
real parts, and for $\phi \in ({\bf \Phi_\mu})_{Q}$, the image
$\phi(\Omega)$ is discrete.
\end{pf}

\vskip 5pt

{\bf Proof of Theorem \ref{thm:mcshane mu-markoff}.} \,\, The
proof now follows essentially the same line as that of Theorem 3,
\cite{bowditch1998plms}, with $\psi(\vec e) = z/xy$ there now
replaced by $\Psi(x,y,z) \in {\mathbf C}$. We write it out with
details as follows, taking care of multiples of $2\pi i$. Note
that this difficulty does not occur in \cite{bowditch1998plms}.
\vskip 3pt

Suppose $\phi \in ({\bf \Phi}_{\mu})_{Q}$ where $\mu \neq 0,4$. We
write ${\mathfrak h} = {\mathfrak h}_{\tau}$ where $\tau=\mu-2$.
Then $\log^{+}|\phi|$ has Fibonacci growth and the sum in
(\ref{eqn:mcshane mu-markoff}) converges absolutely since
${\mathfrak h}(x) = O(|x|^{-2})$ as $|x| \rightarrow \infty$ and
since $\sum_{X\in\Omega}|\phi(X)|^{-2}$ converges absolutely by
Corollary \ref{cor:sum phi^-2 convergence}.  \vskip 3pt

On the other hand, since ${\Omega_\phi}(2)$ is finite, as in the
proof of Theorem \ref{thm:B2}, there is a finite subtree $T$ of
$\Sigma$ with the property that for each edge not in $T$, its
arrow assigned by $\phi$ points towards $T$.  \vskip 3pt


Let $C_n, n \ge 0$ be the set of directed edges $\vec e
=\alpha_\phi(e)$ at a distance $n$ away from $T$ (that is, $e
\notin T$, and the minimal arc in $\Sigma$ from the head of $\vec
e$ to $T$ consists of $n$ edges). It is easy to see that for each
$n \ge 0$, $C_n$ is a circular set.  \vskip 3pt

Now let $\Omega_n, n \ge 0$ be the set of regions $X$ such that $X
\in \Omega^0(e)$ for some edge $\vec e \in C_n$. (Actually, each
such region $X$ intersects $C_n$ in exactly two edges.) Note that
$\Omega_n$ forms a nested sequence with respect to inclusion, that
is, $\Omega_n \subset \Omega_{n+1}$ for all $n \ge 0$ and
furthermore, $\Omega=\bigcup_{n=0}^{\infty}\Omega_n$. On the other
hand, all the $C_n, n \ge 0$ are disjoint. Note also that
$$\sum_{X\in\Omega} {\mathfrak h}(\phi(X)) = \lim_{n \rightarrow
\infty} \sum_{X \in \Omega_n} {\mathfrak h}(\phi(X)).$$ \vskip 3pt

For each ${\vec e}=(X,Y;\rightarrow Z) \in C_n$, we have defined
$$
\psi(\vec e)=\Psi(x,y,z) %
=\log \frac{1+(e^{\nu}-1)(z/xy)}{\sqrt{1-\mu/x^2}\sqrt{1-\mu/y^2}} \in \mathbf C. %
$$
It is easy to see from (\ref{eqn:Psi+Psi+Psi=nu}) and
(\ref{eqn:Psi+Psi=nu}) that $\sum_{{\vec e} \in C_n}\psi(\vec
e)=\nu \mod 2 \pi i$. Furthermore, we have:

\vskip 5pt

{\bf Claim.}\,\, $\tilde{\nu}:=\lim_{n \rightarrow
\infty}\sum_{{\vec e} \in C_n}\psi(\vec e) \in \mathbf C$ exists.

\vskip 5pt

To prove the claim, it suffices to show that $\sum_{{\vec e} \in
C_n}\psi(\vec e)=\sum_{{\vec e} \in C_{n+1}}\psi(\vec e)$ for
sufficiently large $n$. Consider any ${\vec e}=(X,Y; W \rightarrow
Z) \in C_n$. Let ${\vec e_X}=(Y, W; \rightarrow X)$ and ${\vec
e_Y}=(X, W; \rightarrow Y)$. Then ${\vec e_X}, {\vec e_Y} \in C_{n+1}$. %
It is easy to see that when $n$ is sufficiently large, $|w|$ and
one of $|x|,|y|$, say $|x|$, are sufficiently large. It follows
from the definition of the edge weight $\psi(\vec e_Y)$ that
$|\psi(\vec e_Y)|$ is sufficiently small. Since $\psi(\vec
e)=\psi(\vec e_X)+\psi(\vec e_Y) \mod 2\pi i$ by Proposition
\ref{prop:properties of Psi}, we know that in fact $\psi(\vec
e)=\psi(\vec e_X)+\psi(\vec e_Y)$ without the modulo $2\pi i$
restriction, for each $\vec e \in C_n$. This proves the claim.

\vskip 8pt

Now for each ${\vec e}\in C_{n+1}$, we write ${\vec
e}=(X,Y;\rightarrow Z)$, where we assume $X \in \Omega_n$ and $Y
\in \Omega_{n+1}\backslash\Omega_n$. Here $z=\phi(Z)$ satisfies
(\ref{eqn:z=}) since $\vec e$ points towards $Z$. Note that as
$\vec e$ ranges over $C_{n+1}$, each $X \in \Omega_n$ gets counted
twice. Since $\phi \in ({\bf \Phi}_{\mu})_{Q}$, $|y|=|\phi(Y)|$
will be sufficiently large if $n$ is. Thus by (\ref{eqn:estimate})
we have, as $n \rightarrow \infty$,
\begin{eqnarray*}
2\,\Big|\! \sum_{X \in \Omega_n} {\mathfrak h}(\phi(X))-\tilde{\nu} \,\Big| %
&=& \Big| \sum_{X \in \Omega_n} 2\,{\mathfrak h}(x)-\sum_{{\vec e}\,\in C_{n+1}}2\,\psi(\vec e) \,\, \Big| \\
&=& \Big|\sum_{{\vec e}\,\in C_{n+1}}\big(\,{\mathfrak h}(x)-2\,\psi(\vec e)\,\big) \,\, \Big| \\
& &\hspace{-18pt} (\,{\rm where} \,\,\,\, {\vec e}=(X,Y;\rightarrow Z) \,\,\,{\rm with}\,\,\, %
X \in \Omega_n, \, Y \in \Omega_{n+1} \backslash \Omega_n\,) \\
&\le& \sum_{{\vec e}\,\in C_{n+1}} \big| \,\,{\mathfrak h}(x)-2\,\psi(\vec e) \, \big| \\
&\le& \sum_{Y \in \Omega_{n+1}\backslash\Omega_n} \mathrm{constant} \cdot |\phi(Y)|^{-2} %
\,\, \longrightarrow \,\, 0,
\end{eqnarray*}
since $\sum_{Y \in \Omega_{n+1} \backslash \,\Omega_n}|\phi(Y)|^{-2} \rightarrow 0$ %
(as $n \rightarrow \infty$) by the convergence of
$$\sum_{Y\in\Omega} |\phi(Y)|^{-2}=\lim_{n \rightarrow \infty}\sum_{Y \in \,\Omega_n}|\phi(Y)|^{-2}.$$ %
Hence
\begin{eqnarray}\label{eqn:sum=tildenu}
\sum_{X\in\Omega} {\mathfrak h}(\phi(X)) %
= \lim_{n \rightarrow \infty}\sum_{X \in \Omega_n} {\mathfrak h}(\phi(X))=\tilde{\nu}. %
\end{eqnarray}
This proves Theorem \ref{thm:mcshane mu-markoff}. \square

\vskip 10pt


Similar to the estimate (\ref{eqn:estimate}), we have (with the
same assumption there)
\begin{eqnarray}\label{eqn:estimate hat}
|\psi(\vec e)-\hat{\mathfrak h}(x)| \le C |y|^{-2},
\end{eqnarray}
Then Proposition \ref{prop:B3.13 half} can be similarly proved
using (\ref{eqn:estimate hat}).

\vskip 10pt

\noindent {\bf Proof of Proposition \ref{prop:B3.13 half}.}\, By
Proposition \ref{prop:B3.9}, we know that
$\sum_{X\in\Omega^{-}(\vec e)}$ ${\mathfrak h}(\phi(X))$ and
$\sum_{X\in\Omega^{-}(\vec e)}$ $2\,\hat{\mathfrak h}(\phi(X))$
both converge absolutely, and they have difference a multiple of
$2\pi i$ since\, ${\mathfrak h}=2\,\hat{\mathfrak h} \! \mod 2\pi
i$. Thus we only need to evaluate \,
$$\sum_{X\in\Omega^{0}(\vec e)}\hat{\mathfrak h}(\phi(X)) %
+\sum_{X\in\Omega^{-}(\vec e)}2\,\hat{\mathfrak h}(\phi(X)).$$ %
We  apply the same arguments as the proof of Theorem
\ref{thm:mcshane mu-markoff}. \vskip 2pt

Let $C'_n$ be the set of directed edges $\vec {e'}$ in
$\Sigma^{-}(\vec e)$ at a distance $n$ from $\vec e$ (that is, $e'
\subset \Sigma^{-}(\vec e)$ and the minimal arc in $\Sigma$ from
the head of $\vec {e'}$ to the tail of $\vec e$ has exactly $n$
edges). Then $C_n=C'_n \cup \{-\vec e\}$ is a circular set. By the
assumption on $\vec e$, we know that, when $n$ is sufficiently
large, $\vec {e'}=\alpha_\phi(e')$ for every $\vec {e'} \in C'_n$.
\vskip 3pt

Let $\Omega_n:=\Omega^{0-}_n(\vec e)$ be the set of regions $X$ in
$\Omega^{0-}(\vec e)$ such that $X \in \Omega^{0}(e')$ for some
$\vec {e'} \in C'_n$. Then $\Omega_n \subset \Omega_{n+1}$ and
$\Omega^{0-}(\vec e)=\bigcup_{n=0}^{n}\Omega_n$. Thus
$$
\sum_{X \in \Omega^{0-}(\vec e)}{\mathfrak h}(\phi(X)) %
=\lim_{n \rightarrow \infty}\sum_{X \in \Omega_n}{\mathfrak h}(\phi(X)). %
$$

By (\ref{eqn:Psi+Psi+Psi=nu}) and (\ref{eqn:Psi+Psi=nu}), we have
$\psi(\vec e)=\sum_{\vec {e'} \in C'_n}\psi(\vec {e'}) \! \mod
2\pi i$ \,for all $n \ge 0$. Similar to the claim in the proof of
Theorem \ref{thm:mcshane mu-markoff}, we can prove that

\vskip 5pt

{\bf Claim.}\,\, $\tilde{\psi}(\vec e):=\lim_{n \rightarrow \infty} %
\sum_{{\vec e} \in C_n}\psi(\vec e) \in \mathbf C$ exists. \square %

\vskip 5pt

Thus for $n$ sufficiently large we have from (\ref{eqn:estimate hat}) %
\begin{eqnarray*}
& & \Big|\,\tilde{\psi}(\vec e) %
-\sum_{X\in\Omega^{0}(\vec e)}\hat{\mathfrak h}(\phi(X))%
-\sum_{X\in\Omega^{-}(\vec e)}2\,\hat{\mathfrak h}(\phi(X))\,\Big| \\
&=&\Big|\sum_{\vec {e'} \in C'_{n+1}}\psi(\vec {e'}) %
-\sum_{X\in\Omega^{0}(\vec e)}\hat{\mathfrak h}(\phi(X)) %
-\sum_{X\in\Omega^{-}(\vec e)}2\,\hat{\mathfrak h}(\phi(X))\,\Big| \\
&=& \Big|\sum_{\vec {e'} \in C'_{n+1}}\big[\psi(\vec {e'})-\hat{\mathfrak h}(\phi(X))\big]\,\,\Big| \\
& & (\,{\rm where} \,\,\,\, \vec {e'}=(X,Y;\rightarrow Z) \,\,\, %
{\rm with} \,\,\, X \in \Omega_n \,\,\, {\rm and} \,\,\, Y \in \Omega_{n+1} \backslash \Omega_n\,) \\
&\le& \sum_{\vec {e'} \in C_{n+1}} \big| \,\psi(\vec {e'})-\hat{\mathfrak h}(x) \, \big| \\
&\le& \sum_{Y \in\, \Omega_{n+1} \backslash \, \Omega_n} \mathrm{constant} \cdot |\phi(Y)|^{-2} %
\,\, \longrightarrow \,\, 0,
\end{eqnarray*}
since $\sum_{Y \in \Omega_{n+1} \backslash \,\Omega_n}|\phi(Y)|^{-2} \rightarrow 0$ %
(as $n \rightarrow \infty$) by the convergence of
$$
\sum_{Y \in \Omega^{0-}(\vec e)} |\phi(Y)|^{-2} %
=\lim_{n \rightarrow \infty}\sum_{Y \in \,\Omega_n}|\phi(Y)|^{-2}. %
$$
Hence
\begin{eqnarray}\label{eqn:sum=tildepsi}
\tilde{\psi}(\vec e)
&=&\sum_{X\in\Omega^{0}(\vec e)}\hat{\mathfrak h}(\phi(X)) %
+\lim_{n \rightarrow \infty}\sum_{X \in \Omega_n}2\,\hat{\mathfrak h}(\phi(X)) \nonumber \\ %
&=&\sum_{X\in\Omega^{0}(\vec e)}\hat{\mathfrak h}(\phi(X)) %
+\sum_{X\in\Omega^{-}(\vec e)}2\,\hat{\mathfrak h}(\phi(X)). %
\end{eqnarray}
This proves Proposition \ref{prop:B3.13 half} since\,
$\tilde{\psi}(\vec e)={\psi}(\vec e) \,\,\, {\rm mod}\,\, 2\pi i$.
\square

\vskip 20pt
\section{{\bf Variations to once-punctured torus bundles}}\label{s:torus bundles} %
\vskip 10pt

In this section we generalize Bowditch's variations in
\cite{bowditch1997t} of McShane's identity.
%
Bowditch's variations in \cite{bowditch1997t} of McShane's
identity are for complete, finite volume hyperbolic 3-manifolds
fibering over the circle, with fiber the once-punctured torus.
These complete structures can be deformed to incomplete
structures, as shown by Thurston in \cite{thurston1978notes}, and
in certain cases, one can perform hyperbolic Dehn surgery to
obtain closed (complete) hyperbolic 3-manifolds. The main result
of this section is that a further variation of the Bowditch's
identities holds for these deformations satisfying the relative
BQ-conditions (to be defined later), and hence for the closed
hyperbolic 3-manifolds obtained by hyperbolic Dehn surgery on such
manifolds (Theorems \ref{thm:A'}, \ref{thm:B'} and Corollary
\ref{cor:Dehn}). Briefly, this is done as follows. Let $M$ be a
once-punctured torus bundle over the circle with non-trivial
monodromy $H$ which is an hyperbolic element in $\SLtwoZ$
identified with the mapping class group of the once-punctured
torus. Starting with a complete, finite volume hyperbolic
structure on $M$, and by considering the holonomy representation
restricted to the fundamental group of the fiber, Bowditch
\cite{bowditch1997t} constructed a Markoff map which is stabilized
by the cyclic subgroup $\langle H\rangle < {\rm SL}(2, \mathbf
Z)$. His variations of McShane's identity (Theorems \ref{thm:A}
and \ref{thm:B}) are obtained by studying the Markoff map modulo
the action of $\langle H\rangle$. If we deform the complete
hyperbolic structure on $M$ to an incomplete one and look again at
the holonomy representation restricted to the fiber, we get a
generalized Markoff map, with the same stabilizer group $\langle
H\rangle$. We then obtain further variations of the
McShane-Bowditch identities (Theorems \ref{thm:A'} and
\ref{thm:B'}) for these deformations, which also hold for almost
all of the closed hyperbolic 3-manifolds obtained by hyperbolic
Dehn surgery on these manifolds (Corollary \ref{cor:Dehn}).

\vskip 10pt
\subsection{Bowditch's settings for torus bundles}\label{s:settings} %


Let $M$ be an orientable 3-manifold which fibers over the circle,
with the fiber a once-punctured torus, $\mathbb T$. Suppose that
$M$ has a complete finite-volume hyperbolic structure.

Let $\mathscr S$ be the set of closed geodesics in $M$ which are
homotopic in $M$ to simple closed curves in the fiber.
Geometrically, to each $\sigma \in {\mathscr S}$, we associate its
{\bf complex length}, $l(\sigma) \in {\mathbf C} / 2\pi i {\mathbf
Z}$, where $\Re\, l(\sigma)$ is the (real) hyperbolic length of
$\sigma$ and $\Im\, l(\sigma)$ is the rotational component, i.e.
the angle through which a normal vector turns when parallel
transported once around the closed geodesic $\gamma$. In
particular, $e^{l(\sigma)} \in {\mathbf C}$ is well-defined. Then
Bowditch obtained

\begin{thm}\label{thm:A} {\rm (Theorem A in \cite{bowditch1997t})}
\begin{eqnarray*}
\sum_{\sigma \in {\mathscr S}} \frac{1}{1+e^{l(\sigma)}}=0,
\end{eqnarray*}
where the sum converges absolutely.
\end{thm}

The curves in $\mathscr S$ fall naturally into two classes,
${\mathscr S}_{\rm L}$ and ${\mathscr S}_{\rm R}$, as follows.
Recall that $\mathscr C$ is the set of free homotopy classes of
non-trivial non-peripheral simple closed curves on $\mathbb T$ and
can be thought of as the set of rational points in the projective
lamination space, $\mathscr P$, of $\mathbb T$, which in this case
is a circle. The mapping class group $\MCG$ of $\mathbb T$ acts on
$\mathscr P$ preserving the set $\mathscr C$. The monodromy $H \in
\MCG$ of $M$ generates an infinite cyclic subgroup, $\langle H
\rangle$, of the mapping class group. This subgroup, or $H$, has
two fixed points in $\mathscr P$, namely the stable and unstable
laminations, $\mu_{\rm s}$ and $\mu_{\rm u}$, of the monodromy.
These two points separate $\mathscr P$ into two open intervals.
Since $\mu_{\rm s}$ and $\mu_{\rm u}$ are irrational points, this
gives a natural partition of $\mathscr C$ into two subsets,
${\mathscr C}_{\rm L}$ and ${\mathscr C}_{\rm R}$, which in turn
partitions $\mathscr S$ into two subsets, ${\mathscr S}_{\rm L}$
and ${\mathscr S}_{\rm R}$.

If we restrict the sum appearing in Theorem \ref{thm:A} to one or
other of ${\mathscr S}_{\rm L}$ and ${\mathscr S}_{\rm R}$, we
will get the same answer up to change of signs. Bowditch showed
(see Theorem \ref{thm:B}) that this number turns out essentially
to be the modulus of the cusp of $M$.

More precisely, since $M$ has a single parabolic cusp, $M$ is
homeomorphic to the interior of a compact manifold $M \cup
\partial M$, with one toroidal boundary component, $\partial M$.
Then $\partial M$ carries a natural Euclidean structure,
well-defined up to similarity, which arises from identifying
$\partial M$ with a horocycle. If $M$ has positive monodromy, we
may represent $\partial M$ as the quotient of $\mathbf C$ with the
Euclidean metric by the lattice ${\mathbf Z} \oplus \lambda
{\mathbf Z}$, generated by the translations $\zeta \mapsto \zeta +
1$ and $\zeta \mapsto \zeta + \lambda$ corresponding to the
meridian and longitude respectively, where the {\bf longitude} is
defined as the intersection of two leaves, one from each of the
two foliations of $M$ determined by the stable and unstable
laminations. We call $\lambda = \lambda (M)$ the {\bf modulus} of
the cusp. We can suppose that $\Im [\lambda (M)] > 0$. If $M$ has
negative monodromy, we define $\lambda (M)$ as the modulus of the
cusp of the sister of $M$. In this setting Bowditch showed

\vskip 10pt

\begin{thm}\label{thm:B} {\rm (Theorem B in \cite{bowditch1997t})}
\begin{eqnarray*}
\sum_{\sigma \in {\mathscr S_{\rm L}}}
\frac{1}{1+e^{l(\sigma)}}=\pm \lambda (M),
\end{eqnarray*}
where the sign depends only on our conventions of orientation.
\end{thm}

\vskip 5pt

\noindent {\bf Natural ideal triangulations.}\,\, Recall the
setting of \S \ref{s:generalized Markoff maps}. Let us fix an
ordered pair of free generators $a,b$ of $\Gamma$ so that the
algebraic intersection number of the corresponding ordered pair of
oriented simple closed curves $C_a, C_b$ on $\mathbb T$ is equal
to $+1$. Note that the commutator $[a,b]=aba^{-1}b^{-1}$ is
peripheral.

The mapping class group $\MCG$ of $\mathbb T$ may be identified as
$\SLtwoZ$ which acts naturally on $\mathbb H^2, \Sigma, \Omega,
\mathscr P$ and $\mathscr C$. In each case, the kernel is given by
the elliptic involution, so the induced action of $\PSLtwoZ$ is
faithful. An element $H \in \SLtwoZ$ is {\it hyperbolic} if it has
two fixed points $\mu_{\rm u}$ and $\mu_{\rm s}$ in $\mathscr P$,
namely the stable and unstable laminations. The points $\mu_{\rm
u}$ and $\mu_{\rm s}$ are joined by a bi-infinite arc $\beta
\subseteq \Sigma$, which is translated by $H$ in the direction of
$\mu_{\rm s}$. The path $\beta$ can be described combinatorially
in terms of the ``right-left'' decomposition of the matrix $H$.
Note that some conjugate of $H$ can be written as a product of the
matrices
\begin{eqnarray*}
L=\bigg(\begin{matrix} 1 & 1 \\ 0 & 1 \end{matrix}\bigg)
\,\,\,{\rm and}\,\,\, R=\bigg(\begin{matrix} 1 & 0 \\ 1 & 1
\end{matrix}\bigg)
\end{eqnarray*}
which correspond to Dehn twists about $C_b$ and $C_a$
respectively. This decomposition is well-defined up to cyclic
reordering, and the sequence of $L$s and $R$s is the same as the
periodic sequence of left and right turns of $\beta$ in $\Sigma$.
For more details, see \cite{bowditch-maclachlan-reid1995ma}. Note
that $\beta$ partitions $\Omega$ into two subsets, $\Omega_{\rm
L}$ and $\Omega_{\rm R}$, which lie on, respectively, the left and
right of $\beta$. These correspond to the subsets, $\mathscr
C_{\rm L}$ and $\mathscr C_{\rm R}$ of $\mathscr C$ described
earlier.

\vskip 10pt

Now, if we take a homeomorphism, $\eta$, of $\mathbb T$
representing the mapping class $H$, we may form the mapping torus,
$M_H$, which is given by $(\mathbb T \times [0, 1])/ \sim$, where
$\sim$ identifies $(x, 1)$ with $(\eta(x), 0)$ for all $x \in
\mathbb T$. The manifold $M_H$ has a natural compactification by
adjoining a toroidal boundary, $\partial M_H$. The compactified
manifold, $M_H \cup \partial M_H$, as described in
\cite{floyd-hatcher1982ta}, has a natural ideal triangulation
arising from the left-right decomposition of $H$. This can briefly
be summarized as follows. Associated with each vertex of $\Sigma$
is an ideal triangulation of $\mathbb T$ (the edges of which are
dual to the simple closed curves corresponding to the three
complementary regions incident to the vertex). Moving along an
edge in the tree corresponds to performing a (dual) Whitehead
move. If we traverse a period of the path $\beta$, we get a
sequence of Whitehead moves which take us from a given
triangulation to its image under $\eta$. Each Whitehead move gives
rise to an ideal simplex in $\mathbb T \times [0, 1]$, and so,
after identifying $\mathbb T \times \{0\}$ with $\mathbb T \times
\{1\}$ via the relation $\sim$, we obtain an ideal triangulation
of $M_H$. This ideal triangulation gives us, in particular, a
triangulation of the boundary, $\partial M_H$. As in
\cite{floyd-hatcher1982ta}, we may describe the combinatorial
structure of this triangulation lifted to the universal cover,
${\mathbf R}^{2}$, of $\partial M_H$. First, we describe the case
of positive monodromy. To do this, consider the bi-infinite
sequence of vertices of $\Sigma$ lying along the arc $\beta
\subseteq \Sigma$. This sequence is dual to a sequence of ideal
triangles in our tessellation of ${\mathbb H}^{2}$. The union of
these triangles gives a bi-infinite strip invariant under the
transformation $H$. We transfer this strip homeomorphically to the
vertical strip $[0, 1] \times {\mathbf R}$ in ${\mathbf R}^{2}$,
so that the transformation $H$ is conjugated to the map $[(x, y)
\mapsto (x, y + 1)]$. We extend this to a triangulation of
${\mathbf R}^{2}$ by a process of repeated reflection in the pair
of vertical lines which form the boundary of this strip. 
The triangulation of $\partial M_H$ is given by the quotient by
the group generated by $[(x, y) \mapsto (x, y + 1)]$ and $[(x, y)
\mapsto (x + 4, y)]$. It is not hard to see that these
transformations describe, respectively, the longitude and meridian
of $\partial M_H$. Note that the triangulation is in fact
invariant under the map $[(x, y) \mapsto (x + 2, y)]$. Up to this
symmetry, there are two ``vertical'' lines in the triangulation.
The bi-infinite sequence of vertices along one of these lines
corresponds to sequence of regions of $\Omega$ which meet $\beta$
and all lie either in $\Omega_{\rm L}$, or in $\Omega_{\rm R}$.
Two vertices are joined by an edge in this triangulation if and
only if the corresponding regions are adjacent. Thus the
``vertical'' edges correspond to regions meeting on the same side
of $\beta$, whereas all the other edges correspond to regions
meeting on opposite sides of $\beta$. The picture where the
monodromy is negative is similar, except that in this case,
$\partial M_H$ is given as a quotient of ${\mathbf R}^{2}$ by the
group generated by $[(x, y) \mapsto (x + 2, y + 1)]$ and $[(x, y)
\mapsto (x + 4, y)]$. The ``longitude'' might be thought of as a
``half of'' the curve given by $[(x,y) \mapsto (x, y+2)]$. Now, it
follows from the work of Thurston \cite{thurston19??am} that $M =
M_H$ admits a complete finite-volume hyperbolic structure. See
also \cite{otal1994ppt} for an alternative proof and exposition.
This structure is unique by Mostow rigidity. Note that $\partial
M_H$ carries a Euclidean structure, well defined up to similarity,
obtained for example by identifying it with a horocycle in $M_H$.
In this hyperbolic structure, we may realize each tetrahedron in
our ideal triangulation of $M$ as a hyperbolic ideal tetrahedron.
In this way we get a ``hyperbolic ideal triangulation'' of $M$.
This gives rise to a Euclidean realization of the combinatorial
triangulation of $\partial M_H$. Note that Parker
\cite{parker2003lmsln299} has prove that this hyperbolic ideal
triangulation of $M$ obtained above is indeed positively oriented,
hence is a genuine hyperbolic ideal triangulation. Thus we have a
resulting geometric triangulation of $\partial M_H$ which will be
used to prove the variations of McShane's identity.

\vskip 10pt

Regarding $\mathbb T$ as a fibre of $M$, we get an identification
of $\Gamma = \pi_1(\mathbb T)$ as a normal subgroup of $\pi_1(M)$.
In fact, $\pi_1(M)$ is an HNN extension of $\Gamma$ with stable
letter $t$ so that $tgt^{-1}=H_{*}(g)$ for all $g \in \Gamma$,
where $H_{*}$ is the automorphism of $\Gamma$ induced by the
monodromy $H$. We also get an identification of $\mathscr S$ with
the quotient, $\Omega/\langle H \rangle$, of $\Omega$ under the
cyclic group, $\langle H \rangle$, generated by $H$. Clearly, $H$
respects the partition of $\Omega$ as $\Omega_{\rm L} \sqcup
\Omega_{\rm R}$, and we may identify ${\mathscr S}_{\rm L}$ with
$\Omega_{\rm L}/\langle H \rangle$ and ${\mathscr S}_{\rm R}$ with
$\Omega_{\rm R}/\langle H \rangle$. The hyperbolic structure on
$M$ may be described by a representation, ${\hat \rho}: {\pi_1}(M)
\rightarrow \PSLtwoC$. It follows from \cite{culler1986advm} that
$\hat \rho$ lifts to a representation $\rho: {\pi_1}(M)
\rightarrow \SLtwoC$.

Restricting our attention to the fibre subgroup $\Gamma <
{\pi_1}(M)$, we define a Markoff map $\phi: \Omega \rightarrow
{\mathbf C}$ by $\phi(X) = {\rm tr}\,\rho(g)$, where $g \in
\Gamma$ represents the simple closed curve on $\mathbb T$
corresponding to the region $X \in \Omega$. Clearly, $\phi$ is
invariant under the $\langle H \rangle$-action, and so gives rise
to a well-defined map $\Omega/\langle H \rangle \rightarrow
{\mathbf C}$ which we also denote by $\phi$. We write $[X]$ for
the orbit of $X$ under $\langle H \rangle$. If $\sigma \in
{\mathscr S}$ corresponds to $[X] \in \Omega/\langle H \rangle$,
then the complex length, $l(\sigma)$, of $\sigma$ is determined by
the formula $l(\sigma)=l(\phi([X]))=2\cosh^{-1}(\phi([X])/2)$.
Thus $h(\phi([X])) = 1/(1 + e^{l(\sigma)})$, where $h: {\mathbf C}
\backslash [-2,2] \rightarrow {\mathbf C}$ is defined by
$h(x)=\big(1-\sqrt{1-4/x^2}\, \big)\big/2$ as in \S
\ref{s:generalized Markoff maps}. Here we take the square root
with positive real part, corresponding to the fact that $\Re \,
l(\sigma)> 0$. In these terms Bowditch's Theorems \ref{thm:A} and
\ref{thm:B} can be respectively expressed as the identities:
\begin{eqnarray}\label{eqn:restate theorem A}
\textstyle{\sum_{[X] \in \Omega/\langle H \rangle}}\,
h(\phi([X]))=0,
\end{eqnarray}
\begin{eqnarray}\label{eqn:restate theorem B}
\textstyle{\sum_{[X]\in\Omega_{\rm L}/\langle H \rangle}}\,
h(\phi([X]))=\lambda(\partial M).
\end{eqnarray}

\vskip 5pt
\subsection{Incomplete hyperbolic torus bundles}\label{s:incomplete} %

In this subsection we consider  once-punctured torus bundles over
the circle, but now with incomplete hyperbolic structures.

Let $M$ be a once-punctured torus bundle over the circle such that
the holonomy $H \in \MCG$ is pseudo-Anosov. Then as described in
the previous subsection, $M$ can be given a complete, finite
volume hyperbolic structure, and in fact, $M$ can be decomposed
into a collection of ideal hyperbolic tetrahedra. Hence we may
obtain $M$ by gluing a collection of ideal hyperbolic tetrahedra
with suitable edge invariants. Now we consider an incomplete
hyperbolic structure on $M$ obtained by changing the edge
invariants while keeping the consistency conditions as described
by Thurston \cite{thurston1978notes}. We can also regard this
incomplete structure as being obtained by deforming the complete
hyperbolic structure slightly in the representation space of
$\pi_1(M)$.

The developing image of $M$, with this incomplete hyperbolic
structure, in ${\mathbb H}^3$ misses some lines in $\HHH$. We may
assume the $z$-axis $[0,\infty]$ is among them. Consider, for
$\epsilon >0$ sufficiently small, an $\epsilon$-neighborhood $N$
of $[0,\infty]$ in $\HHH$, i.e. a solid cone around the $z$-axis.
Its boundary $\partial N$ is the $\epsilon$-equidistant surface
with center the $z$-axis, hence has a similarity structure.
Actually, we can identify this similarity surface $\partial N$
with $\Cnozero$, where the identification is given by orthogonal
projection in $\HHH$ from $\Cnozero$ onto the $z$-axis via lines
normal to it. Note that $N \backslash [0,\infty]$ projects onto a
neighborhood of $\partial M$ in $M$, hence $\partial N$ projects
onto $\partial M$. In this way $\partial M$ gets a similarity
structure from the (incomplete) hyperbolic structure of $M$.

Recall the combinatorial triangulation of $M$ described in \S
\ref{s:settings}. Now this triangulation of $M$ is realized
similarly by a ``hyperbolic ideal triangulation'' and $\partial M$
gets an induced conformal triangulation. This triangulation can be
thought of as being deformed from the triangulation in the
complete case when we deform the complete structure on $M$ into an
incomplete one.

This incomplete hyperbolic structure on $M$ is given by a
representation ${\hat \rho}: {\pi_1}(M) \rightarrow \PSLtwoC$
which lifts to a representation $\rho: {\pi_1}(M) \rightarrow
\SLtwoC$. Restricted to $\Gamma = \pi_1(\mathbb T) \triangleleft
\pi_1(M)$, we obtain a $\mu$-representation $\rho: \Gamma
\rightarrow \SLtwoC$ for some $\mu \in {\mathbf C}$ with $\mu \neq
4$. Let $\nu=\cosh^{-1}(1-\mu/2)$ as in \S \ref{s:generalized
Markoff maps}.

Lifting the similarity structure on $\partial M$ to $\Cnozero$, we
get a triangulation of $\Cnozero$ which is invariant under the
transformations $\zeta \mapsto e^{\nu} \zeta$ and $\zeta \mapsto
e^{\lambda} \zeta$ if $M$ has positive monodromy. These two
transformations give respectively the meridian and longitude of
$\partial M$, similar to the complete case in \S \ref{s:torus
bundles}. If $M$ has negative monodromy then the two
transformations are given by $\zeta \mapsto e^{\nu} \zeta$ and
$\zeta \mapsto e^{\nu+\lambda} \zeta$.

Note that $\nu \in \CmodTwoPiIZ$ is  half of the complex length of
the $\rho$-image of a peripheral simple closed curved on the
punctured torus $\mathbb T$, the fiber of $M$.

We assume the notation introduced in the beginning of this
section. In particular, for each closed geodesic $\sigma \in
{\mathscr S}$, its complex length \,$l(\sigma) \in \CmodTwoPiIZ$\,
is defined now under the given incomplete hyperbolic structure.

\vskip 10pt

We say that $M$ (with the incomplete hyperbolic structure)  has
{\it discrete length spectrum} on the puncture-torus fiber
$\mathbb T$ if for each $\kappa >0$ there are only finitely many
free homotopy classes of simple closed curves on $\mathbb T$ such
that their geodesic realizations in $M$ have length $\le \kappa$.
We will see later that  $M$ has this property if it is obtained
from the  unique complete finite volume hyperbolic structure by a
sufficiently small deformation (Theorem \ref{thm:openforbundle}).
Then we have the following results corresponding to Bowditch's
Theorems \ref{thm:A} and \ref{thm:B}.

\begin{thm}\label{thm:A'} Suppose $M$ {\rm(}with incomplete hyperbolic
structure{\rm)} has discrete length spectrum on its torus fiber.
Then
\begin{eqnarray}\label{eqn:A'}
\sum_{\sigma \in {\mathscr S}}
\log\bigg(\frac{e^{\,\nu}+e^{\,l(\sigma)}}{e^{-\nu}+e^{\,l(\sigma)}}\bigg)=0
\mod 2 \pi i,
\end{eqnarray}
where the sum converges absolutely.
\end{thm}

\begin{thm}\label{thm:B'} Suppose $M$ {\rm(}with incomplete hyperbolic
structure{\rm)} has discrete length spectrum on its torus fiber.
Then
\begin{eqnarray}\label{eqn:B'}
\sum_{\sigma \in {\mathscr S_{\rm L}}}
\log\bigg(\frac{e^{\,\nu}+e^{\,l(\sigma)}}{e^{-\nu}+e^{\,l(\sigma)}}\bigg)=\pm
\lambda \mod 2 \pi i,
\end{eqnarray}
where the sign depends only on our conventions of orientation and
where $\lambda$ is the complex length of the chosen longitude.
\end{thm}

As a corollary we have

\begin{cor}\label{cor:Dehn} Let $M$ be a once-punctured torus bundle over the circle with
pseudo-Anosov monodromy and let $\overline{M}(p/q)$ be the closed
3-manifold obtained from $M$ by performing Dehn surgery on
$\partial M$ with surgery slope $p/q \in \mathbf Q$. Then, except
for a finite number of surgery slopes, $\overline{M}(p/q)$ has a
hyperbolic structure and, moreover, the following identity holds:
\begin{eqnarray}\label{eqn:Dehn}
\sum_{\sigma \in {\mathscr C_{\rm L}}}
\log\bigg(\frac{e^{\,\nu}+e^{\,l(\sigma)}}{e^{-\nu}+e^{\,l(\sigma)}}\bigg)=\pm
\lambda \mod 2 \pi i,
\end{eqnarray}
where the meanings of $\lambda$, $\nu$ and ${\mathscr C}_{\rm L}$
etc are defined as before for the {\rm(}incomplete{\rm)}
hyperbolic structure on $M$ induced from that of $\overline{M}(p/q)$ %
and the sign depends only on our conventions of orientation.
\end{cor}


\noindent {\it Remark.} We may weaken the condition of discrete
length spectrum in the above results to the condition that the set
$\{\gamma \in {\mathscr S} \mid {\mathfrak R}\,l(\gamma)\le \log
(3+2\sqrt 2)\}$ is finite,  by a simple calculation.

\vskip 5pt
\subsection{Periodic generalized Markoff maps}\label{s:stabilazed} %

In this subsection we consider generalized Markoff maps which are
invariant under a hyperbolic element of the mapping class group
$\MCG$ of $\mathbb T$.

Consider the action of the mapping class group of $\mathbb T$,
$\mathcal{MCG} \cong {\rm SL}(2,{\mathbf Z})$, and its induced
action  on $\Gamma$ and on the space ${\mathcal X}_{\tau}$ of
$\tau$-representations $\rho: \Gamma \rightarrow {\rm
SL}(2,{\mathbf C})$ modulo conjugation. Any $H \in \mathcal{MCG}$
induces an automorphism $H_{\ast}$ of $\Gamma$, and $H$ acts on
${\mathcal X}_{\tau}$ by
$$H(\rho)(g)=\rho(H_{\ast}(g)),$$
for any $\rho \in {\mathcal X}_{\tau}$ and $g \in \Gamma$.
Bowditch \cite{bowditch1997t} studied representations $\rho \in
{\mathcal X}_{\rm tp}={\mathcal X}_{-2}$ stabilized by a cyclic
subgroup $\langle H \rangle < \mathcal{MCG}$ generated by a
hyperbolic element and proved a variation of the McShane's
identity for such representations. This result can be generalized
for $\tau$-representations as follows, and is equivalent to
Theorem \ref{thm:A'}, with slightly weaker assumptions. The case
$\tau=-2$ with $\mathfrak h$ replaced by $h$ is Bowditch's
variation (Theorem A in \cite{bowditch1997t}).

\begin{thm}\label{thm:rep stabilized}
Suppose that a $\tau$-representation $\rho: \Gamma \rightarrow
\SLtwoC$, where $\tau \neq 2$, is stabilized by a hyperbolic
element $H \in {\rm SL}(2,{\mathbf Z}) \cong \mathcal{MCG}$ and
$\rho$ satisfies the BQ-conditions on $\hat \Omega/\langle
H_{\ast}\rangle$, that is, \begin{itemize}

\item[{\rm (i)}] \,\, ${\rm tr}\rho(g) \notin [-2,2]$ for all classes $[g] \in \hat
\Omega/\langle H_{\ast}\rangle$, and

\item[{\rm (ii)}] $|{\rm tr}\rho(g)| \le 2$ for only finitely many classes
$[g] \in \hat \Omega/\langle H_{\ast}\rangle$.
\end{itemize} Then we have
\begin{eqnarray}\label{eqn:rep stabilized}
\sum_{[g] \in \hat \Omega/\langle H_{\ast}\rangle} {\mathfrak
h}\big( {\rm tr}\rho(g)\big)=0,
\end{eqnarray}
where the function $\mathfrak h=\mathfrak h_\tau$ is defined by
{\rm(\ref{eqn:frak h(x)=log I})} and the sum converges absolutely.
\end{thm}

\vskip 5pt

We prove Theorem \ref{thm:rep stabilized} by reformulating it in
terms of generalized Markoff maps as follows.

\vskip 10pt

\begin{thm}\label{thm:stabilized}
Suppose that $\phi \in {\bf \Phi}_{\mu}$ {\rm(}$\mu \neq 4${\rm)}
is invariant under the action of a hyperbolic element $H \in {\rm
SL}(2,{\mathbf Z}) \cong \mathcal{MCG}$ and $\phi$ satisfies the
BQ-conditions on $\Omega/\langle H\rangle$, that is, {\rm (i)}
${\phi}^{-1}([-2,2])=\emptyset$, and {\rm (ii)} $|\phi([X])| \le
2$ for only finitely many classes $[X] \in \Omega/\langle
H\rangle$. Then
\begin{eqnarray}\label{eqn:stabilized}
\sum_{[X] \in \Omega/\langle H\rangle} {\mathfrak h}\big(
\phi([X])\big)=0 \mod 2 \pi i,
\end{eqnarray}
where the sum converges absolutely and the function ${\mathfrak
h}={\mathfrak h}_\tau$ is defined as in {\rm(\ref{eqn:frak
h(x)=log I})} with $\tau=\mu-2$ and $\nu=\cosh^{-1}(-\tau/2)$.
\end{thm}

\vskip 10pt

\begin{pf} The proof is essentially the same as that given by Bowditch in
\cite{bowditch1997t}, except here we use the generalized function
$\mathfrak h$ and the generalized $\phi$-weight $\psi(\vec e)$
(this follows along the same lines as the proof for Theorem
\ref{thm:mcshane mu-markoff}).
\end{pf}

We also have the following openness result for $\mu$-Markoff maps
invariant under a fixed hyperbolic element $H \in \SLtwoZ$ and
satisfying the relative BQ-conditions.

\begin{thm}\label{thm:openforbundle} Let ${\bf \Phi}^H$ be the set of
generalized Markoff maps invariant under a fixed hyperbolic
element $H \in \SLtwoZ$ and ${\bf \Phi}_Q^H$ be the subset
satisfying the BQ-conditions on $\Omega/\langle H\rangle$. Then
${\bf \Phi}_Q^H$ is open in ${\bf \Phi}^H$.

\end{thm}

\begin{pf}
This can be proved similarly as Theorem \ref{thm:B3.16}. For $\phi
\in {\bf \Phi}_Q^H$  we construct the tree $T_{\phi}(t_1)$ which
is invariant under $H$. $T_{\phi}(t_1)$ is not finite but
$T_{\phi}(t_1)/\langle H\rangle $ is a finite graph in the graph
$\Sigma/\langle H\rangle $. Note that $T_{\phi}(t_1)/\langle
H\rangle $ contains the spine of $\Sigma/\langle H\rangle $ which
is the projection of the invariant axis of $H$ in $\Sigma$. Now if
$\phi'\in {\bf \Phi}_Q^H$ is sufficiently close to $\phi$, then
$T_{\phi'}(t_1)/\langle H\rangle \subseteq T_{\phi}(t_2)/\langle
H\rangle$ for some given $t_2>t_1$, hence, $\phi' \in {\bf
\Phi}_Q^H$.
\end{pf}

The above result tells us that $M$ (with the incomplete hyperbolic
structure) would continue to have discrete length spectrum on its
torus fiber for sufficiently small deformations  from the complete
structure.

\vskip 10pt

We may also reformulate Theorem \ref{thm:B'} in terms of
generalized Markoff maps.

Let $H \in {\rm SL}(2,\mathbf Z) \cong \MCG$ be the monodromy of
the once-punctured torus bundle $M$ over the circle. Let $\hat
\rho:\pi_1(M) \rightarrow \PSLtwoC$ be the representation for the
hyperbolic structure on $M$ and $\rho:\pi_1(M) \rightarrow
\SLtwoC$ be a lift of the representation to $\SLtwoC$. Also denote
by $\rho$ its restriction to $\Gamma = \pi_1(\mathbb T)$. Recall
that we regard $\mathbb T$ as the fiber of $M$. Suppose this
$\rho: \Gamma \rightarrow \SLtwoC$ is a $\tau$-representation and
let $\phi \in {\bf \Phi}_{\mu}$ be the corresponding $\mu$-Markoff
map, where $\mu=\tau +2$. Then $\phi$ is invariant under $H$ as
defined in the statement of Theorem \ref{thm:stabilized}. With the
notation above and in \S \ref{s:generalized Markoff maps}, Theorem
\ref{thm:B'} can be reformulated as follows.

\vskip 10pt

\begin{thm}\label{thm:longitude} If $M$ has discrete length
spectrum on its torus fiber then
\begin{eqnarray}\label{eqn:stabilized}
\sum_{[X] \in {\Omega}_{\rm L}/\langle H\rangle} {\mathfrak
h}\big( \phi([X])\big)=\lambda \mod 2 \pi i,
\end{eqnarray}
where the sum converges absolutely.
\end{thm}

\vskip 10pt

\subsection{Proof of Theorem \ref{thm:B'}}\label{s:proof of thm B'} %

In this subsection we give the proof of Theorem \ref{thm:B'}.

\vskip 10pt

In order to prove Theorem \ref{thm:B'}, we need to compute the sum
on the left hand of (\ref{eqn:stabilized}). We shall assume that
the monodromy is positive, the other case is similar. Let us fix
an orientation on the meridian of $\partial M$ consistent with the
orientation of the fibre. We shall use the upper half-space model
of ${\mathbb H}^3$, so that its ideal boundary is identified with
the extended complex plane, ${\mathbf C}_{\infty}={\mathbf C} \cup
\{ \infty \}$, which has $\PSLtwoC$ acting in the usual way. We
can normalize our representation $\rho: \pi_1(M) \rightarrow
\SLtwoC$ so that $0$ and $\infty$ are the fixed points of the
images of the meridian and longitude. Then the developing map of
the similarity structure on $\partial M$ maps to  $\widetilde
\Cnozero$, the holonomy representation is given by the subgroup of
$\rho(\pi_1(M))$ generated by the images of the meridian and the
longitude.

Suppose $a, b$ are free generators of the fibre group $\Gamma$
such that the ordered pair of simple closed curves on $\mathbb T$
that they represent have algebraic intersection number $+1$. The
commutator, $[a,b]=aba^{-1}b^{-1}$ is peripheral in $\mathbb T$,
and so represents a meridian of $\partial M$.
Then $\rho([a, b])$ describes the translation $\zeta \mapsto e^{2\nu}\zeta$. %


\vskip 5pt

Now, to each region $X \in \Omega$, we shall associate an
incomplete semi-infinite geodesic $\Delta(X)$ in the developing
image of $M$ in $\HHH$ so that $\Delta(X)$ is normal to the
$z$-axis and has one end, $q(X)$, lying in the $z$-axis. Note that
$\Delta(X)$ will be well-defined up to the action of $\zeta
\mapsto e^{\nu}\zeta$, indeed all the geometric constructions are
really defined up to the action of $\zeta \mapsto e^{\nu}\zeta$.
Let $p(X) \in \Cnozero$ be the ideal point which is the orthogonal
projection of $q(X)$ on $\Cnozero$ along the direction of
$\Delta(X)$, we want to determine the position of $p(X)$, which
will then determine $\Delta(X)$. This is done as follows. We know
that $X$ corresponds to some simple closed curve $\gamma(X)$ on
$\mathbb T$. Let $\delta(X)$ be an arc on $\mathbb T$ with both
endpoints at the puncture such that $\gamma(X) \cap \delta(X) =
\emptyset$. The homotopy class of $\delta(X)$ relative to its
endpoints is well-defined. We are identifying $\mathbb T$ with a
fibre of $M$, so $\mathbb T$ is naturally homotopy equivalent to
the infinite cyclic cover of $M$. Under this equivalence,
$\delta(X)$ has a unique realization as a geodesic $\Delta(X)$ in
$M$ whose developing images to $\HHH$ are normal to the missing
lines. Choose a lift of $\Delta(X)$ to ${\mathbb H}^3$ which is
normal to the $z$-axis. Any other choice of lift would give us an
image of $\Delta(X)$ (hence an image of $p(X)$) under the cyclic
action generated by $\zeta \mapsto e^{\nu}\zeta$.

\vskip 10pt

Consider the bi-infinite sequence, $(X_j)_{j \in {\mathbf Z}}$, of
all regions of $\Omega_{\rm L}$ adjacent to the invariant path in
$\Sigma$ under $H$. We may choose $\Delta(X_j)$ so that the
transformation $\zeta \mapsto e^{\lambda}\zeta$ acts on them as a
shift. This is guaranteed by the geometric triangulation of
$\partial M$ described above.

It is not hard to see that in the geometric triangulation of
$\Cnozero$, the vertex corresponding to $X_j$ is given by one of
the images of $p(X_j)$ under the action generated by $\zeta
\mapsto e^{\nu}\zeta$, so we may as well assume that it actually
equals $p_j = p(X_j)$. The choice of $p_j$ naturally determines
that of $p_{j-1}$ and $p_{j+1}$, and so, inductively, $p_j$ for
all $j \in {\mathbf Z}$. Now the sequence $(p_j)$ is periodic
under the transformation corresponding to the longitude of
$\partial M$. This transformation is given by $\zeta \mapsto
e^{\lambda}\zeta$, where $\lambda=\lambda(\partial M)$ is the
length of the longitude. This corresponds to the action of
$\langle H \rangle$ on $\Omega$ which has the effect of shifting
the sequence $(X_j)$. Let $m$ be the number of steps through which
this sequence is shifted. Thus, $\lambda = \log p_m -  \log p_0 =
\sum_{\,j=1}^{\,m} (\log p_j - \log p_{j-1}) \mod 2 \pi i$. We
thus want to compute the numbers $\log p_j - \log p_{j-1} \mod 2
\pi i$. Let $\vec e_j$ be the directed edge given by $X_j \cap
X_{j-1}$, whose head lies in $\beta$ and let $C_{\rm L}$ be the
set $\{ \vec e_1, \cdots, \vec e_m \}$ of directed edges.

Fix some $j \in \{ 1, \cdots, m \}$. Let $X = X_{j-1}$, $Y = X_j$,
and let $Z$ be the region at the head of $\vec e_j$. As described
earlier, we can find free generators $a, b$ for $\Gamma$ which
correspond, respectively, to the regions X and Y. Moreover, we can
suppose that $a$ and $b$ are as described earlier and $Z$ is
represented by $ab$. By our discussion of the geometric meaning of
$\Psi$ in Appendix A, especially, by Lemma \ref{lem:geom of Psi},
we have
\begin{eqnarray*}
\log p_j-\log p_{j-1} = \Psi(x,y,z) \mod 2 \pi i,
\end{eqnarray*}
where $x = {\rm tr}\,\rho(a) = \phi(X)$, $y = {\rm tr}\,\rho(b) = \phi(Y)$ %
and $z = {\rm tr}\,\rho(ab) = \phi(Z)$. Thus
\begin{eqnarray*}
\log p_j - \log p_{j-1} = \Psi(\vec{e_j}) \mod 2 \pi i.
\end{eqnarray*}

It follows that in $\CmodTwoPiIZ$ we have
\begin{eqnarray*}
\lambda (\partial M) &=& \log p_m - \log p_0
\,=\,\sum_{j=1}^{m}(\log p_j - \log p_{j-1}) \\
&=&\sum_{\vec e \in C_{\rm L}}\psi(\vec e)
\,=\,\sum_{[X]\in \Omega_{\rm L}/\langle H \rangle} {\mathfrak h}(\phi([X])) \\
&=&\sum_{\sigma \in {\mathscr S}_{\rm
L}}\log\bigg(\frac{e^{\,\nu}+e^{\,l(\sigma)}}{e^{-\nu}+e^{\,l(\sigma)}}\bigg).
\end{eqnarray*}
This proves Theorem \ref{thm:B'} under certain choices of the
orientations involved. \square

\vskip 15pt

{\bf Proof of Corollary \ref{cor:Dehn}.}\,\, By Thurston's
Hyperbolic Dehn Surgery Theorem (see \cite{thurston1978notes}),
$\overline{M}(p/q)$ has a complete hyperbolic structure for all
but a finite number of surgery slopes $p/q \in \mathbf Q$. By
further excluding a finite number of slopes, we may assume that
the induced incomplete hyperbolic structure on $M$ is obtained by
a slight deformation from the unique complete hyperbolic structure
on $M$. By Theorem \ref{thm:openforbundle}, $M$ then has discrete
length spectrum on its torus fiber $\mathbb T$ and the conclusion
follows from Theorem \ref{thm:B'}. \square


\vskip 20pt

\section{{\bf Appendix A - Geometric meanings of the functions
 $\mathfrak h$ and $\Psi$}}
\label{s:geometric meaning of Psi} %
\vskip 10pt

In this appendix we explore the geometric meaning of the gap
function $\mathfrak h(x)=\mathfrak h_\tau(x)$ used earlier as well
as the geometric meaning of the function $\Psi(x,y,z)$ defined and
used in \S \ref{s:proof of theorem mu}. In particular, the
geometric meaning of $\Psi(x,y,z)$ was used in an essential way in
\S \ref{s:proof of thm B'}.

\vskip 10pt

\noindent {\bf Naturally oriented axes.}\,\, To each $A \in {\rm
SL}(2,\mathbf C)\backslash \{\pm I\}$, we associate its {\it
naturally oriented axis} ${\bf a}(A)$ as follows:
\begin{itemize}
\item[(i)] when $A$ represents a loxodromic (including hyperbolic)
isometry of $\mathbb H^3$, the orientation of ${\bf a}(A)$ is
directed from its repelling fixed ideal point to its attracting
fixed ideal point; \item[(ii)] when $A$ represents an elliptic
isometry of $\mathbb H^3$ but not an involution, the orientation
of ${\bf a}(A)$ is defined so that it has rotation angle in $(0,
\pi)$ with respect to ${\bf a}(A)$; and \item[(iii)] when
$A=(A_{ij})_{2 \times 2}$ is an involution, the orientation of
${\bf a}(A)$ is the same as Fenchel defined in
\cite{fenchel1989book}, that is, ${\bf a}(A)$ is directed from
${\rm Fix}^{-}(A)$ to ${\rm Fix}^{+}(A)$, where the ideal fixed
points ${\rm Fix}^{\mp}(A)$ of $A$ are distinguished as follows:
\begin{itemize}
\item[---] ${\rm Fix}^{-}(A):=(A_{11} - i)/A_{21}$ and ${\rm
Fix}^{+}(A):=(A_{11} + i)/A_{21}$ if $A_{21} \neq 0$; \item[---]
${\rm Fix}^{-}(A):=A_{12}i/2$ and ${\rm Fix}^{+}(A):=\infty$ if
$A_{21} = 0, A_{11} = i$; and \item[---] ${\rm
Fix}^{-}(A):=\infty$ and ${\rm Fix}^{+}(A):=-A_{12}i/2$ if $A_{21}
= 0, A_{11} = -i$.
\end{itemize}
\end{itemize}

\vskip 5pt

\noindent {\it Remark.}\,\, If $A$ is an involution, then ${\bf
a}(-A)={\bf a}(A^{-1})$ has the opposite orientation as ${\bf
a}(A)$. If $A$ is not an involution, then ${\bf a}(-A)={\bf a}(A)$
while ${\bf a}(A^{-1})$ has the opposite orientation as ${\bf
a}(A)$.

\vskip 10pt

\noindent {\bf Complex length $\Delta_{{\bf n}}({\bf l},{\bf
m})$.}\,\, Let us recall the complex length between oriented lines
in $\mathbb H^3$ as defined by Fenchel in \S V.3,
\cite{fenchel1989book}. Let ${\bf l,m,n} \in {\rm SL}(2,\mathbf
C)$ be involutions (that is, ${\bf l}^2={\bf m}^2={\bf n}^2=-I $)
representing oriented lines in $\mathbb H^3$ such that ${\bf n}$
is an oriented common normal to ${\bf l}$ and ${\bf m}$. (We say
two oriented lines in ${\mathbb H}^3$ are {\it normal} to each
other if they intersect orthogonally.) Then $\Delta=\Delta_{{\bf
n}}({\bf l},{\bf m}) \in {\mathbf C}/2 \pi i {\mathbf Z}$, the
{\it complex length} from ${\bf l}$ to ${\bf m}$ along ${\bf n}$,
is defined by
\begin{eqnarray}
\cosh \Delta = -(1/2)\,{\rm tr}({\bf ml}), \quad \sinh \Delta = -(i/2)\,{\rm tr}({\bf mnl}). %
\end{eqnarray}
Geometrically, $\Re \Delta \in \mathbf R$ is the signed hyperbolic
length from the point ${\bf l} \cap {\bf n}$ to the point ${\bf m}
\cap {\bf n}$ measured in the direction of ${\bf n}$, and $\Im
\Delta \in (-\pi, \pi]$ is the angle rotated from the the
direction of ${\bf l}$ to that of ${\bf m}$ measured along the
direction of ${\bf n}$.

\vskip 5pt

It is easy to see that, modulo $2\pi i$, $\Delta_{{\bf n}}(-{\bf
l},-{\bf m})=\Delta_{{\bf n}}({\bf l},{\bf m})$ and $\Delta_{{\bf
n}}(-{\bf l},{\bf m})=\Delta_{{\bf n}}({\bf l},-{\bf
m})=\Delta_{{\bf n}}({\bf l},{\bf m})+\pi i$.

\vskip 10pt

\noindent {\bf Geometric meaning of the gap function $\mathfrak
h={\mathfrak h}_{\tau}$.}\,\, The function ${\mathfrak
h}={\mathfrak h}_{\tau}$ has the following geometric
interpretation as the gap function as used in
\cite{tan-wong-zhang2004cone-surfaces}.

\vskip 5pt

Let $\phi \in {\Phi}_{\mu}$ be a $\mu$-Markoff map, where $\mu
\neq 0,4$, and $\rho: \Gamma \rightarrow {\rm SL}(2,\mathbf C)$ be
the corresponding $\tau$-representation, where $\tau=\mu-2$.

\vskip 5pt

Let $a, b$ be an arbitrary pair of generators of $\Gamma$ and let
$X, Y \in \Omega$ be the two regions which correspond to $a, b$
respectively. Let $A=\rho(a), B=\rho(b)$. Then by our convention,
$x=\phi(X)={\rm tr}A, \, y=\phi(Y)={\rm tr}B$. Let $z={\rm tr}BA$.

\vskip 5pt

Consider the commutator $[B^{-1}, A^{-1}]=(B^{-1}A^{-1}B)A$. It
follows from the trace identity (\ref{eqn:2+tr[A,B]=}) that ${\rm
tr}[B^{-1}, A^{-1}]=x^2+y^2+z^2-xyz-2=\mu-2=\tau$.

\vskip 5pt

If $A$ is loxodromic (including hyperbolic), so is
$B^{-1}A^{-1}B$. Let the attracting and repelling fixed points of
$A$ be respectively denoted as ${\rm Fix}^{+}(A)$ and ${\rm
Fix}^{-}(A)$. Similarly we have ${\rm Fix}^{+}(B^{-1}A^{-1}B)$ and
${\rm Fix}^{-}(B^{-1}A^{-1}B)$.

\vskip 5pt

We denote the oriented line in ${\mathbb H}^3$ which is normal to
the axis ${\bf a}(B^{-1}A^{-1}BA)$ and has ${\rm Fix}^{+}(A)$ as
its ending ideal point by $\lfloor {\bf a}(B^{-1}A^{-1}BA), {\rm
Fix}^{+}(A) \rfloor$. Similarly, we have another oriented line
$\lfloor {\bf a}(B^{-1}A^{-1}BA), {\rm Fix}^{-}(B^{-1}A^{-1}B)
\rfloor$.

\vskip 5pt

\begin{figure}
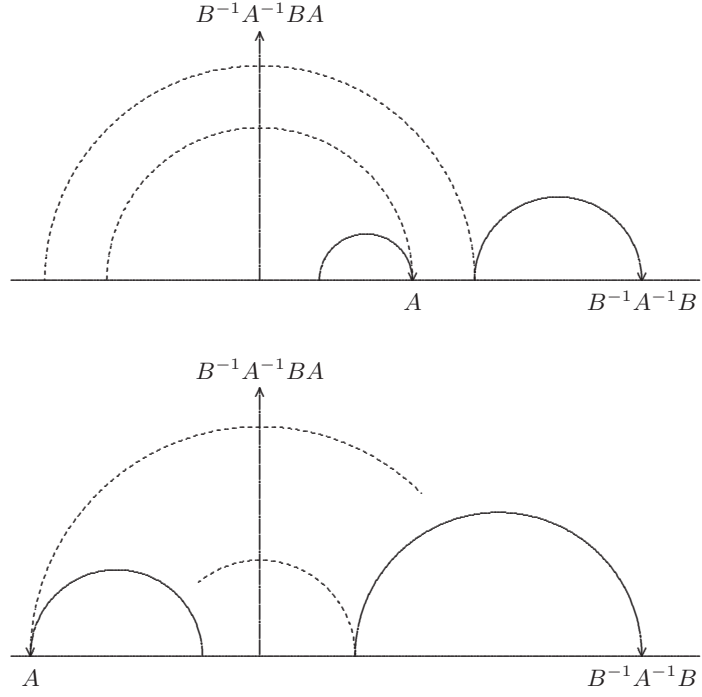



\begin{center}
\mbox{

\beginpicture
\setcoordinatesystem units <0.50in,0.50in>

\setplotarea x from -2.6 to 4.6, y from 0 to 2.8

\plot -2.6 0 4.6 0 /

\plot 0 0 0 2.6 /

\circulararc 180 degrees from 1.6 0 center at 1.1125 0

\circulararc 180 degrees from 4 0 center at 3.125 0

\put {\mbox{\small $B^{-1}A^{-1}BA$}} [ct] <0mm,4mm> at 0 2.6

\put {\mbox{\small $A$}} [cb] <0mm,-4mm> at 1.6 0

\put {\mbox{\small $B^{-1}A^{-1}B$}} [cb] <0mm,-4mm> at 4 0


\plot 0.04 2.52  0 2.6  -0.04 2.52 /

\plot 1.64 0.08  1.6 0 1.56 0.08 /

\plot 4.04 0.08  4 0 3.96 0.08 /

\setdashes<1.5pt>

\circulararc 180 degrees from 1.6 0 center at 0 0

\circulararc 180 degrees from 2.25 0 center at 0 0

\endpicture}

\end{center}

\vskip 15pt

\begin{center}
\mbox{

\beginpicture
\setcoordinatesystem units <0.50in,0.50in>

\setplotarea x from -2.6 to 4.6, y from 0 to 2.8

\plot -2.6 0 4.6 0 /

\plot 0 0 0 2.8 /

\circulararc -180 degrees from -2.4 0 center at -1.5 0

\circulararc 180 degrees from 4 0 center at 2.5 0

\put {\mbox{\small $B^{-1}A^{-1}BA$}} [ct] <0mm,4mm> at 0 2.8

\put {\mbox{\small $A$}} [cb] <0mm,-4mm> at -2.4 0

\put {\mbox{\small $B^{-1}A^{-1}B$}} [cb] <0mm,-4mm> at 4 0

\plot 0.04 2.72  0 2.8  -0.04 2.72 /

\plot -2.44 0.08  -2.4 0 -2.36 0.08 /

\plot 4.04 0.08  4 0 3.96 0.08 /

\setdashes<1.5pt>\circulararc -135 degrees from -2.4 0 center at 0
0

\circulararc 130 degrees from 1 0 center at 0 0
\endpicture}

\end{center}


\caption{The complex gap}\label{fig:gap}
\end{figure}

Then the complex length from the oriented line $\lfloor {\bf
a}(B^{-1}A^{-1}BA), {\rm Fix}^{+}(A) \rfloor$ to $\lfloor {\bf
a}(B^{-1}A^{-1}BA), {\rm Fix}^{-}(B^{-1}A^{-1}B) \rfloor$ along
${\bf a}(B^{-1}A^{-1}BA)$ is exactly given by ${\mathfrak
h}(x)={\mathfrak h}_\tau(x)$. See Figure \ref{fig:gap} for an
illustration in the cases where the isometries $A,B^{-1}A^{-1}B$
and $B^{-1}A^{-1}BA$ have coplanar non-intersecting axes.

\vskip 5pt

\begin{lem}\label{lem:frak h=gap} With the above notation, we have
\begin{eqnarray}\label{eqn:gap=frak h(x)}
{\mathfrak h}_{\tau}(x)\!\!&=&\!\!
{\Delta}_{{\bf a}(B^{-1}A^{-1}BA)} \big(\lfloor {\bf a}(B^{-1}A^{-1}BA), {\rm Fix}^{+}(A) \rfloor, \nonumber \\ %
& &\hspace{64.5pt} \lfloor {\bf a}(B^{-1}A^{-1}BA), {\rm Fix}^{-}(B^{-1}A^{-1}B) \rfloor \big). %
\end{eqnarray}
\end{lem}

\vskip 5pt

\noindent {\it Remark.}\,\, Here it is important to note that the
gap value ${\mathfrak h}(x)$ is independent of the choice of $b$
in the generating pair $a,b$ of $\Gamma$. Note also that ${\rm
Fix}^{-}(B^{-1}A^{-1}B)={\rm Fix}^{+}(B^{-1}AB)$ and that
$B^{-1}AB$ is a conjugate of $A$.

\vskip 5pt

Lemma \ref{lem:frak h=gap} is in fact a special case of the
following general lemma.

\vskip 5pt

\begin{lem}\label{lem:gap general} Let $A,B \in \SLtwoC$. Then the
complex distance from the oriented line $\lfloor{\bf a}(BA), {\rm
Fix}^{+}(A)\rfloor$ to the oriented line $\lfloor{\bf a}(BA), {\rm
Fix}^{-}(B)\rfloor$ measured along the oriented line ${\bf a}(BA)$
is given by
\begin{eqnarray}
2 \tanh^{-1} \left( \frac{ \sinh\frac{l(-BA)}{2} }{ \cosh
\frac{l(-BA)}{2} + \exp \big(\frac{l(A)}{2}+\frac{l(B)}{2}\big) }
\right)& & \nonumber \\ & & \hspace{-120pt} =\log \left(  \frac{
\exp\frac{l(-BA)}{2} + \exp
\big(\frac{l(A)}{2}+\frac{l(B)}{2}\big) }{ \exp
\big(-\frac{l(-BA)}{2}\big) + \exp
\big(\frac{l(A)}{2}+\frac{l(B)}{2}\big) } \right),
\end{eqnarray}
where ${l(A)}/{2},\,{l(B)}/{2},\, {l(-BA)}/{2}$ are the half
translation lengths of $A, B, -BA$ as defined by
{\rm(\ref{eqn:l(A)/2=})} in \S {\rm\ref{s:notation+results}}.
\end{lem}

\vskip 5pt

We omit the proof here. A detailed proof by cross-ratio
calculations can be found in \cite{zhang2004thesis}. For another
proof, one may adapt Mirzakhani's proof in the real case as given
in \cite{mirzakhani2004preprint} to let it work for the general
case here by using Fenchel's cosine rule for oriented right angled
hexagons given in \S VI.2 of \cite{fenchel1989book}.

\vskip 10pt

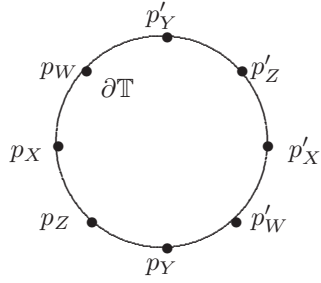
\begin{figure}
\setlength{\unitlength}{1mm} 
\begin{picture}(60,30)
\circulararc 360 degrees from 120 40  center at 80 40 \thicklines
\put(20,20){$\partial {\mathbb T}$} \put(8,12){$p_X$}
\put(45,12){$p'_X$} \put(26,-3){$p_Y$} \put(26,30){$p'_Y$}
\put(12,23){$p_W$} \put(40,3){$p'_W$} \put(12,3){$p_Z$}
\put(40,23){$p'_Z$}

\put(12.3,12.5){ $\bullet$} \put(40.2,12.5){ $\bullet$}

\put(26.8,27){ $\bullet$}\put(26.8,-1){ $\bullet$}

\put(36.8,22.5){ $\bullet$}\put(16.8,2.5){ $\bullet$}

\put(36,2.5){ $\bullet$}\put(16.1,22.5){ $\bullet$}
\end{picture}
\caption{Intersections of $\delta_X, \delta_Y, \delta_Z, \delta_W$
with $\partial {\mathbb T}$} \label{fig:Psi}
\end{figure}

\noindent {\bf Geometric meaning of the function $\Psi$.}\,\, We
first  show the picture and give an informal treatment in the case
where $\phi$ corresponds to the holonomy representation $\rho$ of
a real hyperbolic one-holed torus ${\mathbb T}$ with geodesic
boundary $\partial {\mathbb T}$. For the limiting case where the
boundary is a cusp, see \cite{bowditch1998plms} or
\cite{akiyoshi-miyachi-sakuma2004cm355}. Let $e
\leftrightarrow(X,Y;Z,W)$, and let $\gamma_X, \gamma_Y, \gamma_Z,
\gamma_W$ be the simple closed geodesics on ${\mathbb T}$
corresponding to $X,Y,Z$ and $W$ respectively. $\gamma_X$ and
$\gamma_Y$ intersect each other once, at a Weierstrass point of
${\mathbb T}$, and $\gamma_Z$ and $\gamma_W$ are the two unique
simple closed geodesics which intersect both $\gamma_X$ and
$\gamma_Y$ exactly once. The points of intersection of these
geodesics are the Weierstrass points, there are exactly three, and
each geodesic passes through two and misses the third.
Corresponding to $\gamma_X$ is a unique geodesic segment
$\delta_X$ disjoint from $\gamma_X$, passing through the third
Weierstrass point  and intersecting $\partial  {\mathbb T}$ at
both ends normally. Denote these points of intersection on
$\partial {\mathbb T}$ by $p_X$ and $p'_X$, they are necessarily
on opposite ends of $\partial {\mathbb T}$. Similarly, we have
$\delta_Y$, $\delta_Z$ and $\delta_W$, and the corresponding
points of intersection with $\partial {\mathbb T}$. The eight
points of intersection on $\partial {\mathbb T}$ are shown in
Figure \ref{fig:Psi}. When we consider the intervals $[p_X,p_Y]$
and $[p_X,p'_Y]$ in $\partial {\mathbb T}$, one of them contains
one of the points $p_Z$ or $p'_Z$, the other contains one of the
points $p_W$ or $p'_W$. Then $\Psi(x,y,z)$ is the length of the
interval which contains $p_W$ or $p'_W$. (Note that we get the
same answer if we used $p'_X$ instead of $p_X$). $\Psi(x,y,w)$ is
then the length of the other interval. It is not difficult to see
from this description that $\Psi$ satisfies equations
(\ref{eqn:Psi+Psi+Psi=nu}) and (\ref{eqn:Psi+Psi=nu}).


We now consider the more general case. Given $x,y,z \in {\mathbf
C}$, set
$$\mu:=x^2+y^2+z^2-xyz$$ and assume $\mu \neq 0,4$. As explained
in \S \ref{s:generalized Markoff maps}, there exist $A,B \in
\SLtwoC$, unique up to conjugation in $\SLtwoC$, so that the
traces of $A,B$ and $C:=(BA)^{-1}$ are respectively $x,y,z$, that
is,
$$
{\rm tr}A=x, \quad {\rm tr}B=y, \quad {\rm tr}AB={\rm tr}BA={\rm tr}C=z. %
$$
Note that $CBA=I$. Denote
$[B^{-1},A^{-1}]=B^{-1}A^{-1}BA=B^{-1}A^{-1}C^{-1}$ and let
$\tau:={\rm tr}[B^{-1},A^{-1}]={\rm tr} B^{-1}A^{-1}BA$. Then by
the trace identity (\ref{eqn:2+tr[A,B]=}) in $\SLtwoC$, we have
$\tau=\mu-2 \neq \pm 2$. Let
$\nu=\cosh^{-1}(-\tau/2)=\cosh^{-1}(1-\mu/2)$.

\vskip 3pt

There are matrices $Q,R,P \in {\rm SL}(2,\mathbf C)$ such that
$Q^{2}=R^{2}=P^{2}=-I$ and $A=-RQ, B=-PR, C=-QP$. Hence 
$[B^{-1},A^{-1}]=-(RPQ)^2.$ Here matrices in ${\rm SL}(2,\mathbf
C)$ are considered as orientation-preserving isometries of
$\mathbb H^3$, the upper half-space model of the hyperbolic
3-space. In fact, $Q,R,P$ are involutions representing
appropriately oriented lines in $\HHH$ which are common normals to
${\bf a}(C)$ and ${\bf a}(A)$, to ${\bf a}(A)$ and ${\bf a}(B)$,
and to ${\bf a}(B)$ and ${\bf a}(C)$, respectively, as explained
by Fenchel \cite{fenchel1989book}. These oriented lines
represented by $Q,R,P$ will play the role of the Weierstrass
points used in the real case. Note that the negative signs in this
paragraph are important for later use of Fenchel's cosine and sine
rules for right angled hexagons in this section.

\vskip 5pt

Consider the conjugation $A \mapsto KAK^{-1}$ by a matrix $K \in
{\rm SL}(2,\mathbf C)$. It preserves the natural orientations of
axes defined above in the following sense.

\vskip 10pt

\begin{lem}
For $A \in {\rm SL}(2,\mathbf C)$, we have ${\bf
a}(KAK^{-1})=K{\bf a}(A)K^{-1}$. In particular, when $K$ is an
involution, that is $K^2=-I$, the conjugation is given by $A
\mapsto -KAK$.
\end{lem}
%

Now consider the schematic figure as illustrated in Figure
\ref{fig:alpha+beta+gamma=nu}. It is easy to check that $RPQ
\leftrightarrow QRP \leftrightarrow PQR \leftrightarrow RPQ$ by
conjugation by $Q, P, R$ respectively. Hence ${\bf a}(RPQ)
\leftrightarrow {\bf a}(QRP) \leftrightarrow {\bf a}(PQR)
\leftrightarrow {\bf a}(RPQ)$ by conjugation by $Q, P, R$
respectively.

\vskip 10pt

\noindent {\bf Right angled hexagons ${\mathbf H}({\bf l}, {\bf
m}, {\bf n})$.}\,\, We shall use the following notation for
oriented common normals and right angled hexagons in $\HHH$.

\vskip 3pt

(i) For each ordered pair of oriented lines ${\bf l}$ and ${\bf
m}$ in $\HHH$, let $\lfloor {\bf l}, {\bf m} \rfloor$ denote a
definitely chosen oriented common normal to them, so that $\lfloor
{\bf l}, {\bf m} \rfloor$ and $\lfloor {\bf m}, {\bf l} \rfloor$
always have opposite directions. For example, we may assume
$\lfloor {\bf l}, {\bf m} \rfloor$ is directed from ${\bf l}$ to
${\bf m}$ when ${\bf l}$ and ${\bf m}$ are disjoint.

\vskip 3pt

(ii) Given oriented lines ${\bf l}$, ${\bf m}$, ${\bf n}$ in
$\HHH$, we use ${\mathbf H}({\bf l}, {\bf m}, {\bf n})$ to denote
the right angled hexagon in the sense of Fenchel
\cite{fenchel1989book} formed by the oriented lines

\centerline{ ${\bf l}; \lfloor {\bf l}, {\bf m} \rfloor; {\bf m};
\lfloor {\bf m}, {\bf n} \rfloor;{\bf n}; \lfloor {\bf n}, {\bf l}
\rfloor$ }

\noindent in this cyclic order.

\vskip 10pt


\begin{figure}

\begin{center}\mbox{
\beginpicture
\setcoordinatesystem units <1mm,1mm>

\setplotarea x from 0 to 120, y from 0 to 80

\plot 60 70 110 60 90 20 60 70 10 50 40 10 90 20 / \plot 60 70 40
10 /

\plot 85 65 75 45 65 15 50 40 35 60 25 30 50 40 75 45 100 40 85 65
/

\put {\mbox{\scriptsize $RPQ$}} [cb] <0mm,1.3mm> at 60 70

\put {\mbox{\scriptsize $Q$}} [cb] <0mm,1.3mm> at 85 65

\put {\mbox{\scriptsize $QRP$}} [cb] <0mm,1.3mm> at 110 60

\put {\mbox{\scriptsize $P$}} [lt] <0mm,-1mm> at 100 40

\put {\mbox{\scriptsize $PQR$}} [ct] <0mm,-1.5mm> at 90 20

\put {\mbox{\scriptsize $R$}} [ct] <-2.5mm,-1.5mm> at 75 45

\put {\mbox{\scriptsize $-RQR$}} [ct] <-0.5mm,-1.5mm> at 65 15

\put {\mbox{\scriptsize $-RQRPR$}} [ct] <-0.5mm,-1.5mm> at 40 10

\put {\mbox{\scriptsize $RPRPR$}} [ct] <-5.5mm,-1mm> at 25 30

\put {\mbox{\scriptsize $RPRPQPR$}} [ct] <-4.5mm,3.5mm> at 10 50

\put {\mbox{\scriptsize $RPQPR$}} [ct] <-1.5mm,3.5mm> at 35 60

\put {\mbox{\scriptsize $-RPR$}} [ct] <-6.5mm,1.5mm> at 50 40

\put {\mbox{\scriptsize $\alpha$}} [ct] <3.5mm,-2.5mm> at 60 70

\put {\mbox{\scriptsize $\beta$}} [ct] <0mm,-3.5mm> at 60 70

\put {\mbox{\scriptsize $\gamma$}} [ct] <-3.5mm,-3.5mm> at 60 70

\put {\mbox{\scriptsize $\gamma$}} [ct] <-2.5mm,-0.5mm> at 110 60

\put {\mbox{\scriptsize $\beta$}} [cb] <0mm,2.5mm> at 90 20

\put {\mbox{\scriptsize $\alpha$}} [cb] <-3mm,0.8mm> at 90 20

\put {\mbox{\scriptsize $\beta$}} [cb] <-0.5mm,2.5mm> at 40 10

\put {\mbox{\scriptsize $\gamma$}} [ct] <2.5mm,3mm> at 40 10

\put {\mbox{\scriptsize $\alpha$}} [ct] <3mm,-0.5mm> at 10 50

\endpicture
}\end{center} \caption{} \label{fig:alpha+beta+gamma=nu}
\end{figure}

Thus the oriented lines
\begin{eqnarray*}
& & {\bf a}(RPQ); \lfloor {\bf a}(RPQ), {\bf a}(PQR) \rfloor; \\
& & {\bf a}(PQR); \lfloor {\bf a}(PQR), {\bf a}(QRP) \rfloor; \\
& & {\bf a}(QRP); \lfloor {\bf a}(QRP), {\bf a}(RPQ) \rfloor,
\end{eqnarray*}
in this cyclic order, form the right angled hexagon ${\mathbf
H}({\bf a}(RPQ), {\bf a}(PQR), {\bf a}(QRP))$. Since, as we have
observed above, ${\bf a}(RPQ) \leftrightarrow {\bf a}(QRP)
\leftrightarrow {\bf a}(PQR) \leftrightarrow {\bf a}(RPQ)$ by
conjugation by $Q, P, R$ respectively, we know ${\mathbf H}({\bf
a}(RPQ), {\bf a}(PQR), {\bf a}(QRP))$ has the oriented lines ${\bf
a}(R), {\bf a}(P), {\bf a}(Q)$ as the `midpoints' of its three
sides, that is,
\begin{eqnarray*}
&&\Delta_{\lfloor {\bf a}(RPQ), {\bf a}(PQR) \rfloor} \big({\bf
a}(RPQ), {\bf a}(R) \big) = \Delta_{\lfloor {\bf a}(RPQ), {\bf
a}(PQR) \rfloor} \big({\bf a}(R), {\bf a}(PQR) \big) =: \tilde{c},
\\
&&\Delta_{\lfloor {\bf a}(PQR), {\bf a}(QRP) \rfloor} \big({\bf
a}(PQR), {\bf a}(P) \big) = \Delta_{\lfloor {\bf a}(PQR), {\bf
a}(QRP) \rfloor} \big({\bf a}(P), {\bf a}(QRP) \big) =: \tilde{a},
\\
&&\Delta_{\lfloor {\bf a}(QRP), {\bf a}(RPQ) \rfloor} \big({\bf
a}(QRP), {\bf a}(Q) \big) = \Delta_{\lfloor {\bf a}(QRP), {\bf
a}(RPQ) \rfloor} \big({\bf a}(Q), {\bf a}(RPQ) \big) =: \tilde{b}.
\end{eqnarray*}

\vskip 5pt

\noindent Let the complex lengths of the other three sides of the
right angled hexagon ${\mathbf H}({\bf a}(RPQ), {\bf a}(PQR), {\bf
a}(QRP))$ be denoted as
\begin{eqnarray*}
&&\Delta_{{\bf a}(RPQ)} \big(\lfloor {\bf a}(QRP), {\bf a}(RPQ)
\rfloor, \lfloor {\bf a}(RPQ), {\bf a}(PQR) \rfloor \big) =: {\alpha}, %
\\
&&\Delta_{{\bf a}(PQR)} \big(\lfloor {\bf a}(RPQ), {\bf a}(PQR)
\rfloor, \lfloor {\bf a}(PQR), {\bf a}(QRP) \rfloor \big) =: {\beta}, %
\\
&&\Delta_{{\bf a}(QRP)} \big(\lfloor {\bf a}(PQR), {\bf a}(QRP)
\rfloor, \lfloor {\bf a}(QRP), {\bf a}(RPQ) \rfloor \big) =: {\gamma}. %
\end{eqnarray*}
\noindent Also let
\begin{eqnarray*}
&&\Delta_{{\bf a}(RQ)} \big({\bf a}(Q), {\bf a}(R)\big) =: \tilde{p}, \\
&&\Delta_{{\bf a}(PR)} \big({\bf a}(R), {\bf a}(P)\big) =: \tilde{q}, \\
&&\Delta_{{\bf a}(QP)} \big({\bf a}(P), {\bf a}(Q)\big) =: \tilde{r}. \\
\end{eqnarray*}
Then
\begin{eqnarray*}
& & \cosh\tilde{p}=(1/2)\,{\rm tr}(-RQ)=(1/2)\,{\rm tr}A = x/2, \\
& & \cosh\tilde{q}=(1/2)\,{\rm tr}(-PR)=(1/2)\,{\rm tr}B = y/2, \\
& & \cosh\tilde{r}=(1/2)\,{\rm tr}(-QP)=(1/2)\,{\rm tr}C = z/2. %
\end{eqnarray*}

\vskip 5pt

\noindent Now applying Fenchel's cosine rules as in \S VI.2 of
\cite{fenchel1989book} to the right angled hexagons \vskip 5pt
\centerline{${\mathbf H}({\bf a}(Q), {\bf a}(P), {\bf a}(QRP))$
and ${\mathbf H}({\bf a}(RPQ), {\bf a}(PQR), {\bf a}(QRP))$}\vskip
5pt \noindent respectively we have
\begin{eqnarray}
\cosh \tilde{r} = \cosh \tilde{a} \, \cosh \tilde{b}
+ \sinh \tilde{a} \, \sinh \tilde{b} \, \cosh \gamma, \label{eqn:rab} %
\\
\cosh 2\tilde{c} = \cosh 2\tilde{a} \, \cosh 2\tilde{b} + \sinh
2\tilde{a} \, \sinh 2\tilde{b} \, \cosh \gamma. \label{eqn:2c2a2b}
\end{eqnarray}

\vskip 5pt

{\bf Claim.}
\begin{eqnarray}\label{eqn:ratio=k}
\frac{\cosh \tilde{a}}{\cosh \tilde{p}} = \frac{\cosh
\tilde{b}}{\cosh \tilde{q}} = \frac{\cosh \tilde{c}}{\cosh
\tilde{r}} = \kappa,
\end{eqnarray}
where ${\kappa}^2 = 4/\mu$.

\vskip 5pt

To prove the claim, multiplying \, $4 \cosh\tilde{a}
\cosh\tilde{b}$ \, to both sides of (\ref{eqn:rab}) gives
\begin{eqnarray*}
4 \cosh\tilde{a} \, \cosh\tilde{b} \, \cosh \tilde{r} = 4
\cosh^{2} \tilde{a} \,  \cosh^{2} \tilde{b} + \sinh 2\tilde{a} \,
\sinh 2\tilde{b} \,  \cosh \gamma.
\end{eqnarray*}

\noindent Comparing with (\ref{eqn:2c2a2b}) gives
\begin{eqnarray*}
\hspace{15pt} 4 \cosh\tilde{a} \, \cosh\tilde{b} \, \cosh
\tilde{r} - \cosh 2\tilde{c} = 4 \cosh^{2} \tilde{a} \, \cosh^{2}
\tilde{b} - \cosh 2\tilde{a} \, \cosh 2\tilde{b}.
\end{eqnarray*}

\noindent After simplification we have
\begin{eqnarray}\label{eqn:abr=}
(2\cosh\tilde{a})(2\cosh\tilde{b})(2\cosh\tilde{r})
=(2\cosh\tilde{a})^2+(2\cosh\tilde{b})^2+(2\cosh\tilde{c})^2-4.
\end{eqnarray}

\noindent Similarly we have
\begin{eqnarray}\label{eqn:aqc=}
(2\cosh\tilde{a})(2\cosh\tilde{q})(2\cosh\tilde{c})
=(2\cosh\tilde{a})^2+(2\cosh\tilde{b})^2+(2\cosh\tilde{c})^2-4,
\end{eqnarray}
\begin{eqnarray}\label{eqn:pbc=}
(2\cosh\tilde{p})(2\cosh\tilde{b})(2\cosh\tilde{c})
=(2\cosh\tilde{a})^2+(2\cosh\tilde{b})^2+(2\cosh\tilde{c})^2-4.
\end{eqnarray}

\noindent Hence from (\ref{eqn:abr=})--(\ref{eqn:pbc=}) we have
\begin{eqnarray}\label{eqn:ch a /ch p}
\frac{\cosh \tilde{a}}{\cosh \tilde{p}} = \frac{\cosh
\tilde{b}}{\cosh \tilde{q}} = \frac{\cosh \tilde{c}}{\cosh
\tilde{r}}.
\end{eqnarray}

\vskip 3pt

\noindent Let the common value in (\ref{eqn:ch a /ch p}) be
denoted $\kappa$. Then
\begin{eqnarray*}
& & 2\cosh \tilde{a} = \kappa \, 2\cosh\tilde{p} = \kappa\, x,\\
& & 2\cosh \tilde{b} = \kappa \, 2\cosh\tilde{q} = \kappa\, y,\\
& & 2\cosh \tilde{c} = \kappa \, 2\cosh\tilde{r} = \kappa\, z.
\end{eqnarray*}
Now from (\ref{eqn:abr=}) we have
\begin{eqnarray*}
{\kappa}^2xyz={\kappa}^2(x^2+y^2+z^2)-4,
\end{eqnarray*}
and hence (recalling that $\mu=x^2+y^2+z^2-xyz$)
\begin{eqnarray}\label{eqn:kappa^2=4/mu}
{\kappa}^2=4/\mu.
\end{eqnarray}
This proves the claim.

\vskip 10pt

It follows from the above proof that

\begin{lem}\label{lem:cosh gamma}
\begin{eqnarray}\label{eqn:cosh gamma}
\cosh \gamma = \frac{z/2-xy/\mu}{(x^2/\mu-1)^{1/2}(y^2/\mu-1)^{1/2}}, %
\end{eqnarray}
where the square roots do not necessarily have nonnegative real parts. %
\end{lem}

Actually, from (\ref{eqn:rab}), we have (since ${\kappa}^2=4/\mu$)
\begin{eqnarray*}
\cosh\gamma &=& \frac{\cosh \tilde{r} - \cosh \tilde{a} \cosh \tilde{b}}{\sinh \tilde{a} \sinh \tilde{b}} \\
&=& \frac{z/2-{\kappa}\,(x/2)\,{\kappa}\,(y/2)}{[({\kappa}\, x/2)^2-1]^{1/2}[({\kappa}\, y/2)^2-1]^{1/2}} \\
&=& \frac{z/2-xy/\mu}{(x^2/\mu-1)^{1/2}(y^2/\mu-1)^{1/2}}.
\end{eqnarray*}

\vskip 5pt

As an immediate corollary, we have

\vskip 5pt

\begin{cor}\label{cor:sinh gamma=pm}
\begin{eqnarray}\label{sinh gamma=}
\sinh \gamma = \pm \frac{(\sinh\nu)\,z/\mu}{(x^2/\mu-1)^{1/2}(y^2/\mu-1)^{1/2}}. %
\end{eqnarray}
\end{cor}

\vskip 5pt

The rest of this section is devoted to the determination of the
sign in (\ref{sinh gamma=}). Geometrically, we have

\vskip 10pt

\begin{lem}\label{lem:alpha+beta+gamma=nu}
$\quad (\alpha+\pi i)+(\beta+\pi i)+(\gamma+\pi i)=\nu \mod 2 \pi i$. %
\end{lem}

\vskip 5pt

\begin{pf} We refer to Figure \ref{fig:alpha+beta+gamma=nu}. Consider the following right angled hexagons
\begin{eqnarray*}
& &{\mathcal H_1}:={\mathbf H}({\bf a}(RPQ), {\bf a}(PQR), {\bf a}(QRP)), \\
& &{\mathcal H_2}:={\mathbf H}({\bf a}(PQR), {\bf a}(RPQ), {\bf a}(-RQRPR)), \\
& &{\mathcal H_3}:={\mathbf H}({\bf a}(RPRPQPR), {\bf a}(-RQRPR), {\bf a}(RPQ)). %
\end{eqnarray*}

\noindent It is easy to check that the conjugation by $R$ maps
${\mathcal H_1}$ to ${\mathcal H_2}$, and the conjugation by $RPR$
maps ${\mathcal H_2}$ to ${\mathcal H_3}$; thus the conjugation by
$PR$ maps ${\mathcal H_1}$ to ${\mathcal H_3}$. Since conjugations
preserve the relevant complex lengths, we have
\begin{eqnarray*}
& &\Delta_{{\bf a}(RPQ)} \big(\lfloor {\bf a}(PQR), {\bf a}(RPQ)
\rfloor, \lfloor {\bf a}(RPQ), {\bf a}(-RQRPR) \rfloor \big) = \beta, \\
& &\Delta_{{\bf a}(RPQ)} \big(\lfloor {\bf a}(-RQRPR), {\bf
a}(RPQ) \rfloor, \lfloor {\bf a}(RPQ), {\bf a}(RPRPQPR) \rfloor \big) = \gamma. %
\end{eqnarray*}
Thus
\begin{eqnarray*}
& &(\alpha+\pi i)+(\beta+\pi i)+(\gamma+\pi i) \\
&=& \Delta_{{\bf a}(RPQ)} \big(\lfloor {\bf a}(QRP), {\bf a}(RPQ)
\rfloor,
\lfloor {\bf a}(RPQ), {\bf a}(PQR) \rfloor \big) \\
&+&\Delta_{{\bf a}(RPQ)} \big(\lfloor {\bf a}(PQR), {\bf a}(RPQ)
\rfloor, \lfloor {\bf a}(RPQ), {\bf a}(-RQRPR) \rfloor \big) \\
&+&\Delta_{{\bf a}(RPQ)} \big(\lfloor {\bf a}(-RQRPR), {\bf
a}(RPQ) \rfloor,
\lfloor {\bf a}(RPQ), {\bf a}(RPRPQPR) \rfloor \big) + \pi i \\
&=& \Delta_{{\bf a}(RPQ)} \big(\lfloor {\bf a}(QRP), {\bf a}(RPQ) \rfloor, %
\lfloor {\bf a}(RPQ), {\bf a}(PQR) \rfloor \big) \\
&+&\Delta_{{\bf a}(RPQ)} \big(\lfloor {\bf a}(RPQ), {\bf a}(PQR)
\rfloor, \lfloor {\bf a}(-RQRPR), {\bf a}(RPQ) \rfloor \big) \\
&+&\Delta_{{\bf a}(RPQ)} \big(\lfloor {\bf a}(-RQRPR), {\bf a}(RPQ) \rfloor, %
\lfloor {\bf a}(RPQ), {\bf a}(RPRPQPR) \rfloor \big) + \pi i \\
&=& \Delta_{{\bf a}(RPQ)} \big(\lfloor {\bf a}(QRP), {\bf a}(RPQ)
\rfloor, \lfloor {\bf a}(RPQ), {\bf a}(RPQPRPR) \rfloor \big) + \pi i \\
&=&\Delta_{{\bf a}(RPQ)} \big(\lfloor {\bf a}(QRP), {\bf a}(RPQ)
\rfloor, \lfloor {\bf a}(RPRPQPR), {\bf a}(RPQ) \rfloor \big) \\
&=& \nu \mod 2 \pi i.
\end{eqnarray*}
The last equality follows from Lemmas \ref{lem:K translength} %
and \ref{lem:conjugation} below, since it is easy to see that $\nu
=\cosh^{-1}(-\tau/2)$ is the complex translation length of $RPQ$
(recall that $[B^{-1},A^{-1}]=-(RPQ)^2$ and $\tau={\rm tr}
[B^{-1},A^{-1}]$) and that the conjugation by $RPQ$ maps
$\lfloor {\bf a}(QRP), {\bf a}(RPQ) \rfloor$ %
to $\lfloor {\bf a}(RPRPQPR), {\bf a}(RPQ) \rfloor$. %
\end{pf}

\vskip 5pt

\begin{lem}\label{lem:K translength} The complex translation
length of $K \in \SLtwoC$ is given by
\begin{eqnarray}\label{eqn:l(K)}
l(K)=\cosh^{-1}\big(\textstyle{\frac12}{\rm tr}(K^2)\big).
\end{eqnarray}
\end{lem}

\vskip 5pt

\begin{pf} Recall that $l(K)/2=\cosh^{-1}\big(\frac12{\rm
tr}K\big)$. Hence
\begin{eqnarray*}
\cosh l(K) &=& 2 \cosh^2(l(K)/2) - 1 \\
& =& 2 \big(\textstyle{\frac12}{\rm tr}K\big)^2-1 \\
& =& \textstyle{\frac12}({\rm tr}^2K-2) \\
& =& \textstyle{\frac12}{\rm tr}(K^2), %
\end{eqnarray*}
from which (\ref{eqn:l(K)}) follows since $\Re\, l(K) \ge 0$.
\end{pf}

%
%

\begin{lem}\label{lem:conjugation} If $K \in \SLtwoC$ is non-parabolic
and $L \in \SLtwoC$ is a line matrix such that ${\bf a}(K) \perp
{\bf a}(L)$, then $KLK^{-1}$ is also a line matrix and the complex
translation length $l(K)$ of $K$ is given by
\begin{eqnarray}\label{eqn:l(K)=Delta}
l(K)=\Delta_{{\bf a}(K)}\big({\bf a}(L), {\bf a}(KLK^{-1})\big). %
\end{eqnarray}
\end{lem}

This is Lemma 2.17 in \cite{zhang2004thesis}, with a proof given
there.

\vskip 10pt

Now we can determine the signs in the expressions (\ref{sinh
gamma=}) etc as follows.

\begin{lem}\label{lem:sinh gamma}
\begin{eqnarray}\label{eqn:sinh gamma}
\sinh \gamma = -\frac{(\sinh\nu)\,z}{(x^2/\mu-1)^{1/2}(y^2/\mu-1)^{1/2}} %
\end{eqnarray}
\vskip 3pt \noindent and similarly for $\sinh \alpha$ and $\sinh \beta$. %
\end{lem}

\vskip 5pt

\begin{pf} Let
$\Psi(x,y,z) \in {\mathbf C}$ be defined as in \S \ref{s:proof of
theorem mu} by:
\begin{eqnarray}\label{varPsi(x,y,z)=log}
\Psi(x,y,z)= \log
\frac{[xy+(e^{\nu}-1)z]/\mu}{(x^2/\mu-1)^{1/2}(y^2/\mu-1)^{1/2}}
\end{eqnarray}
or equivalently, by the following two equations:
\begin{eqnarray}
\cosh \Psi(x,y,z)=
\frac{[xy-(\mu/2)z]/\mu}{(x^2/\mu-1)^{1/2}(y^2/\mu-1)^{1/2}},
\end{eqnarray}
\begin{eqnarray}
\sinh \Psi(x,y,z)=
\frac{(\sinh\nu)z/\mu}{(x^2/\mu-1)^{1/2}(y^2/\mu-1)^{1/2}}.
\end{eqnarray}
Similarly for $\Psi(y,z,x), \Psi(z,x,y) \in {\mathbf C}$. We have
from Proposition \ref{prop:properties of Psi}(i) that
\begin{eqnarray*}
\Psi(y,z,x) + \Psi(z,x,y) + \Psi(x,y,z) = \nu \mod 2 \pi i.
\end{eqnarray*}
On the other hand, it follows from Lemma \ref{lem:cosh gamma} that
$-\cosh \gamma =\cosh \Psi(x,y,z)$ etc and hence, modulo $2\pi i$,
\begin{eqnarray}\label{eqn:alpha+pi i}
\alpha+\pi i=\pm \Psi(y,z,x), \quad  \beta+\pi i=\pm \Psi(z,x,y),
\quad \gamma+\pi i=\pm \Psi(x,y,z).
\end{eqnarray}

Now applying Fenchel's sine rule as in \S VI.2 of
\cite{fenchel1989book} to the right angled hexagon \vskip 5pt
\centerline{${\mathbf H}({\bf a}(RPQ), {\bf a}(PQR), {\bf
a}(QRP))$} \vskip 3pt \noindent gives
\begin{eqnarray*}
\frac{\sinh\alpha}{\sinh 2\tilde{a}}=\frac{\sinh\beta}{\sinh
2\tilde{b}}=\frac{\sinh\gamma}{\sinh 2\tilde{c}}.
\end{eqnarray*}
It follows from (\ref{sinh gamma=}) together with
$\kappa\,z=2\cosh\tilde{c}$ that the above common value is given by %
\begin{eqnarray}\label{eqn:by sine rule}
\frac{\sinh\gamma}{\sinh 2\tilde{c}}
&=&\pm \frac{(\sinh\nu)\,\kappa\, z}{\kappa\, \sinh\tilde{a}\,\sinh\tilde{b}\, \sinh 2\tilde{c}} \nonumber \\ %
&=&\pm \frac{\sinh\nu}{\kappa\, \sinh\tilde{a}\,\sinh\tilde{b}\,\sinh\tilde{c}}, %
\end{eqnarray}
(note that in (\ref{eqn:by sine rule}) the $\pm$ is the same as
that in (\ref{sinh gamma=}) etc) and hence the $\pm$ in
expressions (\ref{sinh gamma=}) etc of $\sinh\gamma, \sinh\alpha,
\sinh\beta$ are constant. So are the signs in (\ref{eqn:alpha+pi
i}), that is, we have either, modulo $2\pi i$,
\begin{eqnarray}\label{eqn:alpha+pi i=+}
\alpha+\pi i=\Psi(y,z,x), \quad  \beta+\pi i=\Psi(z,x,y), \quad \gamma+\pi i=\Psi(x,y,z) %
\end{eqnarray}
or, modulo $2\pi i$,
\begin{eqnarray}\label{eqn:alpha+pi i=-}
\alpha+\pi i=-\Psi(y,z,x), \quad  \beta+\pi i=-\Psi(z,x,y), \quad \gamma+\pi i=-\Psi(x,y,z). %
\end{eqnarray}

Since we also have $(\alpha+\pi i)+(\beta+\pi i)+(\gamma+\pi
i)=\nu \mod 2 \pi i$ (Lemma \ref{lem:alpha+beta+gamma=nu}) and
$\nu \neq 0$ we may conclude that (\ref{eqn:alpha+pi i=+}) must
hold. This proves Lemma \ref{lem:sinh gamma}.
\end{pf}

Note that
\begin{eqnarray*}
\alpha+\pi i&=&\Delta_{{\bf a}(RPQ)}\big(\lfloor {\bf a}(RPQ), %
{\bf a}(Q) \rfloor, \lfloor {\bf a}(R), {\bf a}(RPQ) \rfloor \big)+\pi i \\ %
&=&\Delta_{{\bf a}(RPQ)}\big(\lfloor {\bf a}(RPQ), %
{\bf a}(Q) \rfloor, \lfloor {\bf a}(RPQ), {\bf a}(R) \rfloor \big) %
\end{eqnarray*}
and ${\bf a}(RPQ)={\bf a}(-B^{-1}A^{-1}C^{-1})$. %
Then the above argument gives us the following geometric
interpretation of $\Psi(y,z,x)$ etc.

\vskip 8pt

\begin{lem}\label{lem:geom of Psi} Given $x,y,z \in \mathbf C$ so
that $\mu:=x^2+y^2+z^2-xyz \neq 0,4$, there exist line matrices
$Q,R,P \in \SLtwoC$ {\rm(}that is, $Q^2=R^2=P^2=-I${\rm)} such
that $A:=-RQ,\,B:=-PR,\,C:=-QP \in \SLtwoC$ satisfy ${\rm
tr}\,A=x,\, {\rm tr}\,B=y$ and ${\rm tr}\,C=z$. Let
$A':=[B^{-1},A^{-1}], \, B':=[C^{-1},B^{-1}], \, C':=[A^{-1},C^{-1}]$. %
Then
\begin{eqnarray}\label{eqn:geom of Psi(z)}
\Psi(y,z,x)=\Delta_{{\bf a}(-A')}\big(\lfloor {\bf a}(-A'), {\bf
a}(Q) \rfloor, \lfloor {\bf a}(-A'), {\bf a}(R) \rfloor \big) \mod 2\pi i, %
\end{eqnarray}
\begin{eqnarray}\label{eqn:geom of Psi(x)}
\Psi(z,x,y)=\Delta_{{\bf a}(-B')}\big(\lfloor {\bf a}(-B'), {\bf
a}(R) \rfloor, \lfloor {\bf a}(-B'), {\bf a}(P) \rfloor \big) \mod 2\pi i, %
\end{eqnarray}
\begin{eqnarray}\label{eqn:geom of Psi(y)}
\Psi(x,y,z)=\Delta_{{\bf a}(-C')}\big(\lfloor {\bf a}(-C'), {\bf
a}(P) \rfloor, \lfloor {\bf a}(-C'), {\bf a}(Q) \rfloor \big) \mod 2\pi i, %
\end{eqnarray} \vskip 3pt
\noindent where $\Psi(y,z,x),\Psi(z,x,y),\Psi(x,y,z) \in \mathbf
C$ are defined by {\rm(\ref{Eqn:Psi(x,y,z)=log})}. \square
\end{lem}




\vskip 20pt
\section{{\bf Appendix B - Drawing the gaps}}\label{s:drawing the gaps}
\vskip 10pt

Here we explain how to visualize and draw correctly, for a given
representation $\rho: \Gamma \rightarrow {\rm SL}(2,\mathbf C)$
satisfying the BQ-conditions, the gaps in the extended complex
plane which is the ideal boundary of ${\mathbb H}^3$ in the upper
half-space model.

\vskip 10pt

Fix a pair of generators $a,b \in \Gamma$ and let (in this section
only) $c:=[b^{-1},a^{-1}]=b^{-1}a^{-1}ba$. Then each element of
$\Gamma$ can be written uniquely as a reduced word in the letters
$a,b,a^{-1},b^{-1}$. For each $[g] \in {\hat \Omega}$, we wish to
draw the gaps for an infinite set  of representatives $g_i$ in the
class $[g]$,  with respect to the axis ${\bf a}(C)$ of the image
$C=\rho(c)$ of the commutator $c$. These representatives are
essentially the cyclically reduced words, and conjugates of the
cyclically reduced word by the elements of the commutator
subgroup. As discussed in the previous section, for each $g_i$, we
need to produce a word $g_i' \in \Gamma$ which is a conjugate of
$g_i$ so that $(g_i')^{-1}g_i=c$ in $\Gamma$. Hence we must have
$g_i'=g_ic^{-1}$ after cancellation. We need to do this
consistently so that all gaps are in their correct places viewed
in the extended complex plane.

\vskip 10pt

\begin{figure}
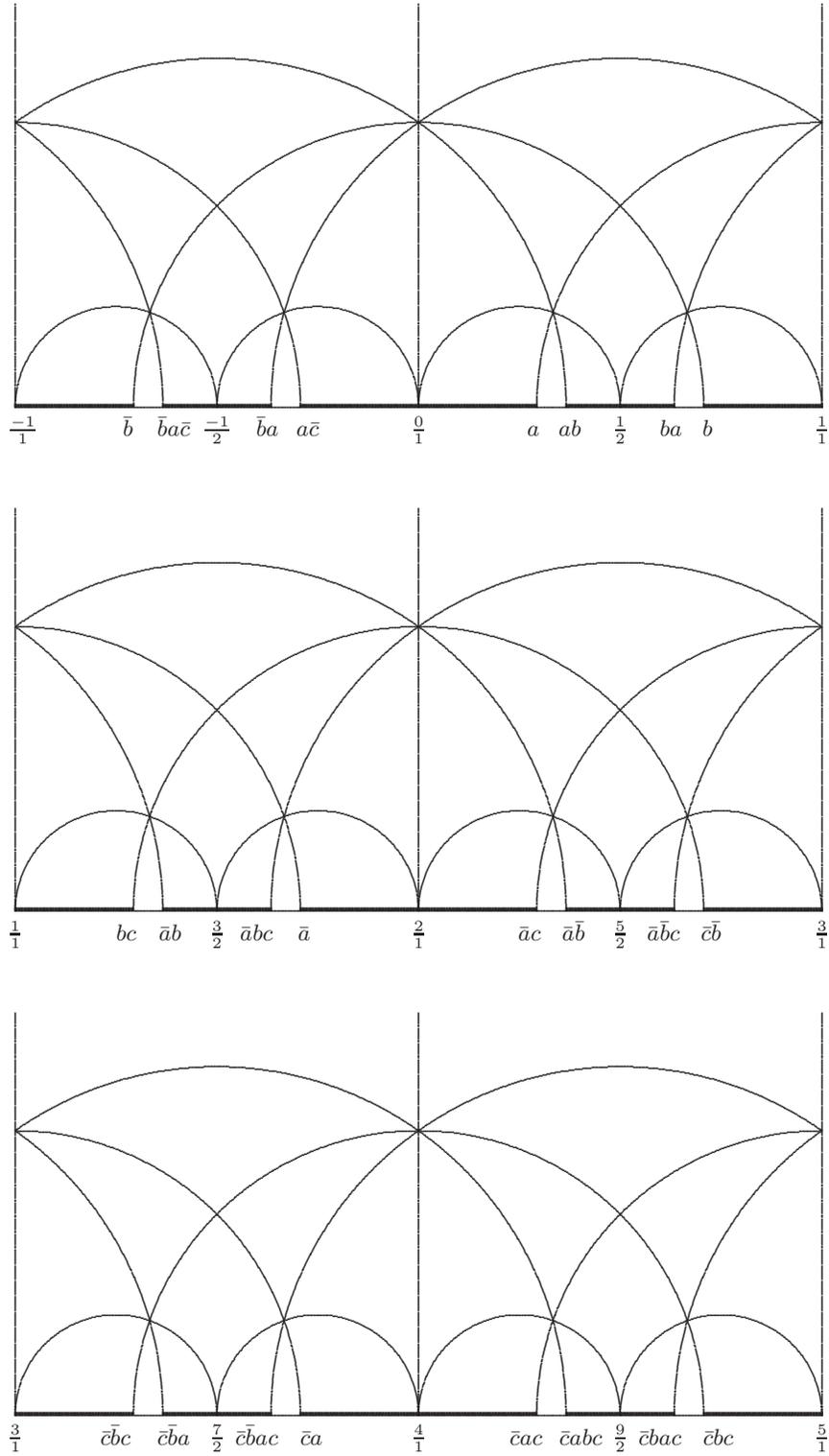


\vskip -25pt
\hskip -25pt \mbox{
\beginpicture

\setcoordinatesystem units <1.10in,1.10in>

\setplotarea x from -2 to 2, y from 0 to 2.5

\plot -2 0 2 0 /


\plot -2 0 -2 2.0 / \plot 0 0 0 2.0 / \plot 2 0 2 2.0 /



\circulararc 180 degrees from -1 0  center at -1.5 0

\circulararc 180 degrees from 0 0  center at -0.5 0

\circulararc 180 degrees from 1 0  center at 0.5 0

\circulararc 180 degrees from 2 0  center at 1.5 0


\circulararc 90 degrees from -0.5858 0  center at -2 0

\circulararc 180 degrees from 1.4142 0  center at 0 0

\circulararc 90 degrees from 2 1.4142  center at 2 0


\circulararc 54.7356 degrees from -1.26795 0 center at -3 0

\circulararc 125.2644 degrees from 0.73205 0 center at -1 0

\circulararc -125.2644 degrees from -0.73205 0 center at 1 0

\circulararc -54.7356 degrees from 1.26795 0 center at 3 0


\put {\mbox{\small $\frac{-1}{1}$}} [cb] <1mm,-5mm> at -2 0

\put {\mbox{\small $\frac{-1}{2}$}} [cb] <0mm,-5mm> at -1 0

\put {\mbox{\small $\frac{0}{1}$}} [cb] <0mm,-5mm> at 0 0

\put {\mbox{\small $\frac{1}{2}$}} [cb] <0mm,-5mm> at 1 0

\put {\mbox{\small $\frac{1}{1}$}} [cb] <0mm,-5mm> at 2 0


\plot -2.0 0.01 -1.4142 0.01 /\plot -2.0 0.015 -1.4142 0.015
/\plot -2.0 0.005 -1.4142 0.005 /

\plot -0.732 0.01 -1.268 0.01 /\plot -0.732 0.015 -1.268 0.015
/\plot -0.732 0.005 -1.268 0.005 /

\plot -0.5858 0.01 0.5858 0.01 /\plot -0.5858 0.015 0.5858 0.015
/\plot -0.5858 0.005 0.5858 0.005 /

\plot 1.4142 0.01 2.0 0.01 /\plot 1.4142 0.015 2.0 0.015 /\plot
1.4142 0.005 2.0 0.005 /

\plot 0.732 0.01 1.268 0.01 /\plot 0.732 0.015 1.268 0.015 /\plot
0.732 0.005 1.268 0.005 /


\put {\mbox{\small ${\bar b}a{\bar c}$}} [cb] <1.5mm,-4mm> at
-1.268 0

\put {\mbox{\small ${\bar b}a$}} [cb] <-0.5mm,-4mm> at -0.732 0

\put {\mbox{\small $a{\bar c}$}} [cb] <1mm,-4mm> at -0.5858 0

\put {\mbox{\small $a$}} [cb] <-0.5mm,-4mm> at 0.5858 0

\put {\mbox{\small $ab$}} [cb] <0.5mm,-4mm> at 0.732 0

\put {\mbox{\small $ba$}} [cb] <-0.5mm,-4mm> at 1.268 0

\put {\mbox{\small $b$}} [cb] <0.5mm,-4mm> at 1.4142 0

\put {\mbox{\small ${\bar b}$}} [cb] <-0.8mm,-4mm> at -1.4142 0

\endpicture
}

\vskip -15pt

\hskip -25pt \mbox{
\beginpicture

\setcoordinatesystem units <1.10in,1.10in>

\setplotarea x from 2 to 6, y from 0 to 2.5

\plot 2 0 6 0 /


\plot 2 0 2 2.0 / \plot 4 0 4 2.0 / \plot 6 0 6 2.0 /


\circulararc 180 degrees from 3 0  center at 2.5 0

\circulararc 180 degrees from 4 0  center at 3.5 0

\circulararc 180 degrees from 5 0  center at 4.5 0

\circulararc 180 degrees from 6 0  center at 5.5 0


\circulararc 90 degrees from 3.4142 0  center at 2 0

\circulararc 180 degrees from 5.4142 0  center at 4 0

\circulararc 90 degrees from 6 1.4142  center at 6 0


\circulararc 54.7356 degrees from 2.732 0 center at 1 0

\circulararc 125.2644 degrees from 4.732 0 center at 3 0

\circulararc -125.2644 degrees from 3.268 0 center at 5 0

\circulararc -54.7356 degrees from 5.268 0 center at 7 0


\put {\mbox{\small $\frac{1}{1}$}} [cb] <0mm,-5mm> at 2 0

\put {\mbox{\small $\frac{3}{2}$}} [cb] <0mm,-5mm> at 3 0

\put {\mbox{\small $\frac{2}{1}$}} [cb] <0mm,-5mm> at 4 0

\put {\mbox{\small $\frac{5}{2}$}} [cb] <0mm,-5mm> at 5 0

\put {\mbox{\small $\frac{3}{1}$}} [cb] <0mm,-5mm> at 6 0


\plot 2.0 0.01 2.5858 0.01 /\plot 2.0 0.015 2.5858 0.015 /\plot
2.0 0.005 2.5858 0.005 /

\plot 2.732 0.01 3.268 0.01 /\plot 2.732 0.015 3.268 0.015 /\plot
2.732 0.005 3.268 0.005 /

\plot 3.4142 0.01 4.5858 0.01 /\plot 3.4142 0.015 4.5858 0.015
/\plot 3.4142 0.005 4.5858 0.005 /

\plot 4.732 0.01 5.268 0.01 /\plot 4.732 0.015 5.268 0.015 /\plot
4.732 0.005 5.268 0.005 /

\plot 5.4142 0.01 6.0 0.01 /\plot 5.4142 0.015 6.0 0.015 /\plot
5.4142 0.005 6.0 0.005 /


\put {\mbox{\small $bc$}} [cb] <-1.0mm,-4mm> at 2.5858 0

\put {\mbox{\small ${\bar a}b$}} [cb] <1.0mm,-4mm> at 2.732 0

\put {\mbox{\small ${\bar a}bc$}} [cb] <-2.0mm,-4mm> at 3.268 0

\put {\mbox{\small ${\bar a}$}} [cb] <0.5mm,-4mm> at 3.4142 0

\put {\mbox{\small ${\bar a}c$}} [cb] <-1.0mm,-4mm> at 4.5858 0

\put {\mbox{\small ${\bar a}{\bar b}$}} [cb] <1.0mm,-4mm> at 4.732
0

\put {\mbox{\small ${\bar a}{\bar b}c$}} [cb] <-1.5mm,-4mm> at
5.268 0

\put {\mbox{\small ${\bar c}{\bar b}$}} [cb] <1mm,-4mm> at 5.4142
0

\endpicture
}


\vskip -15pt


\hskip -25pt \mbox{
\beginpicture

\setcoordinatesystem units <1.10in,1.10in>

\setplotarea x from 6 to 10, y from 0 to 2.5

\plot 6 0 10 0 /


\plot 6 0 6 2.0 / \plot 8 0 8 2.0 / \plot 10 0 10 2.0 /


\circulararc 180 degrees from 7 0  center at 6.5 0

\circulararc 180 degrees from 8 0  center at 7.5 0

\circulararc 180 degrees from 9 0  center at 8.5 0

\circulararc 180 degrees from 10 0  center at 9.5 0


\circulararc 90 degrees from 7.4142 0  center at 6 0

\circulararc 180 degrees from 9.4142 0  center at 8 0

\circulararc 90 degrees from 10 1.4142  center at 10 0


\circulararc 54.7356 degrees from 6.732 0 center at 5 0

\circulararc 125.2644 degrees from 8.732 0 center at 7 0

\circulararc -125.2644 degrees from 7.268 0 center at 9 0

\circulararc -54.7356 degrees from 9.268 0 center at 11 0


\put {\mbox{\small $\frac{3}{1}$}} [cb] <0mm,-5mm> at 6 0

\put {\mbox{\small $\frac{7}{2}$}} [cb] <0mm,-5mm> at 7 0

\put {\mbox{\small $\frac{4}{1}$}} [cb] <0mm,-5mm> at 8 0

\put {\mbox{\small $\frac{9}{2}$}} [cb] <0mm,-5mm> at 9 0

\put {\mbox{\small $\frac{5}{1}$}} [cb] <0mm,-5mm> at 10 0


\plot 6.0 0.01 6.5858 0.01 /\plot 6.0 0.015 6.5858 0.015 /\plot
6.0 0.005 6.5858 0.005 /

\plot 6.732 0.01 7.268 0.01 /\plot 6.732 0.015 7.268 0.015 /\plot
6.732 0.005 7.268 0.005 /

\plot 7.4142 0.01 8.5858 0.01 /\plot 7.4142 0.015 8.5858 0.015
/\plot 7.4142 0.005 8.5858 0.005 /

\plot 8.732 0.01 9.268 0.01 /\plot 8.732 0.015 9.268 0.015 /\plot
8.732 0.005 9.268 0.005 /

\plot 9.4142 0.01 10.0 0.01 /\plot 9.4142 0.015 10.0 0.015 /\plot
9.4142 0.005 10.0 0.005 /


\put {\mbox{\small ${\bar c}{\bar b}c$}} [cb] <-2.5mm,-4mm> at
6.5858 0

\put {\mbox{\small ${\bar c}{\bar b}a$}} [cb] <1.5mm,-4mm> at
6.732 0

\put {\mbox{\small ${\bar c}{\bar b}ac$}} [cb] <-2.0mm,-4mm> at
7.268 0

\put {\mbox{\small ${\bar c}a$}} [cb] <1.5mm,-4mm> at 7.4142 0

\put {\mbox{\small ${\bar c}ac$}} [cb] <-1.5mm,-4mm> at 8.5858 0

\put {\mbox{\small ${\bar c}abc$}} [cb] <2.0mm,-4mm> at 8.732 0

\put {\mbox{\small ${\bar c}bac$}} [cb] <-2.0mm,-4mm> at 9.268 0

\put {\mbox{\small ${\bar c}bc$}} [cb] <2mm,-4mm> at 9.4142 0

\endpicture
}

\vskip 10pt

\caption{The sequence of pairs of words and corresponding
gaps}\label{g:word sequence}
\end{figure}

\noindent {\bf Constructing the `universal cover' of $\hat
\Omega$.}\,\, From the above discussion, to draw the gaps with
respect to a fixed commutator, we really need to construct the
`universal cover' of the rational projective lamination space
$\hat \Omega$ as follows. We shall construct pairs of words
$(L_\frac{p}{q}, R_\frac{p}{q})$ in $\Gamma$, parametrized by
$\frac{p}{q} \in \mathbf Q$, where $\mathbf Q$ is regarded as the
``universal cover'' of ${\mathbf Q}\cup \infty$.  The particular
representation we will draw will be that corresponding to the
once-punctured hyperbolic torus which has lifted holonomy
representation $\rho: \Gamma \rightarrow \SLtwoC$ such that
\begin{eqnarray}\label{eqn:A=,B=}
A=\rho(a)=\bigg(\!\!\begin{array}{cc}0 & \sqrt{2}\big/2 \\
-\sqrt{2} & 2\sqrt{2}\end{array}\!\bigg) \quad {\rm and} \quad
B=\rho(b)=\bigg(\!\!\begin{array}{cc}\sqrt{2} & \sqrt{2}\big/2 \\
\sqrt{2} & \sqrt{2}\end{array}\!\bigg).
\end{eqnarray} %
Note that
$C:=\rho(c)=B^{-1}A^{-1}BA=\bigg(\begin{matrix}-1 & 4 \\
0 & -1\end{matrix}\bigg).$ See Figure \ref{g:word sequence}, where
we use ${\bar a}, {\bar b}, {\bar c}$ to denote $a^{-1}, b^{-1},
c^{-1}$ respectively and we label the attracting fixed points of
$A, A^{-1}$ by $a, {\bar a}$ respectively, and similarly for the
other words.

\vskip 5pt

First, we construct a bi-infinite sequence of pairs of reduced
words $(L_n, R_n)$, $n = \frac{n}{1} \in {\mathbf Z}$, inductively
as follows.

\vskip 4pt
\begin{itemize}

\item[(i)] Set $R_0=a$ and $L_1=b$.

\item[(ii)] We require
\begin{eqnarray*}
(L_n)^{-1}R_n=c
\end{eqnarray*}
for all $n \in {\mathbf Z}$. Hence, for example,
$L_0=b^{-1}ab=ac^{-1}$ and $R_1=a^{-1}ba=bc$.

\item[(iii)] Define
\begin{eqnarray*}
L_{n+2}=(R_{n})^{-1}
\end{eqnarray*}
or equivalently
\begin{eqnarray*}
R_{n-2}=(L_n)^{-1}
\end{eqnarray*}
for $n \in {\mathbf Z}$. Hence $R_{n+2}=(R_{n})^{-1}c$.
\end{itemize}
\vskip 4pt

\noindent It is easy to check that the following relations hold
for all $n \in {\mathbf Z}$:
\begin{eqnarray*}
[(L_{n+1})^{-1},(R_{n})^{-1}]=c,
\end{eqnarray*}
and
\begin{eqnarray*}
L_{n+4}=c^{-1}\,L_n\,c.
\end{eqnarray*}

\vskip 10pt

Next, we generate the general $(L_\frac{p}{q}, R_\frac{p}{q})$ for
all $\frac{p}{q} \in {\mathbf Q}$ by the following rule:

\vskip 5pt
\begin{itemize}
\item[(iv)] For $\frac{p+r}{q+s} \in {\mathbf Q}$, constructed
from Farey neighbors $\frac{p}{q}, \frac{r}{s} \in {\mathbf Q}$,
where $\frac{p}{q} < \frac{r}{s}$, we define
\begin{eqnarray}
L_\frac{p+r}{q+s}=R_\frac{p}{q}L_\frac{r}{s}, \quad
R_\frac{p+r}{q+s}=L_\frac{r}{s}R_\frac{p}{q}.
\end{eqnarray}
\end{itemize}

\noindent {\it Remark.}\,\, We  note that here $\frac{p}{q}$ in
the index is {\it not}  the slope of the free homotopy class
$[L_\frac{p}{q}]=[R_\frac{p}{q}] \in \hat \Omega$. In fact the
rationals in the interval $[0,2)$ in our index set corresponds to
the set ${\mathbf Q} \cup \infty$ of the slopes, via the Farey
triangulation described, and the full set of rationals ${\mathbf
Q}$ in our index set corresponds to the universal cover. The
reason we need to do this is because we need to consider the inner
automorphisms by the commutator subgroup, for a careful treatment,
see \cite{goldman2003gt}.

\vskip 10pt

We can prove by induction that
\begin{itemize}

\item[(I)] for any Farey neighbors $\frac{p}{q}, \frac{r}{s} \in
{\mathbf Q}$, where $\frac{p}{q} < \frac{r}{s}$, there is the
commutator identity:
\begin{eqnarray}
[(L_\frac{r}{s})^{-1},(R_\frac{p}{q})^{-1}]=c;
\end{eqnarray}

\item[(II)] for all $\frac{p}{q} \in {\mathbf Q}$, there is the
conjugation relation:
\begin{eqnarray}
L_{\frac{p}{q}+2}=(R_{\frac{p}{q}})^{-1}.
\end{eqnarray}
and
\begin{eqnarray}
L_{\frac{p}{q}+4}=c^{-1}\,L_{\frac{p}{q}}\,c.
\end{eqnarray}
\end{itemize}

\vskip 10pt

Thus the pair of words $(L_{\frac{p}{q}}, R_{\frac{p}{q}})$ we
constructed above are indeed conjugate to each other as the
following lemma shows.

\vskip 10pt

\begin{lem}
For any Farey neighbors $\frac{p}{q}, \frac{r}{s} \in {\mathbf
Q}$, where $\frac{p}{q} < \frac{r}{s}$, we have
\begin{eqnarray}
L_{\frac{p}{q}}=(L_{\frac{r}{s}})^{-1}\,R_{\frac{p}{q}}\,L_{\frac{r}{s}}.
\end{eqnarray} 
\end{lem}

\vskip 5pt

\begin{figure}
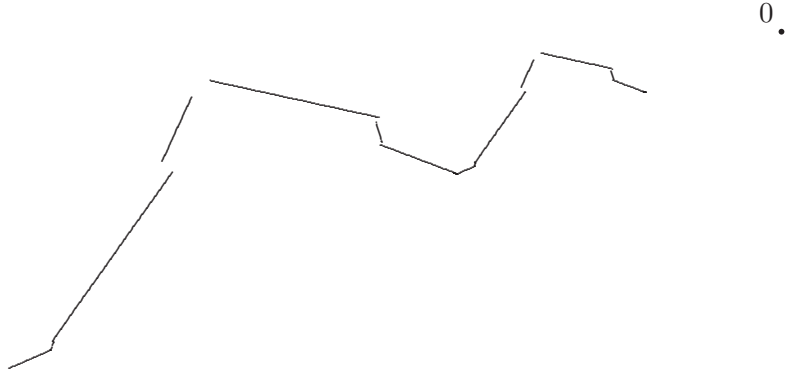

\begin{center}
\mbox{
\beginpicture
\setcoordinatesystem units <4.0in,4.0in>

\setplotarea x from -1.1 to 0.4, y from -0.5 to 0.1

\plot  -1.010780664 -0.4388992616 -0.9562910360 -0.4152388786 /
\plot  -0.9562878185 -0.4152295084 -0.9561816310 -0.4149207342 /
\plot  -0.9561824865 -0.4149191756 -0.9525007090 -0.4045700153 /
\plot  -0.9525109190 -0.4045669982 -0.9533564425 -0.4028036382 /
\plot  -0.9537428650 -0.4027252413 -0.7971707590 -0.1820163898 /
\plot  -0.8104781465 -0.1677672250 -0.7718412725 -0.0834555917 /
\plot  -0.7473032095 -0.0621087595 -0.5267954040 -0.1102367685 /
\plot  -0.5301290745 -0.1170619722 -0.5228116280 -0.1432953496 /
\plot  -0.5249237695 -0.1459811350 -0.4278646656 -0.1826768884 /
\plot  -0.4278652317 -0.1826765194 -0.4278702598 -0.1826781480 /
\plot  -0.4278702596 -0.1826781476 -0.4278901062 -0.1827183858 /
\plot  -0.4278901064 -0.1827183858 -0.4278901224 -0.1827183766 /
\plot  -0.4278901224 -0.1827183764 -0.4277622286 -0.1830772583 /
\plot  -0.4277622286 -0.1830772584 -0.4277622288 -0.1830772580 /
\plot  -0.4277622290 -0.1830772582 -0.4277628620 -0.1830784528 /
\plot  -0.4277628620 -0.1830784528 -0.4277628666 -0.1830784539 /
\plot  -0.4277628666 -0.1830784541 -0.4248235065 -0.1843870472 /
\plot  -0.4248235064 -0.1843870472 -0.4248235072 -0.1843870482 /
\plot  -0.4248235067 -0.1843870482 -0.4248234010 -0.1843873678 /
\plot  -0.4248234008 -0.1843873677 -0.4248234010 -0.1843873678 /
\plot  -0.4248234010 -0.1843873678 -0.4247368061 -0.1844260940 /
\plot  -0.4247368062 -0.1844260938 -0.4247368062 -0.1844260944 /
\plot  -0.4247368062 -0.1844260942 -0.4247342568 -0.1844272344 /
\plot  -0.4247342568 -0.1844272344 -0.4247341821 -0.1844272682 /
\plot  -0.4247341795 -0.1844272690 -0.4018374148 -0.1744850794 /
\plot  -0.4018360626 -0.1744811420 -0.4017914421 -0.1743513936 /
\plot  -0.4017918016 -0.1743507388 -0.4002447038 -0.1700019792 /
\plot  -0.4002489939 -0.1700007112 -0.4006042866 -0.1692597401 /
\plot  -0.4007666628 -0.1692267979 -0.3349744218 -0.07648403290 /
\plot -0.3405662406 -0.07049647610 -0.3243308678 -0.03506838171 /
\plot -0.3140198740 -0.02609835519 -0.2213615888 -0.04632194208 /
\plot -0.2227624112 -0.04918992051 -0.2196875901 -0.06021329350 /
\plot -0.2205751206 -0.06134187130 -0.1797904872 -0.07676157745 /
\plot -0.1797907251 -0.07676142240 -0.1797928379 -0.07676210670 /
\plot -0.1797928378 -0.07676210660 -0.1798011775 -0.07677901475 /
\plot -0.1798011774 -0.07677901475 -0.1798011842 -0.07677901090 /
\plot -0.1798011842 -0.07677901085 -0.1797474428 -0.07692981455 /
\plot -0.1797474428 -0.07692981455 -0.1797474429 -0.07692981440 /
\plot -0.1797474430 -0.07692981445 -0.1797477089 -0.07693031640 /
\plot -0.1797477089 -0.07693031645 -0.1797477109 -0.07693031695 /
\plot -0.1797477108 -0.07693031695 -0.1785125798 -0.07748019315 /
\plot -0.1785125797 -0.07748019320 -0.1785125800 -0.07748019355 /
\plot -0.1785125798 -0.07748019360 -0.1785125353 -0.07748032785 /
\plot -0.1785125354 -0.07748032790 -0.1785125354 -0.07748032790 /
\plot -0.1785125353 -0.07748032790 -0.1784761479 -0.07749660080 /
\plot -0.1784761479 -0.07749660080 -0.1784761478 -0.07749660090 /
\plot -0.1784761479 -0.07749660090 -0.1784750768 -0.07749708005 /
\plot -0.1784750767 -0.07749708000 -0.1784750451 -0.07749709415 /
\put {\mbox{\large $0$}} [ct] <0mm,1mm> at -0.02 0.03 \put
{\mbox{\Huge $\cdot$}} [cc] <0mm,0mm> at 0 0
\endpicture
}\end{center}

\caption{Gaps for a $(-0.8)$-Markoff map with $x=6,y=3+i$}
\label{fig:drawing_gaps_1}
\end{figure}

\noindent {\bf Drawing the gaps.}\,\, Now for each pair of words
$(L_\frac{p}{q}, R_\frac{p}{q})$, $\frac{p}{q} \in {\mathbf Q}$,
constructed as above, we draw a line segment in the extended
complex plane connecting the attracting fixed points of
$\rho(L_\frac{p}{q})$ and $\rho(R_\frac{p}{q})$ to indicate the
gap corresponding to the pair. Of course these gaps should be
measured against the commutator, in the case where the commutator
is parabolic, we may normalize the picture so that the commutator
fixes $\infty$,  see Figure \ref{g:word sequence}. When the
commutator is loxodromic we may normalize so that the fixed points
of the commutator $C$ are 0 and $\infty$. The picture obtained for
a particular generalized Markoff map $\phi \in ({\bf
\Phi}_\mu)_{Q}$ is given in Figure \ref{fig:drawing_gaps_1} (with
this normalization), and we have drawn some of the gaps for
$\frac{p}{q}$ between $0$ and $4$. The complete picture for
$\frac{p}{q} \in \mathbf Q$ can be obtained by applying the
M\"{o}bius transformation $C=\rho(c)$ and its inverse on the given
picture repeatedly.

\vskip 10pt

\noindent {\bf The attracting fixed points.}\,\, Finally we
explain the formula to find the attracting fixed point for a given
loxodromic (hyperbolic) $A \in {\rm SL}(2,\mathbf C)$.

Let $A=\begin{pmatrix} A_{11} & A_{12} \\ A_{21} & A_{22}
\end{pmatrix}$. The two fixed points of $A$ are the roots of
quadratic equation $A_{21}z^2+(A_{22}-A_{11})z-A_{12}=0$, hence
given by
\begin{eqnarray*}
z=[(A_{11}-A_{22})\pm\sqrt{(A_{11}+A_{22})^2-4}\,]\big/2.
\end{eqnarray*}
Recall that the square root here has been assumed to have positive
real part.

In this form, however, it is {\it not} true that one sign gives
always the attracting or always the repelling fixed point. We can
rewrite it in the following form:
\begin{eqnarray*}
z=\big[ (A_{11}-A_{22}) \pm (A_{11}+A_{22})
\sqrt{1-4(A_{11}+A_{22})^{-2}}\, \big] \big/2.
\end{eqnarray*}
Then corresponding to the plus and minus signs are respectively
the attracting and repelling fixed points.

\begin{lem}\label{lem:Fix} Suppose $A \in \SLtwoC$ is loxodromic/hyperbolic. Then
\begin{eqnarray}
& &\hspace{-20pt}{\rm Fix}^{+}(A)=\big[ (A_{11}-A_{22}) +
(A_{11}+A_{22})\sqrt{1-4(A_{11}+A_{22})^{-2}}\, \big]\big/2, \\
& &\hspace{-20pt}{\rm Fix}^{-}(A)=\big[ (A_{11}-A_{22}) -
(A_{11}+A_{22})\sqrt{1-4(A_{11}+A_{22})^{-2}}\, \big]\big/2. %
\end{eqnarray}
\end{lem}

\begin{pf}
This is true since we know that they correspond respectively the
two eigenvalues of the matrix $A$:
\begin{eqnarray}
& &{\lambda}^{+}=(A_{11}+A_{22})\,[1 + \sqrt{1-4(A_{11}+A_{22})^{-2}}\,]\big/2, \\
& &{\lambda}^{-}=(A_{11}+A_{22})\,[1 - \sqrt{1-4(A_{11}+A_{22})^{-2}}\,]\big/2  %
\end{eqnarray}
which have respectively norms greater and less than $1$.
\end{pf}


\end{document}